\documentclass[a4paper]{article}
\usepackage{amsmath,amssymb,amsthm,amscd,amsfonts}
\allowdisplaybreaks 

\def\bc{\begin{center}}
\def\ec{\end{center}} \DeclareMathOperator{\Ad}{Ad} \DeclareMathOperator{\ad}{ad}\DeclareMathOperator{\degr}{deg}
\DeclareMathOperator{\Diff}{Diff}
 \DeclareMathOperator{\Tr}{Tr}

\DeclareMathOperator{\id}{id} \DeclareMathOperator{\im}{Im}
 \DeclareMathOperator{\Ann}{Ann}
\DeclareMathOperator{\SO}{SO}\DeclareMathOperator{\Spin}{Spin}
\DeclareMathOperator{\Simp}{Sp} \DeclareMathOperator{\U}{U}
\DeclareMathOperator{\uu}{u} 
 \DeclareMathOperator{\SU}{SU}

\DeclareMathOperator{\RE}{Re}\newcommand{\D}{\bar D}

\newcommand{\x}{{\bf x}}
\newcommand{\y}{{\bf y}}
\newcommand{\z}{{\bf z}}
\newcommand{\ii}{{\bf i}}
\newcommand{\jj}{{\bf j}}
\newcommand{\kk}{{\bf k}} \newtheorem{theore}{Theorem}

\newtheorem{proposit}{Proposition}
\newtheorem{Rem}{Remark}
\newtheorem{Lem}{Lemma}
 \def\pd#1#2{\frac{\partial{#1}}{\partial{#2}}}
\def\pd1#1{\frac{\partial}{\partial#1}}
\def\d1#1{\frac{d}{d#1}}

\begin{document} \author{A.V. Shchepetilov\\ Department of Physics, Moscow State University,
\\ 119992 Moscow, Russia\footnote{email:
alexey@quant.phys.msu.su}}

\title{Algebras of invariant differential operators on unit sphere bundles over two-point homogeneous Riemannian spaces}
\date{}
\maketitle

\begin{abstract} Let $G$ be the identity component of the isometry group for an arbitrary curved two-point homogeneous space
$M$. We consider algebras of $G$-invariant differential operators
on bundles of unit spheres over $M$. The generators of this
algebra and the corresponding relations for them are found. The
connection of these generators with two-body problem on two-point
homogeneous spaces is discussed.
\end{abstract}

Keywords: two point-homogeneous spaces, invariant differential
operators, two-body problem

PACS number: 02.40.Vh, 02.40.Ky, 02.20.Qs, 03.65.Fd

Mathematical Subject Classification: 16S32, 43A85, 22F30, 70F05.
\newpage

\section{Introduction} \label{Introduction}\markright{\ref{Introduction}
Introduction}

The property of a differential operator on a smooth manifold $M$
to be invariant with respect to an action of some group $G$
(especially a Lie group) on $M$ plays a great role in mathematical
physics since it helps select physically significant operators.
The algebra $\Diff(M)$ of all $G$-invariant differential operators
with complex or real coefficients on $M$ gives the material for
constructing $G$-invariant physical theories on $M$. Properties of
such theory are in close connection with properties of the algebra
$\Diff(M)$.

A homogeneous smooth manifold $M$ of the Lie group $G$ is called
{\it commutative space}, if the algebra $\Diff(M)$ is commutative.
The well known example of a commutative space is the symmetric
space of the rank $l$. Recall that the {\it rank of a symmetric
space} is the dimension of its maximal flat completely geodesic
submanifold. The commutative algebra $\Diff(M)$ for this space is
generated by $l$ independent commutative generators
\cite{Hel1973}. Particularly, for symmetric spaces of the rank one
(which are the same as two-point homogeneous spaces) the algebra
$\Diff(M)$ is generated by the Laplace-Beltrami operator. Also,
the class of commutative spaces contains weakly symmetric spaces
\cite{Vinberg}.

There are known only some sporadic examples of noncommutative
algebras $\Diff(M)$ (see, for example Ch.2, \cite{Hel84}). One of
these example is the noncommutative algebra $\Diff(M_{1})$ for
$M_{1}=\mathbf{O}_{0}(1,n)/\mathbf{SO}(n-1)$ studied in
\cite{Reimann}, where $\mathbf{O}_{0}(1,n)$ is the identity
component for the group $\mathbf{O}(1,n)$. In that paper the space
$M_{1}$ was interpreted as the total space for the bundles of unit
spheres over the hyperbolic space $\mathbf{H}^{n}(\mathbb{R})$.
Denote the total space of the bundle of unit spheres over a
Riemannian space $M$ by $M_{\mathbb{S}}$.

The space $\mathbf{H}^{n}(\mathbb{R})$ is a representative of the
class of two-point homogeneous Riemannian spaces (TPHRS) for which
any pair of points can be transformed by means of appropriate
isometry to any other pair of points with the same distance
between them. Equivalently, these spaces are characterized by the
property that the natural action of the isometry group on the
bundle of unit spheres over them are transitive. Thus the natural
problem arises: "{\it describe the algebras
$\Diff(M_{\mathbb{S}})$ for the bundle of unit spheres over all
TPHRS $M$}".

From the point of view of the two-body problem in classical and
quantum mechanics, TPHRS are characterizes by the property that
the distance between particles is the only invariant of the
isometry group $G$ in the configuration space. The space
$M_{\mathbb{S}}$ is isomorphic to an orbit in general position for
the symmetry group $G$ of the two-body problem on the TPHRS $M$
acting in the configuration space of this problem. Due to the
two-point homogeneity of $M$ the codimension of these orbits is
one. Thus, for the two-body problem in TPHRS there is the degree
of freedom corresponding to the distance between particles; other
degrees of freedom correspond to the homogeneous manifold and can
be described in terms of the symmetry group $G$.

In the paper \cite{Shchep2002} the polynomial expression for the
Hamiltonian $\hat H$ of the quantum mechanical two-body problem on
an arbitrary TPHRS  $M$ was found through the radial differential
operator and elements of $\Diff(M_{\mathbb{S}})$. In the present
paper the generators of algebras $\Diff(M_{\mathbb{S}})$ and the
corresponding relations for them are found for all curved
two-point homogeneous spaces $M$. Some properties of these
generators are discussed.

This paper is organized as follows. In section \ref{InvOperators}
we perform the necessary information on invariant differential
operators on homogeneous spaces. We recall the classification of
TPHRS in section \ref{twopoint}. In section \ref{spes} we specify
the construction of invariant differential operators for the space
$M_{\mathbb{S}}$, where $M$ is a two-point homogeneous Riemannian
space. There are found some generators common for all
$\Diff(M_{\mathbb{S}})$. In section \ref{QuatMod} the model of the
quaternion projective space $\mathbf{P}^{n}(\mathbb{H})$ is
described. In section \ref{Quat} the generators for the algebra
$\Diff(\mathbf{P}^{n}(\mathbb{H})_{\mathbb{S}})$ are calculated
and by the formal correspondence the analogous generators of the
algebra $\Diff(\mathbf{H}^{n}(\mathbb{H})_{\mathbb{S}})$ is
obtained. The corresponding relations for these algebras are found
in section \ref{Relations}. In section \ref{Complex} we consider
from the same point of view algebras
$\Diff(\mathbf{P}^{n}(\mathbb{C})_{\mathbb{S}})$,
$\Diff(\mathbf{H}^{n}(\mathbb{C})_{\mathbb{S}})$ and in section
\ref{Real} the algebras
$\Diff(\mathbf{P}^{n}(\mathbb{R})_{\mathbb{S}}),\Diff(\mathbf{S}^{n}_{\mathbb{S}})$.
In the section \ref{octo} there is a description of the Cayley
plane $\mathbf{P}^{2}(\mathbb{C}a)$ through the exceptional Jordan
algebra $ \mathfrak{h}_{3}(\mathbb{C}a)$. In section \ref{CaS} the
generators for the algebra
$\Diff(\mathbf{P}^{2}(\mathbb{C}a)_{\mathbb{S}})$ are calculated
and by the formal correspondence the analogous generators of the
algebra $\Diff(\mathbf{H}^{2}(\mathbb{C}a)_{\mathbb{S}})$ is
obtained. The corresponding relations for these algebras are found
in section \ref{RelCa}.

The connection of constructed generators with two-body problem on
two-point homogeneous spaces is discussed in section
\ref{Connect}.

In appendix A we describe the technique for calculating the
commutative relations for algebras of differential operators under
consideration and in appendix B one interesting fact for an
arbitrary TPHRS is proved.

\section{Invariant differential operators on
homogeneous spaces} \label{InvOperators}
\markright{\ref{InvOperators} Invariant differential operators}
Let $G$ be a Lie group, $M$ be a Riemannian $G$-homogeneous left
space, $x_{0}\in M, K\subset G$ be the stationary subgroup of the
point $x_{0}\in M$, and $\mathfrak{k}\subset \mathfrak{g}\equiv
T_{e}G$ be the corresponding Lie algebras. Choose a subspace $
\mathfrak{p}\subset\mathfrak{g}$ such that
$\mathfrak{g}=\mathfrak{p}\oplus\mathfrak{k}$ (a direct sum of
linear spaces). The space $ \mathfrak{p}$ can be identified with
the tangent space $T_{x_{0}}M$.

The stationary subgroup $K$ is compact, since it is also the
subgroup of the group $\mathbf{SO}(n)$. By the averaging on the
group $K$ we can define a $\Ad_{K}$-invariant scalar product
$\langle\cdot,\cdot\rangle$ on $\mathfrak{g}$ and choose the
subspace $\mathfrak{p}$ orthogonal to $\mathfrak{k}$ with respect
to this product \cite{Hel84}, \cite{KoNo2}. In this case we have
the inclusion $\Ad_{K}(\mathfrak{p})\subset \mathfrak{p}$, i.e.
the space $M$ is reductive.

Identify the space $M$ with the factor space of left conjugate
classes of the group $G$ with respect to the subgroup $K$. Let
$\pi:G\to G/K$ be the natural projection.

Let $S(V)$ be a graded symmetric algebra over a finite dimensional
complex space $V$, i.e. a free commutative algebra over the field
$\mathbb{C}$, generated by elements of any basis of $V$. The
adjoint action of the group $G$ on $\mathfrak{g}$ can be naturally
extended to the action of $G$ on the algebra $S(\mathfrak{g})$
according to the formula:
$$
\Ad_{q}:\,Y_{1}\cdot\ldots\cdot Y_{i}\rightarrow
\Ad_{q}(Y_{1})\cdot\ldots\cdot\Ad_{q}(Y_{i}),\,Y_{1},\dots,Y_{i}\in
\mathfrak{g}.
$$
Denote by $\mathfrak{g}^{K}$ the set of all $\Ad$-invariants in
$S(\mathfrak{g})$.

The main result of the invariant differential operators theory is
the existence of the one to one correspondence between the algebra
$\Diff(M)$ and the set $\mathfrak{p}^{K}$ of all
$\Ad_{K}$-invariants in $S(\mathfrak{p})$ \cite{Hel84}. For our
purpose the next version \cite{Vinberg} of this result is more
convenient. Let $U(\mathfrak{g})$ be the universal enveloping
algebra with the standard filtration for the Lie algebra
$\mathfrak{g}$ and $U(\mathfrak{g})^{K}$ be it's subalgebra,
consisting of all $\Ad_{K}$-invariant elements in
$U(\mathfrak{g})$. Let $\mu$ be the linear mapping of
$S(\mathfrak{p})$ into $U(\mathfrak{g})$, according to the formula
$$
\mu(Y_{1}\cdot\ldots\cdot Y_{p})=\frac{1}{p!}\sum\limits_{
\sigma\in \mathfrak{S} _{p}}Y_{\sigma(1)}\cdot\ldots\cdot
Y_{\sigma(p)},
$$
where on the left side the element $Y_{1}\cdot\ldots\cdot Y_{p}$
is supposed to be in $S(\mathfrak{p})$ and on the right side it is
supposed to be in $U(\mathfrak{g})$. Here $\mathfrak{S}_{p}$ is
the permutation group of $p$ elements. Obviously, $\mu$ is
injective.

Let $U(\mathfrak{g})\mathfrak{k}$ be the left ideal in
$U(\mathfrak{g})$ generated by $\mathfrak{k}$ and
$(U(\mathfrak{g})\mathfrak{k})^{K}$ be the set of all
$\Ad_{K}$-invariant elements in $U(\mathfrak{g})\mathfrak{k}$. The
set $(U(\mathfrak{g})\mathfrak{k})^{K}$ is a two-sided ideal in
$U(\mathfrak{g})^{K}$  since for elements $f\in \mathfrak{k}$ and
$g\in U(\mathfrak{g})^{K}$ we have $fg=\ad_{f}g+gf=gf$. Also
$\mu(\mathfrak{p}^{K})\subset U(\mathfrak{g})^{K}$, because $M$ is
reductive. Hence we can define the factor algebra
$U(\mathfrak{g})^{K}/(U(\mathfrak{g})\mathfrak{k})^{K}$. Let
$\eta:\, U(\mathfrak{g})^{K}\to
U(\mathfrak{g})^{K}/(U(\mathfrak{g})\mathfrak{k})^{K}$ be the
canonical projection.

\begin{theore}[\cite{Vinberg}]
The algebras $\Diff(M)$ and
$U(\mathfrak{g})^{K}/(U(\mathfrak{g})\mathfrak{k})^{K}$ are
isomorphic.
\end{theore}
Every element of
$U(\mathfrak{g})^{K}/(U(\mathfrak{g})\mathfrak{k})^{K}$ has the
unique representative from $\mathfrak{p}^{K}$ or equivalently from
$\mu(\mathfrak{p}^{K})$. We can get the relations in $\Diff(M)$
operating in $U(\mathfrak{g})^{K}$ modulo
$(U(\mathfrak{g})\mathfrak{k})^{K}$. This approach leads to
simpler calculations than the operations through local coordinates
on $M$ (like in \cite{Reimann}) that gives quite cumbersome
calculations even in the relatively simple case of
$M=\mathbf{H}^{n}(\mathbb{R})_{ \mathbb{S}}$.

Below we are interested in the representation of the associative
algebra $\Diff(M)$ by it's generators and corresponding relations.
Let $\{g_{i}\}$ be a set of generators of the commutative
subalgebra $ \mathfrak{p}^{K}\subset S(\mathfrak{p})$. Not loosing
generality we can suppose that all $g_{i}$ are homogeneous
elements w.r.t. the grading of $S(\mathfrak{p})$. Then the
elements $\eta\circ\mu (g_{i})$ generate the algebra $\Diff(M)$.

Relations for the elements $\eta\circ\mu (g_{i})$ are of two
types. First type consists of relations induced by relations in
$U(\mathfrak{g})$. Due to the universality of $U( \mathfrak{g})$
all these relation are commutative ones, induced by the Lie
operation in $ \mathfrak{g}$. They are reduced to {\it commutative
relations} of the simplest form: $[D_{1},D_{2}]=\tilde D$, where
the operators $D_{1},D_{2}\in\Diff(M)$ have degrees $m_{1}$ and
$m_{2}$ respectively and the degree of $\tilde D\in\Diff(M)$ is
less or equal $m_{1}+m_{2}-1$.

Suppose now that there is a relation in $U(\mathfrak{g})$ of the
form:
\begin{equation*}
P(\eta\circ\mu( g_{1}),\ldots,\eta\circ\mu (g_{k}))=0
\end{equation*}
or equivalently
\begin{equation}\label{rel2}
P(\mu( g_{1}),\ldots,\mu (g_{k}))=\tilde D,
\end{equation}
where $P$ is a polynomial and $\tilde D\in
(U(\mathfrak{g})\mathfrak{k})^{K}$. Using the commutative
relations for $\mu(g_{i}),\,i=1,\ldots,k$ we can reduce the
polynomial $P$ to the polynomial $P_{1}$, symmetric w.r.t. all
permutations of it's arguments, and equation (\ref{rel2}) becomes:
\begin{equation}\label{rel22}
P_{1}(\mu(g_{1}),\ldots,\mu(g_{k}))=D^{*},
\end{equation}
$D^{*}\in (U(\mathfrak{g})\mathfrak{k})^{K}$. After this reduction
the relation (\ref{rel22}) may be trivial: $P_{1}=0,\,D^{*}=0$. It
means that (\ref{rel2}) is the commutative relation. Suppose that
relation (\ref{rel22}) is nontrivial. Consider the sum
$P_{2}(t_{1},\ldots,t_{k})$ of monomials with the highest total
degree from polynomial $P_{1}(t_{1},\ldots,t_{k})$ with
commutative variables $t_{1},\ldots,t_{k}$. Due to the symmetry of
$P_{1}$ the polynomial $P_{2}(t_{1},\ldots,t_{k})$ is nontrivial.
On the other hand from (\ref{rel22}) we obtain that
$P_{2}(g_{1},\ldots,g_{k})=0$ due to the expansion $
\mathfrak{g}=\mathfrak{p}\oplus\mathfrak{k}$. Thus every relation
in the algebra
$U(\mathfrak{g})^{K}/(U(\mathfrak{g})\mathfrak{k})^{K}$ w.r.t.
it's generators $\eta\circ\mu(g_{1}),\ldots,\eta\circ\mu(g_{k})$
modulo commutative relations is in one to one correspondence with
the relations for homogeneous generators $g_{1},\ldots,g_{k}$ of
commutative algebra $ \mathfrak{p}^{K}$. We call such relations
{\it the relations of the second type}.

The filtration of the algebra $U(\mathfrak{g})$ induces the
filtration of the algebra $\Diff(M)$ which coincides with the
natural filtration of $\Diff(M)$ as the algebra of differential
operators.

For simplicity throughout the whole paper we consider invariance
of differential operators only w.r.t. the identity component of a
whole isometry group.

\section{Two-point homogeneous Riemannian spaces}
\label{twopoint} \markright{\ref{twopoint} Two point homogeneous
Riemannian spaces}

The classification of two-point homogeneous Riemannian spaces can
been found in \cite{Tits}, \cite{Wang}, (see also
\cite{Matsumoto}, \cite{Wolf}, \cite{Chav}), and is as follows:
\begin{enumerate}
\item the Euclidean space $\mathbf{ E}^{n}$;
\item the sphere $\mathbf{ S}^{n}$;
\item the real projective space $\mathbf{P}^{n}(\mathbb{R})$;
\item the complex projective space $\mathbf{P}^{n}(\mathbb{C})$;
\item the quaternion projective space $\mathbf{P}^{n}(\mathbb{H})$;
\item the Cayley projective plane $\mathbf{P}^{2}(\mathbb{C}a)$;
\item the real hyperbolic space (Lobachevski space) $\mathbf{H}^{n}(\mathbb{R})$;
\item the complex hyperbolic space $\mathbf{H}^{n}(\mathbb{C})$;
\item the quaternion hyperbolic space  $\mathbf{H}^{n}(\mathbb{H})$;
\item the Cayley hyperbolic plane $\mathbf{H}^{2}(\mathbb{C}a)$.
\end{enumerate}
The isometry groups for all these spaces except Cayley planes are
classical and for the Cayley planes they are two real forms of the
complex special simple group $ \mathbf{F}_{4}$.

For the Lie algebra $\mathfrak{g}$ of the isometry group of the
compact two-point homogeneous space $M$ there is the following
general expansion \cite{Shchep2002},\cite{Hel78},\cite{Loos},
which is the specification of the expansion
$\mathfrak{g}=\mathfrak{p}\oplus\mathfrak{k}$ from section
\ref{InvOperators}.
\begin{proposit}\label{prop1}
The algebra $\mathfrak{g}$ admits the following expansion into the
direct sum of subspaces:
\begin{equation}\label{expan}
\mathfrak{g}=\mathfrak{a}\oplus\mathfrak{k}_{0}\oplus\mathfrak{k}_{\lambda}\oplus
\mathfrak{k}_{2\lambda}\oplus\mathfrak{p}_{\lambda}\oplus
\mathfrak{p}_{2\lambda},
\end{equation}
where $\dim\mathfrak{a}=1$, $\lambda$ is a nontrivial linear form
on the space $ \mathfrak{a}$,
$\dim\mathfrak{k}_{\lambda}=\dim\mathfrak{p}_{\lambda}=q_{1}$,
$\dim\mathfrak{k}_{2\lambda}=\dim\mathfrak{p}_{2\lambda}=q_{2}$,
$\mathfrak{p}=\mathfrak{a}\oplus\mathfrak{p}_{\lambda}\oplus\mathfrak{p}_{2\lambda}$,
$\mathfrak{k}=\mathfrak{k}_{0}\oplus\mathfrak{k}_{\lambda}\oplus\mathfrak{k}_{2\lambda};\,
q_{1},q_{2}\in \{0\}\cup \mathbb{N}$, the subalgebra
$\mathfrak{a}$ is the maximal commutative subalgebra in the
subspace $\mathfrak{p}$. Here $\mathfrak{k}$ is a stationary
subalgebra, corresponding to some point $x_{0}\in M$. All summands
in (\ref{expan}) are $\ad_{\mathfrak{k}_{0}}$-invariant; in
particular $\mathfrak{k}_{0}$ is a subalgebra of $\mathfrak{k}$.

Identified the space $ \mathfrak{p}$ with the space $T_{x_{0}}M$.
Under this identification the restriction of the Killing form for
the algebra $\mathfrak{g}$ onto the space $\mathfrak{p}$ and the
scalar product on $T_{x_{0}}M$ are proportional. In particular,
the decomposition $\mathfrak{g}=\mathfrak{p}\oplus\mathfrak{k}$ is
uniquely determined by the point $x_{0}$. Let
$\langle\cdot,\cdot\rangle$ be a scalar product on the algebra
$\mathfrak{g}$ such that it is proportional to the Killing form
and its restriction onto the subspace $ \mathfrak{p}\cong
T_{x_{0}}M$ coincides with the Riemannian metric $g$ on
$T_{x_{0}}M$. The spaces $\mathfrak{a},\,\mathfrak{k}_{0},\,
\mathfrak{k}_{\lambda},\,\mathfrak{p}_{\lambda},\,
\mathfrak{k}_{2\lambda},\,\mathfrak{p}_{2\lambda}$ are pairwise
orthogonal.

 Besides, the following inclusions are
valid:
\begin{align}\label{inclusions}
[\mathfrak{a},\mathfrak{p}_{\lambda}]\subset\mathfrak{k}_{\lambda},\;
[\mathfrak{a},\mathfrak{k}_{\lambda}]\subset\mathfrak{p}_{\lambda},\;
[\mathfrak{a},\mathfrak{p}_{2\lambda}]\subset\mathfrak{k}_{2\lambda},\;
[\mathfrak{a},\mathfrak{k}_{2\lambda}]\subset\mathfrak{p}_{2\lambda},\;
[\mathfrak{a},\mathfrak{k}_{0}]=0,\;
[\mathfrak{k}_{\lambda},\mathfrak{p}_{\lambda}]\subset\mathfrak{p}_{2\lambda}
\oplus\mathfrak{a},\nonumber \\
[\mathfrak{k}_{\lambda},\mathfrak{k}_{\lambda}]\subset\mathfrak{k}_{2\lambda}\oplus
\mathfrak{k}_{0},\;
[\mathfrak{p}_{\lambda},\mathfrak{p}_{\lambda}]\subset\mathfrak{k}_{2\lambda}\oplus
\mathfrak{k}_{0},\;
[\mathfrak{k}_{2\lambda},\mathfrak{k}_{2\lambda}]\subset\mathfrak{k}_{0},\;
[\mathfrak{p}_{2\lambda},\mathfrak{p}_{2\lambda}]\subset\mathfrak{k}_{0},\;
[\mathfrak{k}_{2\lambda},\mathfrak{p}_{2\lambda}]\subset\mathfrak{a},\nonumber\\
[\mathfrak{k}_{\lambda},\mathfrak{k}_{2\lambda}]\subset\mathfrak{k}_{\lambda},\;
[\mathfrak{k}_{\lambda},\mathfrak{p}_{2\lambda}]\subset\mathfrak{p}_{\lambda},\;
[\mathfrak{p}_{\lambda},\mathfrak{k}_{2\lambda}]\subset\mathfrak{p}_{\lambda},\;
[\mathfrak{p}_{\lambda},\mathfrak{p}_{2\lambda}]\subset\mathfrak{k}_{\lambda}.
\end{align} Moreover, for any basis $e_{\lambda,i},\;i=1,\dots,q_{1}$ in the
space $\mathfrak{p}_{\lambda}$ and any basis
$e_{2\lambda,i},\;i=1,\dots,q_{2}$ in the space
$\mathfrak{p}_{2\lambda}$ there are the basis
$f_{\lambda,i},\;i=1,\dots,q_{1}$ in the space
$\mathfrak{k}_{\lambda}$ and the basis
$f_{2\lambda,j},\;j=1,\dots,q_{2}$ in the space
$\mathfrak{k}_{2\lambda}$ such that:
\begin{align}\label{commSpesial1}
[\Lambda,e_{\lambda,i}]&=-\frac12f_{\lambda,i},\;[\Lambda,f_{\lambda,i}]=\frac12e_{\lambda,i},\;
[e_{\lambda,i},f_{\lambda,i}]=-\frac12\Lambda,\nonumber\\ \langle
e_{\lambda,i},e_{\lambda,j}\rangle&=\langle
f_{\lambda,i},f_{\lambda,j}\rangle=\delta_{ij}R^{2},\;
i,j=1,\dots,q_{1},\nonumber\\
[\Lambda,e_{2\lambda,i}]&=-f_{2\lambda,i},\;[\Lambda,f_{2\lambda,i}]=
e_{2\lambda,i},\; [e_{2\lambda,i},f_{2\lambda,i}]=-\Lambda,\\
\langle e_{2\lambda,i},e_{2\lambda,j}\rangle&=\langle
f_{2\lambda,i},f_{2\lambda,j}\rangle=\delta_{ij}R^{2},\;i,j=1,\dots,q_{2},
\;\langle\Lambda,\Lambda\rangle=R^{2}\nonumber,
\end{align}
where the vector $\Lambda\in\mathfrak{a}$ satisfies the conditions
$\langle\Lambda,\Lambda\rangle=R^{2},\,|\lambda(\Lambda)|=\dfrac12$.
Here the positive constant $R$ is connected with the maximal
sectional curvature $\varkappa_{m}$ of the space $M$ by the
formula $\varkappa_{m}=R^{-2}$.
\end{proposit}

Nonnegative integers $q_{1}$ and $q_{2}$ are said to be {\it
multiplicities of the space} $M$. They characterize $M$ uniquely
up to the exchange
$\mathbf{S}^{n}\leftrightarrow\mathbf{P}^{n}(\mathbb{R})$. For the
spaces $\mathbf{S}^{n}$ and $\mathbf{P}^{n}(\mathbb{R})$ we have
$q_{1}=0,\;q_{2}=n-1$; for the space
$\mathbf{P}^{n}(\mathbb{C}):\; q_{1}=2n-2,\;q_{2}=1$; for the
space $\mathbf{P}^{n}(\mathbb{H}):\; q_{1}=4n-4,\;q_{2}=3$; and
for the space $\mathbf{P}^{2}(\mathbb{C}a):\; q_{1}=8,\;q_{2}=7$.
Conversely, for the spaces $\mathbf{S}^{n}$ and
$\mathbf{P}^{n}(\mathbb{R})$ we could reckon that
$q_{1}=n-1,q_{2}=0$. Our choice corresponds to the isometries
$\mathbf{P}^{1}(\mathbb{C})\cong\mathbf{S}^{2},\,
\mathbf{P}^{1}(\mathbb{H})\cong\mathbf{S}^{4}$.

\begin{Rem}\label{transform}
Noncompact two-point homogeneous spaces of types 7,8,9,10 are
analogous to the compact two-point homogeneous spaces of types
2,4,5,6, respectively. In particular, it means that Lie algebras
$\mathfrak{g}$ of symmetry groups of analogous spaces are
different real forms of a simple complex Lie algebra. The
transition from one such real form to another can be done by
multiplying the subspace $\mathfrak{p}$ from the decomposition
$\mathfrak{g}=\mathfrak{k}\oplus\mathfrak{p}$ by the imaginary
unit $\mathbf{i}$ (or by $-\mathbf{i}$).
\end{Rem}

The information concerning the representation theory of symmetry
groups of two-point homogeneous spaces can be found in
\cite{Coll}.

\section{Invariant differential operators on $M_{\mathbb{S}}$} \label{spes}\markright{\ref{spes}
Invariant differential operators on $M_{\mathbb{S}}$}

Here we shall specify the construction from section
\ref{InvOperators} for the space $M_{\mathbb{S}}$, where $M$ is a
two-point homogeneous compact Riemannian space, using Proposition
\ref{prop1}.

Let $G$ be the identity component of the isometry group for $M$
and $K$ be its stationary subgroup, corresponding to the point
$x_{0}\in M$. The group $G$ naturally acts on the space
$M_{\mathbb{S}}$ and this action is transitive \cite{Wolf} (lemma
8.12.1). In particular $K$ acts transitively on the unit sphere
$\mathbb{S}_{x_{0}}\subset T_{x_{0}}M$. Identify the space
$\mathfrak{p}$ from Proposition \ref{prop1} with the space
$T_{x_{0}}M$. After this identification the action of $K$ on
$T_{x_{0}}M$ becomes the adjoint action $\Ad_{K}$ on
$\mathfrak{p}$. Let $K_{0}$ be the subgroup of $K$, corresponding
to the subalgebra $\mathfrak{k}_{0}\subset\mathfrak{k}$. Due to
relations (\ref{inclusions}),(\ref{commSpesial1}) $K_{0}$ is the
stationary subgroup of the group $G$, corresponding to the
point\footnote{By $(\alpha,\dots,\omega)$ we denote the set of
objects $\alpha,\dots,\omega$.} $y:=(x_{0},\Lambda')\in
M_{\mathbb{S}}$, where $\Lambda':=\dfrac1R\Lambda$.

Let
$\tilde{\mathfrak{p}}:=\mathfrak{a}\oplus\mathfrak{p}_{\lambda}\oplus\mathfrak{p}_{2\lambda}
\oplus\mathfrak{k}_{\lambda}\oplus\mathfrak{k}_{2\lambda}$. Since
$[\mathfrak{k}_{0},\tilde{\mathfrak{p}}]\subset\tilde{\mathfrak{p}}$
the expansion
$\mathfrak{g}=\tilde{\mathfrak{p}}\oplus\mathfrak{k}_{0}$ is
reductive. One has $T_{y}M_{\mathbb{S}}=T_{x_{0}}M\oplus
T_{\Lambda'}\mathbb{S}_{x_{0}}$. Due to Proposition \ref{prop1} we
obtain
$\Lambda'\bot\left(\mathfrak{p}_{\lambda}\oplus\mathfrak{p}_{2\lambda}\right)$
and $[f_{\lambda,i},\Lambda']=-(2R)^{-1}e_{\lambda,i},\,
i=1,\dots,q_{1},\,[f_{2\lambda,j},\Lambda']=-R^{-1}
e_{2\lambda,j},\, j=1,\dots,q_{2}$. Therefore the space
$\mathfrak{k}_{\lambda}\oplus\mathfrak{k}_{2\lambda}$ is
identified through the $K$-action on $T_{x_{0}}M$ with the space
$T_{\Lambda'}\mathbb{S}_{x_{0}}$ and the $K_{0}$-action on the
space
$T_{y}M_{\mathbb{S}}\simeq\tilde{\mathfrak{p}}=\mathfrak{a}\oplus\mathfrak{p}_{\lambda}\oplus\mathfrak{p}_{2\lambda}
\oplus\mathfrak{k}_{\lambda}\oplus\mathfrak{k}_{2\lambda}$ is
again adjoint.

From Proposition \ref{prop1} we see that $\Ad_{K_{0}}$ conserves
all summand in the last expansion. On the other hand the
$K_{0}$-action on $T_{\Lambda'}\mathbb{S}_{x_{0}}$ is the
differential of $K_{0}$-action on $(\Lambda')^{\bot}\subset
T_{x_{0}}M$. Since the last action is linear we obtain that the
$\Ad_{K_{0}}$-action in $\mathfrak{p}_{\lambda}$ is equivalent to
it's action in $\mathfrak{k}_{\lambda}$ and the
$\Ad_{K_{0}}$-action in $\mathfrak{p}_{2\lambda}$ is equivalent to
it's action in $\mathfrak{k}_{2\lambda}$.  Let
$\chi_{\lambda}:\,\mathfrak{k}_{\lambda}\to\mathfrak{p}_{\lambda},\,
\chi_{2\lambda}:\,\mathfrak{k}_{2\lambda}\to\mathfrak{p}_{2\lambda}$
be isomorphisms of linear spaces such that
$\left.\Ad_{K_{0}}\right|_{\mathfrak{p}_{\lambda}}\circ\chi_{\lambda}=\chi_{\lambda}\circ
\left.\Ad_{K_{0}}\right|_{\mathfrak{k}_{\lambda}}$ and
$\left.\Ad_{K_{0}}\right|_{\mathfrak{p}_{2\lambda}}\circ\chi_{2\lambda}=\chi_{2\lambda}\circ
\left.\Ad_{K_{0}}\right|_{\mathfrak{k}_{2\lambda}}$.

After the substitution
$\mathfrak{p}\to\tilde{\mathfrak{p}},\,\mathfrak{k}\to\mathfrak{k}_{0}$
we can apply the construction from section \ref{InvOperators} for
calculating generators and relations for the algebra
$\Diff(M_{\mathbb{S}})$. Let $g_{i}$ are independent invariants of
the $\Ad_{K_{0}}$-action in $\tilde{\mathfrak{p}}$. Then elements
$\eta\circ\mu(g_{i})$ generate the algebra
$\Diff(M_{\mathbb{S}})$.

There are some obvious invariants of the $\Ad_{K_{0}}$-action in
$\tilde{\mathfrak{p}}$. First of all it is
$\Lambda\in\mathfrak{a}$, since
$[\mathfrak{k}_{0},\mathfrak{a}]=0$. Secondly the
$\Ad_{K_{0}}$-action is isometric w.r.t. the Killing form which is
proportional to the scalar product $\langle\cdot,\cdot\rangle$, so
the $\Ad_{K_{0}}$-action conserves the restrictions of
$\langle\cdot,\cdot\rangle$ onto spaces
$\mathfrak{p}_{\lambda},\,\mathfrak{p}_{2\lambda},\,
\mathfrak{k}_{\lambda},\,\mathfrak{k}_{2\lambda}$. Thirdly the
$\Ad_{K_{0}}$-action conserves functions \newline
$\langle\chi_{\lambda}X_{1},Y_{1}\rangle,\,X_{1}\in\mathfrak{k}_{\lambda},\,
Y_{1}\in\mathfrak{p}_{\lambda}$ and
$\langle\chi_{2\lambda}X_{2},Y_{2}\rangle,\,X_{2}\in\mathfrak{k}_{2\lambda},\,
Y_{2}\in\mathfrak{p}_{2\lambda}$.

Let define $\chi_{\lambda}$ and $\chi_{2\lambda}$ by
\begin{equation}\label{chi}
\chi_{\lambda}X=2[\Lambda,X],\,X\in \mathfrak{k}_{\lambda},\;
\chi_{2\lambda}Y=[\Lambda,Y],\,Y\in \mathfrak{k}_{2\lambda}.
\end{equation}
It is obvious that
$\chi_{\lambda}:\mathfrak{k}_{\lambda}\mapsto\mathfrak{p}_{\lambda}$
and
$\chi_{2\lambda}:\mathfrak{k}_{2\lambda}\mapsto\mathfrak{p}_{2\lambda}$
are bijections. For any $k\in K_{0}$ it holds
$\Ad_{k}\Lambda=\Lambda$ since $[\mathfrak{k}_{0},
\mathfrak{a}]=0$, therefore  from (\ref{chi}) we obtain:
$\Ad_{k}\circ\chi_{\lambda}X=2[\Lambda,\Ad_{k}X]=\chi_{\lambda}\circ\Ad_{k}X,\,X\in\mathfrak{k}_{\lambda}$
and
$\Ad_{k}\circ\chi_{2\lambda}Y=[\Lambda,\Ad_{k}Y]=\chi_{2\lambda}\circ\Ad_{k}Y,\,Y\in\mathfrak{k}_{2\lambda}$.
It is clear that
$\chi_{\lambda}f_{\lambda,i}=e_{\lambda,i},\,i=1,\dots,q_{1}$ and
$\chi_{2\lambda}f_{2\lambda,j}=e_{2\lambda,j},\,j=1,\dots,q_{2}$.

The bases
$$
\{\frac1R e_{\lambda,i}\}_{i=1}^{q_{1}},\{\frac1R
f_{\lambda,i}\}_{i=1}^{q_{1}},\,\{\frac1R
e_{2\lambda,j}\}_{j=1}^{q_{2}},\{\frac1R
f_{2\lambda,j}\}_{j=1}^{q_{2}}
$$
respectively in spaces
$\mathfrak{p}_{\lambda},\mathfrak{k}_{\lambda},\mathfrak{p}_{2\lambda},\mathfrak{k}_{2\lambda}$
are orthonormal, so one has the following generators of
$\Diff(M_{\mathbb{S}})$:
\begin{align}\label{Ggen}
D_{0}&=\eta(\Lambda),\,D_{1}=\eta\left(\sum\limits_{i=1}^{q_{1}}e^{2}_{\lambda,i}\right),\,
D_{2}=\eta\left(\sum\limits_{i=1}^{q_{1}}f^{2}_{\lambda,i}\right),\,
D_{3}=\eta\left(\frac12\sum\limits_{i=1}^{q_{1}}\{e_{\lambda,i},f_{\lambda,i}\}\right),\,
\nonumber\\
D_{4}&=\eta\left(\sum\limits_{j=1}^{q_{2}}e^{2}_{2\lambda,j}\right),\,
D_{5}=\eta\left(\sum\limits_{j=1}^{q_{2}}f^{2}_{2\lambda,j}\right),\,
D_{6}=\eta\left(\frac12\sum\limits_{j=1}^{q_{2}}\{e_{2\lambda,j},f_{2\lambda,j}\}\right),
\end{align}
where $\{\cdot,\cdot\}$ means anticommutator. For brevity wee
shall omit the symbol $\eta$ below. From (\ref{commSpesial1}) one
easily obtains:
\begin{align*}
[D_{0},D_{1}]&=-D_{3},\,[D_{0},D_{2}]=D_{3},\,[D_{0},D_{3}]=\frac12
(D_{1}-D_{2}),\\ [D_{0},D_{4}]&=-2D_{6},\,[D_{0},D_{5}]=2D_{6},\,
[D_{0},D_{6}]=D_{4}-D_{5}.
\end{align*}

In order to find full system of invariants and relations in
$\Diff(M_{\mathbb{S}})$ we need more detailed information about
the $\Ad_{K_{0}}$-action in $\tilde{\mathfrak{p}}$ and commutators
in $\mathfrak{g}$. This information will be extracted in the
following sections from the models of two-point homogeneous
compact Riemannian spaces.

It is easily seen that every automorphism of Lie algebra
$\mathfrak{g}$, conserving its subalgebra $\mathfrak{k}_{0}$,
generates an automorphism of $\Diff(M_{\mathbb{S}})$. From
relations (\ref{inclusions}) one obtains that the map
$\sigma:\,\mathfrak{g}\to\mathfrak{g},\,\left.\sigma\right|_{\mathfrak{k}}=\id,\,
\left.\sigma\right|_{\mathfrak{p}}=-\id$ is the automorphism of
$\mathfrak{g}$. It generates the automorphism of
$\Diff(M_{\mathbb{S}})$: $D_{0}\to -D_{0},\,D_{1}\to D_{1},\,
D_{2}\to D_{2},\,D_{3}\to -D_{3},\,D_{4}\to D_{4},\, D_{5}\to
D_{5},\,D_{6}\to -D_{6}$. We shall denote it by the same symbol
$\sigma$.

Another obvious automorphism is the one parametric group
$\zeta_{\alpha}$ of internal automorphisms, generated by the
$\ad_{\Lambda}$-action. From (\ref{commSpesial1}) one obtains:
\begin{align*}
\zeta_{\alpha}(\Lambda)&=\Lambda,\,\zeta_{\alpha}(e_{\lambda,i})=\cos(\alpha/2)e_{\lambda,i}-
\sin(\alpha/2)f_{\lambda,i},\\
\zeta_{\alpha}(f_{\lambda,i})&=\sin(\alpha/2)e_{\lambda,i}+
\cos(\alpha/2)f_{\lambda,i},\,i=1,\dots,q_{1},\\
\zeta_{\alpha}(e_{2\lambda,j})&=\cos(\alpha)e_{2\lambda,j}-
\sin(\alpha)f_{2\lambda,j},\,
\zeta_{\alpha}(f_{2\lambda,i})=\sin(\alpha)e_{2\lambda,j}+
\cos(\alpha)f_{2\lambda,j},\,j=1,\dots,q_{2}.
\end{align*}
Therefore
\begin{align*}
\zeta_{\alpha}(D_{0})&=D_{0},\,
\zeta_{\alpha}(D_{1})=\cos^{2}(\alpha/2)D_{1}+\sin^{2}(\alpha/2)D_{2}-\sin(\alpha)
D_{3},\\
\zeta_{\alpha}(D_{2})&=\sin^{2}(\alpha/2)D_{1}+\cos^{2}(\alpha/2)D_{2}+\sin\alpha
D_{3},\\ \zeta_{\alpha}(D_{3})&=\frac12\sin(\alpha)
(D_{1}-D_{2})+\cos(\alpha)D_{3},\\
\zeta_{\alpha}(D_{4})&=\cos^{2}(\alpha)D_{4}+\sin^{2}(\alpha)
D_{5}-\sin(2\alpha) D_{6},\\
\zeta_{\alpha}(D_{5})&=\sin^{2}(\alpha)D_{4}+\cos^{2}(\alpha)
D_{5}+\sin(2\alpha) D_{6},\\
\zeta_{\alpha}(D_{6})&=\frac12\sin(2\alpha)
(D_{4}-D_{5})+\cos(2\alpha) D_{6}.
\end{align*}
In particular
$\zeta_{\pi}(D_{1})=D_{2},\,\zeta_{\pi}(D_{2})=D_{1},\,\zeta_{\pi}(D_{3})=-D_{3},\,
\zeta_{\pi}(D_{i})=D_{i},\,i=0,4,5,6$.

Proposition \ref{prop1} imply that the base
$$
\frac1R\Lambda,\frac1R e_{\lambda,i},\frac1R f_{\lambda,i},\frac1R
e_{2\lambda,j},\frac1R
f_{2\lambda,j},\,i=1,\dots,q_{1},j=1,\dots,q_{2}
$$
in the space $\tilde{\mathfrak{p}}$ is orthonormal, therefore the
operator $D^{*}=D_{0}^{2}+D_{1}+D_{2}+D_{4}+D_{5}$ corresponds to
the Casimir operator in $U(\mathfrak{g})$ \cite{Post} (lecture
18). It means that $D^{*}$ lies in the centre of the algebra
$\Diff(M_{\mathbb{S}})$.

Let $\pi_{1}:\, M_{\mathbb{S}}\to M$ be the canonical projection
and $\tilde\pi_{1}$ is the map of a function $f$ on $M$ to the
function $f\circ\pi_{1}$ on $M_{\mathbb{S}}$. Due to the
identification $\mathfrak{p}\simeq T_{x_{0}}M$ it is clear that
the operator $(D_{0}^{2}+D_{1}+D_{2})\circ\tilde\pi_{1}$ is the
Laplace-Beltrami operator on $M$.

\section{The model for the space $\mathbf{P}^{n}(\mathbb{H})$} \label{QuatMod}\markright{\ref{QuatMod}
The model for the space $\mathbf{P}^{n}(\mathbb{H})$}

Let $\mathbb{H}$ be the quaternion algebra over the field $
\mathbb{R}$ with the base $1,\ii,\jj,\kk$, where
$\ii\jj=-\jj\ii=\kk,\,\jj\kk=-\kk\jj=\ii,\,\kk\ii=-\ii\kk=\jj$.
The conjugation acts as follows:
$\overline{x+y\ii+z\jj+t\kk}=x-y\ii-z\jj-t\kk,\,x,y,z,t\in\mathbb{R}$.

Let $ \mathbb{H}^{n+1}$ be the right quaternion space and
$(z_{1},\ldots,z_{n+1})$ be coordinates on it. Let
$\mathbf{P}^{n}(\mathbb{H})$ be a factor space of the space $
\mathbb{H}^{n+1}\backslash\{0\}$ with respect to the right action
of the multiplicative group
$\mathbb{H}^{*}=\mathbb{H}\backslash\{0\}$. The set
$(z_{1}:\ldots:z_{n+1})$ up to the multiplication from the right
by an arbitrary element from the group $\mathbb{H}^{*}$ is the set
of homogeneous coordinates for the element\footnote{To distinguish
the point $x\in M$ from their coordinates we shall single out it
by the bold type} $\pi(\z)$ on the space
$\mathbf{P}^{n}(\mathbb{H})$, where $\pi :
\mathbb{H}^{n+1}\backslash\{0\}\rightarrow
\mathbf{P}^{n}(\mathbb{H})$ is the canonical projection. Let
$\langle\x,\y\rangle:=\sum\limits_{i=1}^{n+1}\bar x_{i}
y_{i},\,\x=(x_{1},\ldots,x_{n+1}),\,\y=(y_{1},\ldots,y_{n+1})\in
\mathbb{H}^{n+1}$ be the standard scalar product in the space
$\mathbb{H}^{n+1}$. Let
$\z\in\mathbb{H}^{n+1}\backslash\{0\},\,{\bf\xi}_i\in
T_{\z}\mathbb{H}^{n+1},\,{\bf \zeta}_{i}=\pi_{*}{\bf\xi}_{i}\in
T_{\pi(\z)}\left(\mathbf{P}^{n}(\mathbb{H})\right),\,i=1,2$. A
metric on the space $\mathbf{P}^{n}(\mathbb{H})$
\begin{equation}\label{metricPH}
 \tilde g|_{\z}({\bf\zeta}_{1},{\bf\zeta}_{2})=\left(\langle{\bf\xi}_{1},{\bf\xi}_{2}\rangle\langle\z,\z\rangle-
\langle{\bf\xi}_{1},\z\rangle\langle\z,{\bf\xi}_{2}\rangle\right)/\langle\z,\z\rangle^{2},
\end{equation}
is the analog for the metric with a constant sectional curvature
on the space $\mathbf{P}^{n}(\mathbb{R})$ and the metric with a
constant holomorphic sectional curvature on the space
$\mathbf{P}^{n}(\mathbb{C})$. The real part of the metric
(\ref{metricPH}) is a Riemannian metric on the space
$\mathbf{P}^{n}(\mathbb{H})$:
\begin{equation}\label{metricPRH}
g=4R^{2}\RE\tilde g.
\end{equation}

The normalizing factor in (\ref{metricPRH}) is chosen due to the
following reasons. The space $\mathbf{P}^{1}(\mathbb{H})$ with
this metric is the sphere ${\bf S}^{4}$ with the standard metric
of the constant sectional curvature $R^{-2}$. To see this we can
consider a homeomorphism
$\nu:\mathbf{P}^{1}(\mathbb{H})\rightarrow
\bar{\mathbb{H}}\cong{\bf S}^{4},\,
\nu(z_{1},z_{2})=z_{1}\left(z_{2}\right)^{-1}=z\in\overline{\mathbb{H}}$,
where $\overline{\mathbb{H}}$ is the quaternion space completed
with the point at infinity. For $n=1$ the formula
(\ref{metricPRH}) has the form
\begin{equation} \label{n1}
g=4R^{2}\frac{(d\bar z_{1}dz_{1}+d\bar
z_{2}dz_{2})(|z_{1}|^{2}+|z_{2}|^{2})-(d\bar z_{1}\cdot
z_{1}+d\bar z_{2}\cdot z_{2}) (\bar z_{1}dz_{1}+\bar
z_{2}dz_{2})}{(|z_{1}|^{2}+|z_{2}|^{2})^{2}}.
\end{equation}
Using the formula $|z_{2}|^{2}dz_{1}-z_{1}\bar
z_{2}dz_{2}=|z_{2}|^{2}(dz)z_{2}$ by direct calculations we can
reduce the expression (\ref{n1}) to the form:
$$
g=\frac{4R^{2}dzd\bar z}{\left(1+|z|^{2}\right)^{2}}
$$
which is the metric with the constant sectional curvature $R^{-2}$
on the sphere $\mathbf{ S}^{4}$.

The left action of the group $\U_{\mathbb{H}}(n+1)$, consisting of
quaternion matrices $A$ of the size $(n+1)\times (n+1)$ such that
$\bar A^{T}A=E$, conserves the scalar product
$\langle\cdot,\cdot\rangle$ in the space $
\mathbb{H}^{n+1},\,\dim_{\mathbb{R}}\U_{\mathbb{H}}(n+1)=(2n+3)(n+1)$.
If we write every quaternion coordinates in $\mathbb{H}^{n+1}$ as
a pair of complex numbers, then the group $\U_{\mathbb{H}}(n+1)$
becomes the symplectic group $\Simp(n+1)$.

Left and right multiplications always commute, so the left action
of the group $\U_{\mathbb{H}}(n+1)$ is correctly defined also on
the space $\mathbf{P}^{n}(\mathbb{H})$. Obviously, it is
transitive and conserves the metric $g$. The stationary subgroup
of the point from the space $\mathbf{P}^{n}(\mathbb{H})$ with the
homogeneous coordinates $(1,0,\ldots,0)$ is the group
$\U_{\mathbb{H}}(n)\U_{\mathbb{H}}(1)$, where the group
$\U_{\mathbb{H}}(n)$ acts onto the last $n$ coordinates, and the
group $\U_{\mathbb{H}}(1)$ acts by the left multiplication of all
homogeneous coordinates by quaternions with the unit norm. All
stationary subgroups on a homogeneous space are conjugated and
hence isomorphic. Therefore
$\mathbf{P}^{n}(\mathbb{H})=\U_{\mathbb{H}}(n+1)/(\U_{\mathbb{H}}(n)\U_{\mathbb{H}}(1))$.

The Lie algebra $\mathfrak{u}_{\mathbb{H}}(n+1)$ consist of
quaternion matrices $A$ of the size $(n+1)\times(n+1)$ such that
$\bar A^{T}=-A$. Let $E_{kj}$ be the matrix of the size
$(n+1)\times(n+1)$ with the unique nonzero element equals $1$,
locating at the intersection of the $k$-th row and the $j$-th
column. Choose the base for the algebra
$\mathfrak{u}_{\mathbb{H}}(n+1)$ as:
\begin{align}
\label{BasisUH} \Psi_{kj}=\frac12(E_{kj}-E_{jk}),\,1\leqslant
k<j\leqslant n+1,\,
\Upsilon_{kj}=\frac{\ii}2(E_{kj}+E_{jk}),\nonumber\\
\Omega_{kj}=\frac{\jj}2(E_{kj}+E_{jk}),
\,\Theta_{kj}=\frac{\kk}2(E_{kj}+E_{jk}),\,1\leqslant k\leqslant
j\leqslant n+1.
\end{align}
The commutative relations for these elements are:
\begin{align}
\label{commutators1}
[\Psi_{kj},\Psi_{ml}]&=\frac12\left(\delta_{jm}\Psi_{kl}-\delta_{km}\Psi_{jl}+\delta_{kl}
\Psi_{jm}-\delta_{jl}\Psi_{km}\right),\nonumber \\
\left[\Psi_{kj},\Upsilon_{ml}\right]&=\frac12\left(\delta_{jm}\Upsilon_{kl}-\delta_{km}
\Upsilon_{jl}+\delta_{lj}\Upsilon_{km}-\delta_{lk}\Upsilon_{jm}\right),\nonumber
\\
[\Upsilon_{kj},\Upsilon_{ml}]&=\frac12\left(\delta_{jm}\Psi_{lk}+\delta_{km}\Psi_{lj}+
\delta_{kl}\Psi_{mj}+\delta_{jl}\Psi_{mk}\right),\nonumber\\
[\Upsilon_{kj},\Omega_{ml}]&=\frac12\left(\delta_{jm}\Theta_{lk}+\delta_{km}\Theta_{lj}+
\delta_{kl}\Theta_{mj}+\delta_{jl}\Theta_{mk}\right),
\end{align}
plus the analogous equalities, obtaining from the latter three
relations by the cyclic permutation $\Upsilon\rightarrow
\Omega\rightarrow \Theta\rightarrow \Upsilon $, where
$\Psi_{kj}=-\Psi_{jk},\,\Psi_{kk}=0,\,\Upsilon_{kj}=\Upsilon_{jk},\,\Omega_{kj}=\Omega_{jk},\,
\Theta_{kj}=\Theta_{jk}$.

\section{Algebras $\Diff(\mathbf{P}^{n}(\mathbb{H})_{\mathbb{S}})$ and
$\Diff(\mathbf{H}^{n}(\mathbb{H})_{\mathbb{S}})$}
\label{Quat}\markright{\ref{Quat} Algebras
$\Diff(\mathbf{P}^{n}(\mathbb{H})_{\mathbb{S}})$ and
$\Diff(\mathbf{H}^{n}(\mathbb{H})_{\mathbb{S}})$}

Consider now the total space of unit spheres bundle
$\mathbf{P}^{n}(\mathbb{H})_{\mathbb{S}}$ over the space
$\mathbf{P}^{n}(\mathbb{H})$. Let $(\z,{\bf \zeta})$, where
$\z\in\mathbf{P}^{n}(\mathbb{H}),\,{\bf \zeta}\in
T_{\z}\mathbf{P}^{n}(\mathbb{H})$ be a general point of the space
$\mathbf{P}^{n}(\mathbb{H})_{\mathbb{S}}$. Due to the isomorphism
$\mathbf{P}^{1}(\mathbb{H})\cong{\bf S}^{4}$ we assume here
$n\geqslant 2$.

Let $\tilde \z_{0}=(1,0,\ldots,0)\in \mathbb{H}^{n+1}$, an element
${\bf\xi}_{0}\in T_{\tilde \z_{0}}\mathbb{H}^{n+1}\cong
\mathbb{H}^{n+1}$ has coordinates $(0,1,0,\ldots,0)$. Put
$\z_{0}=\pi\tilde \z_{0},\,{\bf \zeta}_{0}=\pi_{*}{\bf \xi}_{0}\in
T_{\z_{0}}\mathbf{P}^{n}(\mathbb{H})$.

The stationary subgroup $K_{0}$ of the group
$\U_{\mathbb{H}}(n+1)$, corresponding to the point
$(\z_{0},\zeta_{0})$ is generated by the group
$K_{1}=\U_{\mathbb{H}}(n-1)$, acting onto the last $(n-1)$-th
homogeneous coordinates and by the group
$K_{2}=\U_{\mathbb{H}}(1)$, acting by the left multiplication of
all homogeneous coordinates by quaternions with the unit norm. The
algebra $\mathfrak{k}_{0}$ of the group $K_{0}$ (corresponding to
Proposition \ref{prop1}) is $(2n^{2}-3n+4)$-dimensional and is
generated by elements (\ref{BasisUH}) with $3\leqslant k \leqslant
j\leqslant n+1$ and the elements:
$$
\sum_{k=1}^{n+1}\Upsilon_{kk},\,\sum_{k=1}^{n+1}\Omega_{kk},\sum_{k=1}^{n+1}\Theta_{kk}.
$$
Choose the complimentary subspace $\tilde{\mathfrak{p}}$ to the
subalgebra $\mathfrak{k}_{0}$ in the algebra
$\mathfrak{g}=\uu_{\mathbb{H}}(n+1)$ as the linear hull of
elements:
\begin{align}\label{subspacePH}
\Psi_{1k},\,\Upsilon_{1k},\,\Omega_{1k},\,\Theta_{1k},\,2\leqslant
k\leqslant n+1,\,
\Psi_{2k},\,\Upsilon_{2k},\,\Omega_{2k},\,\Theta_{2k},\,3\leqslant
k\leqslant n+1,\nonumber \\
\Upsilon_{*}=\frac{\ii}2\left(E_{11}-E_{22}\right),\,
\Omega_{*}=\frac{\jj}2\left(E_{11}-E_{22}\right),\,
\Theta_{*}=\frac{\kk}2\left(E_{11}-E_{22}\right).
\end{align}
Taking into account relations (\ref{commutators1}) it is easily
obtained that the expansion
$\mathfrak{u}_{\mathbb{H}}(n+1)=\tilde{\mathfrak{p}}\oplus\mathfrak{k}_{0}$
is reductive, i.e.
$[\tilde{\mathfrak{p}},\mathfrak{k}_{0}]\subset\tilde{\mathfrak{p}}$.

It is readily seen from (\ref{commutators1}) that setting:
\begin{align}\label{reperH}
\Lambda&=-\Psi_{12},\,e_{\lambda,k-2}=\Psi_{1k},\,e_{\lambda,n-3+k}=\Upsilon_{1k},\,
e_{\lambda,2n-4+k}=\Omega_{1k},\,e_{\lambda,3n-5+k}=\Theta_{1k},\notag
\\
f_{\lambda,k-2}&=-\Psi_{2k},\,f_{\lambda,n-3+k}=-\Upsilon_{2k},\,
f_{\lambda,2n-4+k}=-\Omega_{2k},\,f_{\lambda,3n-5+k}=-\Theta_{2k},\,k=3,\dots,n+1,\notag
\\
e_{2\lambda,1}&=\Upsilon_{12},\,e_{2\lambda,2}=\Omega_{12},\,e_{2\lambda,3}=\Theta_{12},\,
f_{2\lambda,1}=\Upsilon_{*},\,f_{2\lambda,2}=\Omega_{*},\,f_{2\lambda,3}=\Theta_{*},
\end{align}
we obtain the base from Proposition \ref{prop1} for
$q_{1}=4n-4,q_{2}=3$.

Now we are to find the full set of independent
$\Ad_{K_{0}}$-invariants in $S(\tilde{\mathfrak{p}})$. According
to section \ref{spes}, the expansion
$\tilde{\mathfrak{p}}=\mathfrak{a}\oplus\mathfrak{k}_{\lambda}\oplus
\mathfrak{k}_{2\lambda}\oplus\mathfrak{p}_{\lambda}\oplus
\mathfrak{p}_{2\lambda}$ is invariant w.r.t. the
$\Ad_{K_{0}}$-action. In the space $ \mathfrak{a}$ the
$K_{0}$-action is trivial that gives the invariant
$D_{0}=\Lambda\in\mu(\mathfrak{p}^{K_{0}})$, already found in
section \ref{spes}.

From formulas (\ref{reperH}) we see that the space
$\mathfrak{p}_{\lambda}\cong \mathbb{H}^{n-1}$ consists of
matrices of the form
$$
\left(\begin{array}{cc} 0 & - a^{*} \\ a & 0
\end{array}\right)\equiv\left(\begin{array}{ccccc} 0 & 0 & -\bar a_{1} & \ldots & -\bar
a_{n-1} \\ 0 & 0 & 0 & \ldots & 0 \\ a_{1} & 0 & 0 & \ldots & 0 \\
\vdots & \vdots & \vdots & \ddots & \vdots \\ a_{n-1} & 0 & 0 &
\ldots & 0
\end{array}\right),\,a_{1},\ldots,a_{n-1}\in\mathbb{H}.
$$
Likewise, the space $\mathfrak{k}_{\lambda}\cong \mathbb{H}^{n-1}$
consist of matrices of the form
$$
\left(\begin{array}{ccc} 0 & 0 & 0 \\ 0 & 0 & - b^{*} \\ 0 & b & 0
\end{array}\right)\equiv\left(\begin{array}{ccccc} 0 & 0 & 0 & \ldots & 0 \\ 0 & 0 & -\bar
b_{1} & \ldots & -\bar b_{n-1} \\ 0 & b_{1} & 0 & \ldots & 0
\\ \vdots & \vdots & \vdots & \ddots & \vdots \\ 0 & b_{n-1} & 0 &
\ldots & 0
\end{array}\right),\,b_{1},\ldots,b_{n-1}\in\mathbb{H}.
$$
Due to the formula
$$
\left(\begin{array}{cc} 1 & 0 \\ 0 & U
\end{array}\right)\left(\begin{array}{cc} 0 & - a^{*} \\ a & 0
\end{array}\right)\left(\begin{array}{cc} 1 & 0 \\ 0 & U^{*}
\end{array}\right)=\left(\begin{array}{cc} 0 & - (Ua)^{*} \\ Ua & 0
\end{array}\right),\, U\in\U_{\mathbb{H}}(n-1), a\in\mathbb{H}^{n-1}
$$
the action of the group $K_{1}$ in the space
$\mathfrak{p}_{\lambda}$ is equivalent to the standard action of
the group $U_{\mathbb{H}}(n-1)$ in the space $\mathbb{H}^{n-1}:\,
a\rightarrow Ua$. In the space $\mathfrak{k}_{\lambda}$ the action
of $K_{1}$ is similar: $b\rightarrow Ub$.

The standard action of the group $\U_{\mathbb{H}}(n-1)$ in the
space $\mathbb{H}^{n-1}$ has one independent real invariant:
$\langle\z,\z\rangle,\, \z\in \mathbb{H}^{n-1}$, and the diagonal
action of $\U_{\mathbb{H}}(n-1)$ in the space
$\mathfrak{p}_{\lambda}\oplus\mathfrak{k}_{\lambda}\cong
\mathbb{H}^{n-1}\oplus\mathbb{H}^{n-1}$ has six (independent iff
$n\geqslant 3$) real invariants:
\begin{equation}\label{ComInv}
\langle\z_{1},\z_{1}\rangle\in
\mathbb{R},\,\langle\z_{2},\z_{2}\rangle\in
\mathbb{R},\,\langle\z_{1},\z_{2}\rangle\in
\mathbb{H}\cong\mathbb{R}^{4},\,\z_{1},\z_{2}\in\mathbb{H}^{n-1}.
\end{equation}
Denote the corresponding elements from
$\mu(\tilde{\mathfrak{p}}^{K_{1}})\in U(\mathfrak{g})^{K_{1}}$ in
the following way:
\begin{align}\label{def1}
D_{1}&=\sum\limits_{k=3}^{n+1}\left(\Psi_{1k}^{2}+\Upsilon_{1k}^{2}+\Omega_{1k}^{2}+\Theta_{1k}^{2}\right),\,
D_{2}=\sum\limits_{k=3}^{n+1}\left(\Psi_{2k}^{2}+\Upsilon_{2k}^{2}+\Omega_{2k}^{2}+\Theta_{2k}^{2}\right),\nonumber\\
D_{3}&=-\frac12\sum\limits_{k=3}^{n+1}\left(\left\{\Psi_{1k},\Psi_{2k}\right\}+\left\{\Upsilon_{1k},\Upsilon_{2k}\right\}+
\left\{\Omega_{1k},\Omega_{2k}\right\}+\left\{\Theta_{1k},\Theta_{2k}\right\}\right),\nonumber\\
\square_{1}&=\frac12\sum\limits_{k=3}^{n+1}\left(-\left\{\Psi_{1k},\Upsilon_{2k}\right\}+
\left\{\Psi_{2k},\Upsilon_{1k}\right\}+
\left\{\Theta_{1k},\Omega_{2k}\right\}-\left\{\Theta_{2k},\Omega_{1k}\right\}\right),\\
\square_{2}&=\frac12\sum\limits_{k=3}^{n+1}\left(-\left\{\Psi_{1k},\Omega_{2k}\right\}+
\left\{\Psi_{2k},\Omega_{1k}\right\}+
\left\{\Upsilon_{1k},\Theta_{2k}\right\}-\left\{\Upsilon_{2k},\Theta_{1k}\right\}\right),\nonumber\\
\square_{3}&=\frac12\sum\limits_{k=3}^{n+1}\left(-\left\{\Psi_{1k},\Theta_{2k}\right\}+
\left\{\Psi_{2k},\Theta_{1k}\right\}+
\left\{\Omega_{1k},\Upsilon_{2k}\right\}-\left\{\Omega_{2k},\Upsilon_{1k}\right\}\right)\nonumber.
\end{align}

If $n=2$, then there is the unique independent relation between
invariants (\ref{ComInv}):
\begin{equation}\label{EulerIdentity}
|\langle \z_{1},\z_{2}\rangle|^{2}=|\bar
z_{1}z_{2}|^{2}=|z_{1}|^{2}\,|z_{2}|^{2}=\langle
\z_{1},\z_{1}\rangle\langle
\z_{2},\z_{2}\rangle,\,\z_{1}=z_{1},\z_{2}=z_{2}\in \mathbb{H}.
\end{equation}
If we write this identity in coordinates, then we will obtain the
well known Euler identity which is the key ingredient in the proof
of the Lagrange theorem from number theory: {\it if two integers
have the form $a^{2}+b^{2}+c^{2}+d^{2},\,a,b,c,d\in \mathbb{Z}$,
then their product has the same form}.

The elements $D_{1},D_{2},D_{3}$, already found in section
\ref{spes}, are invariant w.r.t. the action of the whole group
$K_{0}$, therefore they correspond to operators of the second
order from $\Diff(\mathbf{P}^{n}(\mathbb{H})_{\mathbb{S}})$. The
elements $\square_{1},\square_{2},\square_{3}$ are not invariant
w.r.t. the action of the group $K_{2}\cong \U_{ \mathbb{H}}(1)$.
Obviously, the $K_{2}$-action on the linear hull of elements
$\square_{1},\square_{2},\square_{3}$ is equivalent to the well
known action of the group $\SO (3)\cong \U_{\mathbb{H}}(1)/(1,-1)$
in the space $\mathbb{H}'$ of pure imaginary quaternions:
$$
x\rightarrow qx\bar q,\,x\in\mathbb{H}',\,q\in\U_{\mathbb{H}}(1),
$$
after the identification
$\square_{1}\leftrightarrow\ii,\,\square_{1}\leftrightarrow\jj,\,
\square_{3}\leftrightarrow\kk$.

The $K_{2}$-action on 3-dimensional spaces
$\mathfrak{p}_{2\lambda},\,\mathfrak{k}_{2\lambda}$ are the same
after the identification
$\Upsilon_{12},\Upsilon_{*}\leftrightarrow\ii;\,\Omega_{12},\Omega_{*}\leftrightarrow\jj;\,
\Theta_{12},\Theta_{*}\leftrightarrow\kk$; while the
$K_{1}$-action in these spaces is trivial. Thus we are to find
invariants of diagonal action of the group $\SO (3)$ in the space
$\mathbb{R}^{3}\oplus\mathbb{R}^{3}\oplus\mathbb{R}^{3}$. It is
clear that there are $6=9-3$ such independent invariants:
$$\langle \x,\x\rangle,\,\langle \y,\y\rangle,\,\langle \z,\z\rangle,\,\langle \x,\y\rangle,\,\langle \x,\z\rangle,
\,\langle \z,\y\rangle,\,\x,\y,\z\in \mathbb{R}^{3}$$ and
invariant $\langle \x,\y,\z\rangle\equiv\langle \x,\y\times
\z\rangle$ algebraically connected with the first six:
\begin{align}\label{ComRel}
\langle \x,\y,\z\rangle^{2}&=\x^{2}\y^{2}\z^{2}+2\langle
\x,\y\rangle\langle \x,\z\rangle\langle \y,\z\rangle-\x^{2}\langle
\y,\z\rangle^{2}\nonumber\\&-\y^{2}\langle
\x,\z\rangle^{2}-\z^{2}\langle \x,\y\rangle^{2},
\end{align}
where $\y\times\z$ is the standard vector product in
$\mathbb{R}^{3}$. Relation (\ref{ComRel}) can be verified using
the well known formulas: $\langle
\x,\y\rangle^{2}=\x^{2}\y^{2}-\langle\x\times\y\rangle^{2}$ and
$\x\times (\y\times\x)=\langle \x,\z\rangle \y-\langle
\x,\y\rangle \z$.

It gives the following invariants from  $U(\mathfrak{g})^{K_{0}}$:
\begin{align}\label{def2}
D_{4}&=\Upsilon_{12}^{2}+\Omega_{12}^{2}+\Theta_{12}^{2},\,
D_{5}=\Upsilon_{*}^{2}+\Omega_{*}^{2}+\Theta_{*}^{2},\nonumber\\
D_{6}&=\frac12\left(\{\Upsilon_{12},\Upsilon_{*}\}+\{\Omega_{12},\Omega_{*}\}+\{\Theta_{12},\Theta_{*}\}\right),
\nonumber\\
D_{7}&=\frac12\left(\{\square_{1},\Upsilon_{12}\}+\{\square_{2},\Omega_{12}\}+\{\square_{3},\Theta_{12}\}\right),
\\
D_{8}&=\frac12\left(\{\square_{1},\Upsilon_{*}\}+\{\square_{2},\Omega_{*}\}+\{\square_{3},\Theta_{*}\}\right),\,
D_{9}=\square_{1}^{2}+\square_{2}^{2}+\square_{3}^{2},\nonumber\\
D_{10}&=\square_{1}\Omega_{12}\Theta_{*}-\square_{1}\Omega_{*}\Theta_{12}+
\square_{2}\Upsilon_{*}\Theta_{12}-\square_{2}\Upsilon_{12}\Theta_{*}+
\square_{3}\Omega_{*}\Upsilon_{12}-\square_{3}\Omega_{12}\Upsilon_{*}\nonumber.
\end{align}
Here we took into account that every three factors from all
summand in the last expression pairwise commutate. The invariants
$D_{4}, D_{5}, D_{6}$ correspond to the general case, considered
in section \ref{spes}.

In fact invariants $D_{7},D_{8},D_{9}$ and $D_{10}$ are not in
$\mu\left(\mathfrak{p}^{K_{0}}\right)$ because they are not
symmetric w.r.t. all transposition of their factors of the first
degree. After complete symmetrization we can obtain invariants
from $\mu\left(\mathfrak{p}^{K_{0}}\right):\, \tilde D_{k}\equiv
D_{k}+D_{k}^{*}
\mod\left(U(\mathfrak{g})\mathfrak{k}\right)^{K_{0}},\,k=7,8,9,10$,
where $D_{k}^{*}$ are elements from $U(\mathfrak{g})^{K_{0}}$ with
$\deg D_{k}^{*}<\deg D_{k}$. For convenience we will use elements
$D_{k}$ instead of $\tilde D_{k},\,k=7,8,9,10$.

Thus operators $D_{0},\ldots,D_{10}$ generate the algebra
$\Diff(\mathbf{P}^{n}(\mathbb{H})_{\mathbb{S}})$.

The degrees of the generators are as follows:
\begin{align}\label{degrees}
\degr(D_{0})&=1,\,\degr(D_{1})=\degr(D_{2})=\degr(D_{3})=\degr(D_{4})=\degr(D_{5})=\degr(D_{6})=2,\nonumber\\
\degr(D_{7})&=\degr(D_{8})=3,\,\degr(D_{9})=\degr(D_{10})=4.
\end{align}

In the model of the space $\mathbf{P}^{n}(\mathbb{H})$ we can
transpose the coordinates $z_{1}$ and $z_{2}$. The operators
$D_{3},D_{4},D_{5},D_{8},D_{9},D_{10}$  are symmetric (invariant)
w.r.t. this transposition and the operators
$D_{0},\square_{1},\square_{2},\square_{3},D_{6},D_{7}$ are skew
symmetric. The operators $D_{1}$ and $D_{2}$ turn into each other
under this transposition.

It is easily verified that automorphisms $\zeta_{\alpha},\sigma$
acts on $\square_{i},D_{7},\dots,D_{10},i=1,2,3$ as
\begin{align*}
\zeta_{\alpha}(\square_{i})&=\square_{i},i=1,2,3,\,\zeta_{\alpha}(D_{7})=\cos(\alpha)
D_{7}-\sin(\alpha) D_{8},\,\zeta_{\alpha}(D_{8})=\sin(\alpha)
D_{7}+\cos(\alpha)D_{8},\nonumber\\
\zeta_{\alpha}(D_{9})&=D_{9},\,\zeta_{\alpha}(D_{10})=D_{10},
\sigma(\square_{i})=-\square_{i},i=1,2,3,\,\sigma(D_{7})=D_{7},
\,\sigma(D_{8})=-D_{8},\nonumber\\
\sigma(D_{9})&=D_{9},\,\sigma(D_{10})=D_{10}.
\end{align*}
Taking into account their action on other generators (see section
\ref{spes}) we obtain that the transposition of $z_{1}$ and
$z_{2}$ is equivalent to the composition $\sigma\circ\zeta_{\pi}$.

In order to get the generators of the algebra
$\Diff(\mathbf{H}^{n}(\mathbb{H})_{\mathbb{S}})$ one can use
Remark \ref{transform}, formula (\ref{reperH}) and make the formal
substitution:
\begin{align*}
\Lambda &\rightarrow
\ii\Lambda,\,\Psi_{1k}\rightarrow\ii\Psi_{1k},\,\Upsilon_{1k}\rightarrow\ii\Upsilon_{1k},\,\Omega_{1k}\rightarrow\ii
\Omega_{1k},\,\Theta_{1k}\rightarrow\ii\Theta_{1k},\,\Upsilon_{12}\rightarrow\ii\Upsilon_{12},\\
\Omega_{12}&\rightarrow \ii\Omega_{12},\,
\Theta_{12}\rightarrow\ii\Theta_{12},\,\Psi_{2k}\rightarrow\Psi_{2k},\,\Upsilon_{2k}\rightarrow\Upsilon_{2k},\,\Omega_{2k}
\rightarrow\Omega_{2k},\,\Theta_{2k}\rightarrow\Theta_{2k},\\
\Upsilon_{*}&\rightarrow\Upsilon_{*},\,
\Omega_{*}\rightarrow\Omega_{*},\,\Theta_{*}\rightarrow\Theta_{*},\,k=3,\ldots,n+1.
\end{align*}
This substitution produces the following substitution for the
generators $D_{0},\ldots,D_{10}$:
\begin{align}\label{GT}
D_{0}&\rightarrow\ii \bar D_{0},\,D_{1}\rightarrow -\bar
D_{1},\,D_{2}\rightarrow \bar D_{2},\,D_{3}\rightarrow\ii \bar
D_{3},\,D_{4}\rightarrow -\bar D_{4},\,D_{5}\rightarrow \bar
D_{5},\nonumber\\ D_{6}&\rightarrow \ii\bar
D_{6},\,D_{7}\rightarrow -\bar D_{7},\,D_{8}\rightarrow \ii\bar
D_{8},\,D_{9}\rightarrow
 -\bar D_{9},\,D_{10}\rightarrow -\bar D_{10}.
\end{align}
The operators $\bar D_{0},\ldots,\bar D_{10}$ generate the algebra
$\Diff(\mathbf{H}^{n}(\mathbb{H})_{\mathbb{S}})$.

\section{Relations in algebras $\Diff(\mathbf{P}^{n}(\mathbb{H})_{\mathbb{S}})$ and
$\Diff(\mathbf{H}^{n}(\mathbb{H})_{\mathbb{S}})$}
\label{Relations}\markright{\ref{Relations} Relations in
$\Diff(\mathbf{P}^{n}(\mathbb{H})_{\mathbb{S}})$ and
$\Diff(\mathbf{H}^{n}(\mathbb{H})_{\mathbb{S}})$}

Here we are to find the independent relations in
$\Diff(\mathbf{P}^{n}(\mathbb{H})_{\mathbb{S}})$ for it's
generators $D_{0},\ldots, \newline D_{10}$. They are of two types
(see section \ref{InvOperators}). First type is commutative
relations, because a commutator of two differential operator of
orders $m_{1}$ and $m_{2}$ is an operator of an order
$m_{3}\leqslant m_{1}+m_{2}-1$. It gives $11(11-1)/2=55$
relations. If $n\geqslant 3$ due to (\ref{ComRel}) the second type
consists of only one independent relation of the form:
\begin{equation}\label{NonComRel}
D_{10}^{2}-D_{4}D_{5}D_{9}-2D_{6}D_{7}D_{8}+D_{9}D_{6}^{2}+D_{4}D_{8}^{2}+D_{5}D_{7}^{2}=D',
\end{equation}
where $D'$ is an operator of an order $\leqslant 7$, which is
polynomial in $D_{0},\ldots, D_{10}$. If $n=2$ the formula
(\ref{EulerIdentity}) gives another independent relation of the
form:
\begin{equation*}
\frac12\{D_{1},D_{2}\}-D_{3}^{2}-D_{9}=D'',
\end{equation*}
where $D''$ is an operator of an order $\leqslant 3$, polynomial
in $D_{0},\ldots, D_{8}$. The direct calculations gives
$D''=D_{1}+D_{2}$, therefore in the case $n=2$ we have the
additional relations:
\begin{equation}\label{additional}
\frac12\{D_{1},D_{2}\}-D_{3}^{2}-D_{9}=D_{1}+D_{2}.
\end{equation}
For $n=2$ using this relation we can exclude the element $D_{9}$
from the list of generators.

In principle all relations can  be obtained by straightforward
calculations in $U(\mathfrak{g})$ modulo
$\left(U(\mathfrak{g})\mathfrak{k}\right)^{K_{0}}$, but these
calculations became too cumbersome to write all of them here. In
Appendix A there is an example of deriving some commutative
relation. After getting some commutative relations by direct
calculations it is possible to get some other ones (see Appendix
A) using the Jacobi identity:
$$
[D_{i},[D_{j},D_{k}]]+[D_{k},[D_{i},D_{j}]]+[D_{j},[D_{k},D_{i}]]=0,
$$
which is valid, in particular, in every associative algebra. This
identity gives also a tool for checking the commutative relations
already found. Below there are all 55 commutative relations in
lexicographic order. The relation (\ref{NonComRel}) became too
difficult to obtain in a similar way. May be we need a computer
algebra calculations to obtain the explicit expression for $D'$.

\begin{align}\label{CompComm}
[D_{0},D_{1}]&=-D_{3},\,[D_{0},D_{2}]=D_{3},\,[D_{0},D_{3}]=\frac12
(D_{1}-D_{2}),\,[D_{0},D_{4}]=-2D_{6},\,\nonumber\\
[D_{0},D_{5}]&=2D_{6},\,[D_{0},D_{6}]=D_{4}-D_{5},
[D_{0},D_{7}]=-D_{8}, [D_{0},D_{8}]=D_{7},\,[D_{0},D_{9}]=0,\,
\nonumber\\
[D_{0},D_{10}]&=0,\,[D_{1},D_{2}]=-\{D_{0},D_{3}\}-2D_{7},\,
[D_{1},D_{3}]=-\frac12\{D_{0},D_{1}\}+D_{8}+n(n-1)D_{0},\nonumber\\
[D_{1},D_{4}]&=2D_{7},\,[D_{1},D_{5}]=0,\,[D_{1},D_{6}]=D_{8},\,
[D_{1},D_{7}]=-\frac12\{D_{3},D_{6}\}-\frac12\{D_{1},D_{4}\}\nonumber\\
&+\frac38(D_{1}-D_{2})+D_{9}+D_{10}+n(n-1)D_{4},\,
[D_{1},D_{8}]=-\frac12\{D_{3},D_{5}\}-\frac12\{D_{1},D_{6}\}\nonumber\\
&+\frac34D_{3}+n(n-1)D_{6},\,
[D_{1},D_{9}]=-\{D_{3},D_{8}\}-\{D_{1},D_{7}\}-\frac34\{D_{0},D_{3}\}\nonumber\\
&+2(n-\frac32)(n+\frac12)D_{7},\,
[D_{1},D_{10}]=\frac12\{D_{6},D_{8}\}-\frac12\{D_{5},D_{7}\}+\frac38\{D_{0},D_{3}\}+\frac12D_{7},\nonumber\\
[D_{2},D_{3}]&=\frac12\{D_{0},D_{2}\}+D_{8}-n(n-1)D_{0},\,
[D_{2},D_{4}]=-2D_{7},\,[D_{2},D_{5}]=0,\nonumber\\
[D_{2},D_{6}]&=-D_{8},\,
[D_{2},D_{7}]=-\frac12\{D_{3},D_{6}\}+\frac12\{D_{2},D_{4}\}+\frac38(D_{1}-D_{2})-D_{9}-D_{10}\nonumber\\
&-n(n-1)D_{4},\,
[D_{2},D_{8}]=-\frac12\{D_{3},D_{5}\}+\frac12\{D_{2},D_{6}\}+\frac34
D_{3}-n(n-1)D_{6},\nonumber\\
[D_{2},D_{9}]&=-\{D_{3},D_{8}\}+\{D_{2},D_{7}\}+\frac34\{D_{0},D_{3}\}-2(n-\frac32)(n+\frac12)D_{7},
\\
[D_{2},D_{10}]&=-\frac12\{D_{6},D_{8}\}+\frac12\{D_{5},D_{7}\}-\frac38\{D_{0},D_{3}\}-\frac12D_{7},\,
[D_{3},D_{4}]=0,\,[D_{3},D_{5}]=2D_{8},\nonumber\\
[D_{3},D_{6}]&=D_{7},\,[D_{3},D_{7}]=-\frac14\{D_{1}+D_{2},D_{6}\}
+n(n-1)D_{6},\,[D_{3},D_{8}]=-\frac14\{D_{1}+D_{2},D_{5}\}\nonumber\\
&+n(n-1)D_{5}+D_{9}+D_{10},\,
[D_{3},D_{9}]=-\frac12\{D_{1}+D_{2},D_{8}\}+\frac38\{D_{0},D_{1}-D_{2}\}\nonumber\\
&+2(n-\frac32)(n+\frac12)D_{8},\,
[D_{3},D_{10}]=\frac12\{D_{6},D_{7}\}-\frac12\{D_{4},D_{8}\}-\frac3{16}\{D_{0},D_{1}-D_{2}\}\nonumber\\
&+\frac12D_{8},\,[D_{4},D_{5}]=-2\{D_{0},D_{6}\},\,[D_{4},D_{6}]=-\{D_{0},D_{4}\}+\frac32D_{0},\nonumber\\
[D_{4},D_{7}]&=\frac12\{D_{1}-D_{2},D_{4}\}+\frac34(D_{2}-D_{1}),\,
[D_{4},D_{8}]=\frac12\{D_{1}-D_{2},D_{6}\}-\{D_{0},D_{7}\},\nonumber\\
[D_{4},D_{9}]&=\{D_{1}-D_{2},D_{7}\},\,[D_{4},D_{10}]=0,\,
[D_{5},D_{6}]=\{D_{0},D_{5}\}-\frac32D_{0},\,\nonumber\\
[D_{5},D_{7}]&=\{D_{3},D_{6}\}+\{D_{0},D_{8}\},\,[D_{5},D_{8}]=\{D_{3},D_{5}\}-\frac32D_{3},\,
[D_{5},D_{9}]=2\{D_{3},D_{8}\},\,\nonumber\\ [D_{5},D_{10}]&=0,\,
[D_{6},D_{7}]=\frac14\{D_{1}-D_{2},D_{6}\}+\frac12\{D_{3},D_{4}\}+\frac12\{D_{0},D_{7}\}-\frac34D_{3},\nonumber\\
[D_{6},D_{8}]&=\frac14\{D_{1}-D_{2},D_{5}\}+\frac12\{D_{3},D_{6}\}-\frac12\{D_{0},D_{8}\}
+\frac38(D_{2}-D_{1}),\,\nonumber\\
[D_{6},D_{9}]&=\frac12\{D_{1}-D_{2},D_{8}\}+\{D_{3},D_{7}\},\,[D_{6},D_{10}]=0,\,
[D_{7},D_{8}]=\frac14\{D_{1}-D_{2},D_{8}\}\nonumber\\
&-\frac12\{D_{3},D_{7}\}+\frac3{16}\{D_{0},D_{1}+D_{2}\}
-\frac12\{D_{0},D_{9}+D_{10}\}-\frac34n(n-1)D_{0},\,\nonumber\\
[D_{7},D_{9}]&=\frac14\{D_{3},D_{6}\}+\frac18\{D_{1}-D_{2},D_{4}\}+\frac12\{D_{1}-D_{2},D_{9}+D_{10}\}-
\frac38(D_{1}^{2}-D_{2}^{2})\nonumber\\
&+\frac34(n^{2}-n-\frac14)(D_{1}-D_{2}),\,
[D_{7},D_{10}]=\frac14\{D_{2}-D_{1},D_{6}^{2}\}-\frac14\{\{D_{0},D_{7}\},D_{6}\}\nonumber\\
&+\frac14\{\{D_{0},D_{4}\},D_{8}\}+\frac18\{\{D_{1}-D_{2},D_{5}\},D_{4}\}-\frac14\{D_{3},D_{6}\}\nonumber\\
&+\frac18\{D_{2}-D_{1},3D_{4}+D_{5}\}-\frac12\{D_{0},D_{8}\}+\frac{15}{32}(D_{1}-D_{2}),\nonumber\\
[D_{8},D_{9}]&=\frac18\{D_{1}-D_{2},D_{6}\}+\frac14\{D_{3},D_{5}\}-\frac38\{D_{3},D_{1}+D_{2}\}
+\{D_{3},D_{9}+D_{10}\}\nonumber\\
&+\frac32(n^{2}-n-\frac14)D_{3},\,[D_{8},D_{10}]=-\frac14\{\{D_{3},D_{6}\},D_{6}\}
+\frac14\{\{D_{0},D_{6}\},D_{8}\}\nonumber\\
&-\frac14\{\{D_{0},D_{5}\},D_{7}\}+\frac14\{\{D_{3},D_{5}\},D_{4}\}-
\frac12\{D_{3},D_{5}\}-\frac14\{D_{3},D_{4}\}+\frac14\{D_{0},D_{7}\}\nonumber\\
&+\frac{9}{16}D_{3},\,[D_{9},D_{10}]=\frac14\{-\{D_{6},D_{8}\}+\{D_{5},D_{7}\},D_{1}-D_{2}\}+
\frac12\{\{D_{3},D_{8}\},D_{4}\}\nonumber\\
&-\frac12\{\{D_{3},D_{6}\},D_{7}\}+\frac14\{D_{2}-D_{1},D_{7}\}-
\frac12\{D_{3},D_{8}\}.\nonumber
\end{align}

It is interesting that the operators $D_{9}$ and $D_{10}$ arise in
the right hand sides of these relations only in the combination
$D_{9}+D_{10}$.

Using relations (\ref{CompComm}) it is not difficult to verify
that the operator $D^{*}=D_{0}^{2}+D_{1}+D_{2}+D_{4}+D_{5}$ lies
in the centre of the algebra
$\Diff(\mathbf{P}^{n}(\mathbb{H})_{\mathbb{S}})$ in accordance
with the section \ref{spes}.

Using substitution (\ref{GT}) one can obtains from
(\ref{CompComm}) the commutative relations for the algebra
$\Diff(\mathbf{H}^{n}(\mathbb{H})_{\mathbb{S}})$.

The analog for the operator $D^{*}$ from the centre of the algebra
now becomes $\bar D^{*}=\bar D_{0}^{2}+\bar D_{1}-\bar D_{2}+\bar
D_{4}-\bar
D_{5}\in\Diff(\mathbf{H}^{n}(\mathbb{H})_{\mathbb{S}})$. In the
case $n=2$ the additional relations (\ref{additional}) becomes:
\begin{equation*}
\frac12\{\bar D_{1},\bar D_{2}\}-\bar D_{3}^{2}-\bar D_{9}=\bar
D_{1}-\bar D_{2}.
\end{equation*}

\section{Algebras $\Diff(\mathbf{P}^{n}(\mathbb{C})_{\mathbb{S}})$ and
$\Diff(\mathbf{H}^{n}(\mathbb{C})_{\mathbb{S}})$}
\label{Complex}\markright{\ref{Complex} Algebras
$\Diff(\mathbf{P}^{n}(\mathbb{C})_{\mathbb{S}})$ and
$\Diff(\mathbf{H}^{n}(\mathbb{C})_{\mathbb{S}})$}

\subsection{The model for the space $\mathbf{P}^{n}(\mathbb{C})$}

Taking the factor space of $\mathbb{C}^{n+1}\backslash\{0\}$ w.r.t
the action of the multiplicative group
$\mathbb{C}^{*}=\mathbb{C}\backslash\{0\}$ (due to the
commutativity of the complex multiplication it makes no difference
left or right), we obtain {\it the complex projective space}
$\mathbf{P}^{n}(\mathbb{C})$. Let $\pi :
\mathbb{C}^{n+1}\backslash\{0\}\rightarrow
\mathbf{P}^{n}(\mathbb{C})$ be the canonical projection. Let now
$\langle\x,\y\rangle:=\sum\limits_{i=1}^{n+1}\bar x_{i}
y_{i},\,\x=(x_{1},\ldots,x_{n+1}),\,\y=(y_{1},\ldots,y_{n+1})\in
\mathbb{C}^{n+1}$ be the standard scalar product in the space
$\mathbb{C}^{n+1}$.

The metric $\tilde g$ of the constant holomorphic sectional
curvature on the space $\mathbf{P}^{n}(\mathbb{C})$ is defined by
the same formula (\ref{metricPH}) as on the space
$\mathbf{P}^{n}(\mathbb{H})$, where now
$\z\in\mathbb{C}^{n+1}\backslash\{0\},\,{\bf\xi}_i\in
T_{\z}\mathbb{C}^{n+1},\,{\bf \zeta}_{i}=\pi_{*}{\bf\xi}_{i}\in
T_{\pi(\z)}\left(\mathbf{P}^{n}(\mathbb{C})\right),\,i=1,2$.

The Riemannian metric $g$ on the space
$\mathbf{P}^{n}(\mathbb{C})$ is:
\begin{equation}\label{metricPRC}
g=4R^{2}\RE\tilde g.
\end{equation}
If $n=2$ it is not difficult to verify (like in section
\ref{QuatMod}) that the homeomorphism
$\tau:\mathbf{P}^{1}(\mathbb{C})\rightarrow
\bar{\mathbb{C}}\cong{\bf S}^{2},\,
\tau(z_{1},z_{2})=z_{1}\left(z_{2}\right)^{-1}=z\in\overline{\mathbb{C}}$,
transforms (\ref{metricPRC}) into the metric
$$
g=\frac{4R^{2}dzd\bar z}{\left(1+|z|^{2}\right)^{2}},
$$
of the sectional curvature $R^{-2}$ on the sphere $\mathbf{
S}^{2}$.

The left action of the group $G=\SU(n+1)$ on the space
$\mathbb{C}^{n+1}$ conserves the scalar product
$\langle\cdot,\cdot\rangle$ and induces the action in the space $
\mathbf{P}^{n}(\mathbb{C})$, conserving metrics $\tilde g$ and
$g$.

The stationary subgroup, corresponding to the point of the space
$\mathbf{P}^{n}(\mathbb{C})$ with homogeneous coordinates
$(1:0:\ldots:0)$, are the group $\U(n)=\SU(n)\U(1)$, where the
factor $\SU(n)$ acts in the standard way onto the last $n$
coordinates, and the factor $\U(1)$ acts by the multiplication of
the first coordinate by $e^{\ii\phi}$ and the second one by
$e^{-\ii\phi},\,\phi\in R\mod 2\pi$. Thus
$\mathbf{P}^{n}(\mathbb{C})=\SU(n+1)/\U(n)$.

Choose a base of the algebra $\mathfrak{su}(n+1)$ in the form:
\begin{align}
\label{BasisSU1} \Psi_{kj}=\frac12(E_{kj}-E_{jk}),\,
\Upsilon_{kj}=\frac{\ii}2(E_{kj}+E_{jk}),1\leqslant k<j\leqslant
n+1,
\end{align}
\begin{align}
\label{BasisSU2}
\Upsilon_{k}=\frac{\ii}2(E_{11}-E_{kk}),\,2\leqslant k\leqslant
n+1.
\end{align}
The commutators for these elements are easily extracted from
(\ref{commutators1}), taking into account that
$\Upsilon_{k}=\frac12(\Upsilon_{11}-\Upsilon_{kk})$ using the
notations from (\ref{BasisUH}).

\subsection{Algebras
$\Diff(\mathbf{P}^{n}(\mathbb{C})_{\mathbb{S}})$ and
$\Diff(\mathbf{H}^{n}(\mathbb{C})_{\mathbb{S}})$}

Consider now the space $\mathbf{P}^{n}(\mathbb{C})_{\mathbb{S}}$.
Due to the isomorphism $\mathbf{P}^{1}(\mathbb{C})\cong{\bf
S}^{2}$ we again assume that $n\geqslant 2$.

Let $\tilde \z_{0}=(1,0,\ldots,0)\in \mathbb{C}^{n+1}$, an element
$\xi_{0}\in T_{\tilde\z_{0}}\mathbb{C}^{n+1}\cong
\mathbb{C}^{n+1}$ has coordinates $(0,1,0,\ldots,0)$. Put
$\z_{0}=\pi\tilde\z_{0},\,\zeta_{0}=\pi_{*}\xi_{0}\in
T_{\z_{0}}\mathbf{P}^{n}(\mathbb{C})$.

The stationary subgroup $K_{0}$ of the group $\SU(n+1)$,
corresponding to the point $(\z_{0},\zeta_{0})$, is generated by
the group $K_{1}=\SU(n-1)$, acting onto the last $n-1$-th
coordinates and by the group $K_{2}=\U(1)$, acting onto the
homogeneous coordinates of $\mathbf{P}^{n}(\mathbb{C})$ as:
\begin{equation}\label{actK2}
(x_{1}:\ldots:x_{n+1})\rightarrow
(e^{\ii\phi}x_{1}:e^{\ii\phi}x_{2}:e^{-2\ii\phi}x_{3}:x_{4}:\ldots:x_{n+1}),
\end{equation}
$\dim_{\mathbb{R}}K_{0}=(n-1)^{2}$ and we obtain
$K_{0}\cong\U(n-1)$.

The algebra $\mathfrak{k}_{0}$ of the group $K_{0}$ is the linear
hull of elements (\ref{BasisSU1}) as $3\leqslant k<j\leqslant n+1$
and elements:
$$
\Upsilon_{j}-\Upsilon_{3}=\frac{\ii}2(E_{33}-E_{jj}),\,3<j\leqslant
n+1,\,2\Upsilon_{3}-\Upsilon_{2}=\frac{\ii}2(E_{11}+E_{22}-2E_{33}).
$$

Choose the complimentary subspace $\tilde{\mathfrak{p}}$ to the
subalgebra $\mathfrak{k}_{0}$ in the algebra
$\mathfrak{g}=\mathfrak{su}(n+1)$ as the linear hull of elements:
\begin{align}\label{subspacePC}
\Psi_{1k},\,\Upsilon_{1k},\,2\leqslant k\leqslant n+1,\,
\Psi_{2k},\,\Upsilon_{2k},\,3\leqslant k\leqslant
n+1,\,\Upsilon_{*}=\Upsilon_{2}.
\end{align}
Taking into account relations (\ref{commutators1}) it is easily
obtained that the expansion
$\mathfrak{su}(n+1)=\tilde{\mathfrak{p}}\oplus\mathfrak{k}_{0}$ is
reductive, i.e.
$[\tilde{\mathfrak{p}},\mathfrak{k}_{0}]\subset\tilde{\mathfrak{p}}$.

We will obtain the particular case of Proposition \ref{prop1} for
$q_{1}=2n-2,q_{2}=1$  setting:
\begin{align}\label{reper2}
\Lambda=-\Psi_{12},\,e_{\lambda,k-2}&=\Psi_{1k},\,e_{\lambda,n-3+k}=\Upsilon_{1k},\,
f_{\lambda,k-2}=-\Psi_{2k},\,f_{\lambda,n-3+k}=-\Upsilon_{2k},\,\notag
\\ e_{2\lambda,1}&=\Upsilon_{12},\,
f_{2\lambda,1}=\Upsilon_{*},\, k=3,\dots,n+1.
\end{align}

Now we are to find the generators of $\Ad_{K_{0}}$-invariants in
$S(\tilde{\mathfrak{p}})$. The expansion
$\tilde{\mathfrak{p}}=\mathfrak{a}\oplus\mathfrak{k}_{\lambda}\oplus
\mathfrak{k}_{2\lambda}\oplus\mathfrak{p}_{\lambda}\oplus
\mathfrak{p}_{2\lambda}$ is invariant w.r.t. the
$\Ad_{K_{0}}$-action. In the spaces
$\mathfrak{a},\,\mathfrak{p}_{2\lambda},\,\mathfrak{k}_{2\lambda}$
the $K_{0}$-action is trivial that gives the invariants
$D_{0}=\Lambda,\,D_{4}=\Upsilon_{12},\,D_{5}=\Upsilon_{*}\in\mu(\mathfrak{p}^{K_{0}})$.
Operators $D_{4},D_{5}$ are square roots of their analogs from
section \ref{spes}.

From formulas (\ref{reper2}) we see that the space
$\mathfrak{p}_{\lambda}\cong \mathbb{C}^{n-1}$ consists of
matrices of the form
$$
\left(\begin{array}{cc} 0 & - a^{*} \\ a & 0
\end{array}\right)\equiv\left(\begin{array}{ccccc} 0 & 0 & -\bar a_{1} & \ldots & -\bar
a_{n-1} \\ 0 & 0 & 0 & \ldots & 0 \\ a_{1} & 0 & 0 & \ldots & 0 \\
\vdots & \vdots & \vdots & \ddots & \vdots \\ a_{n-1} & 0 & 0 &
\ldots & 0
\end{array}\right),\,a_{1},\ldots,a_{n-1}\in\mathbb{C}.
$$
Similarly, the space $\mathfrak{k}_{\lambda}\cong
\mathbb{C}^{n-1}$ consist of matrices of the form
$$
\left(\begin{array}{ccc} 0 & 0 & 0 \\ 0 & 0 & - b^{*} \\ 0 & b & 0
\end{array}\right)\equiv\left(\begin{array}{ccccc} 0 & 0 & 0 & \ldots & 0 \\ 0 & 0 & -\bar
b_{1} & \ldots & -\bar b_{n-1} \\ 0 & b_{1} & 0 & \ldots & 0
\\ \vdots & \vdots & \vdots & \ddots & \vdots \\ 0 & b_{n-1} & 0 &
\ldots & 0
\end{array}\right),\,b_{1},\ldots,b_{n-1}\in\mathbb{C}.
$$
The action of the group $K_{1}$ in the spaces
$\mathfrak{p}_{\lambda}$ and $\mathfrak{k}_{\lambda}$ is
equivalent to the standard action of the group $\SU(n-1)$ in the
space $\mathbb{C}^{n-1}:\, a\rightarrow Ua,\,U\in \SU(n-1)$,
likewise in section \ref{Quat}. It is easy to verify that the
action (\ref{actK2}) generates the action $a_{1}\rightarrow
\exp^{-3\ii\phi}a_{1},\,a_{i}\rightarrow
\exp^{-\ii\phi}a_{i},b_{1}\rightarrow
\exp^{-3\ii\phi}b_{1},\,b_{i}\rightarrow
\exp^{-\ii\phi}b_{i},i=2,\ldots,n-1$. Therefore the $K_{0}$-action
in spaces $\mathfrak{p}_{\lambda}$ and $\mathfrak{k}_{\lambda}$ is
equivalent to the standard $\U(n-1)$-action in $\mathbb{C}^{n-1}$.

This action has one independent real invariant:
$\langle\z,\z\rangle,\, \z\in \mathbb{C}^{n-1}$, and the diagonal
action of $\U(n-1)$ in the space
$\mathfrak{p}_{\lambda}\oplus\mathfrak{k}_{\lambda}\cong
\mathbb{C}^{n-1}\oplus\mathbb{C}^{n-1}$ has four (independent iff
$n\geqslant 3$) real invariants:
\begin{equation}\label{ComInvC}
\langle\z_{1},\z_{1}\rangle\in
\mathbb{R},\,\langle\z_{2},\z_{2}\rangle\in
\mathbb{R},\,\langle\z_{1},\z_{2}\rangle\in
\mathbb{C}\cong\mathbb{R}^{2},\,\z_{1},\z_{2}\in\mathbb{C}^{n-1}.
\end{equation}
Denote the corresponding elements from
$\mu(\tilde{\mathfrak{p}}^{K_{0}})\in U(\mathfrak{g})^{K_{0}}$ in
the following way:
\begin{align}\label{def1C}
D_{1}&=\sum\limits_{k=3}^{n+1}\left(\Psi_{1k}^{2}+\Upsilon_{1k}^{2}\right),\,
D_{2}=\sum\limits_{k=3}^{n+1}\left(\Psi_{2k}^{2}+\Upsilon_{2k}^{2}\right),\\
D_{3}&=-\frac12\sum\limits_{k=3}^{n+1}\left(\left\{\Psi_{1k},\Psi_{2k}\right\}+\left\{\Upsilon_{1k},\Upsilon_{2k}\right\}
\right),\,
\square=\frac12\sum\limits_{k=3}^{n+1}\left(-\left\{\Psi_{1k},\Upsilon_{2k}\right\}+
\left\{\Psi_{2k},\Upsilon_{1k}\right\}\right)\nonumber.
\end{align}
In this case only operator $\square$ is new w.r.t. section
\ref{spes}.

If $n=2$, then there is the unique independent relation between
invariants (\ref{ComInvC}):
\begin{equation}\label{AddInv}
|\langle \z_{1},\z_{2}\rangle|^{2}=|\bar
z_{1}z_{2}|^{2}=|z_{1}|^{2}\,|z_{2}|^{2}=\langle
\z_{1},\z_{1}\rangle\langle
\z_{2},\z_{2}\rangle,\,\z_{1}=z_{1},\z_{2}=z_{2}\in \mathbb{C}.
\end{equation}

Thus operators $D_{0},\ldots,D_{5},\square$ generate the algebra
$\Diff(\mathbf{P}^{n}(\mathbb{C})_{\mathbb{S}})$.

The degrees of the generators are as follows:
\begin{align}\label{degreesC}
\degr(D_{0})&=\degr(D_{4})=\degr(D_{5})=1,\,\degr(D_{1})=\degr(D_{2})=\degr(D_{3})=\degr(\square)=2.
\end{align}

The operators $D_{3},D_{4}$ are symmetric and the operators
$D_{0},\square,D_{5}$ are skew symmetric w.r.t. the transposition
of coordinates $z_{1}$ and $z_{2}$. The operators $D_{1}$ and
$D_{2}$ turn into each other under this transposition.

In order to get the generators of the algebra
$\Diff(\mathbf{H}^{n}(\mathbb{C})_{\mathbb{S}})$ we can use the
formal substitution:
\begin{align*}
\Lambda &\rightarrow
\ii\Lambda,\,\Psi_{1k}\rightarrow\ii\Psi_{1k},\,\Upsilon_{1k}\rightarrow\ii\Upsilon_{1k},\,
\Upsilon_{12}\rightarrow\ii\Upsilon_{12},\\
\Psi_{2k}&\rightarrow\Psi_{2k},\,\Upsilon_{2k}\rightarrow\Upsilon_{2k},\,
\Upsilon_{*}\rightarrow\Upsilon_{*},\, \,k=3,\ldots,n+1.
\end{align*}
This substitution produces the following substitution for the
generators $D_{0},\ldots,D_{5},\square$:
\begin{align}\label{GTC}
D_{0}&\rightarrow\ii \bar D_{0},\,D_{1}\rightarrow -\bar
D_{1},\,D_{2}\rightarrow \bar D_{2},\,D_{3}\rightarrow\ii \bar
D_{3},\,D_{4}\rightarrow \ii\bar D_{4},\,\square\rightarrow\ii
\bar \square,\, D_{5}\rightarrow \bar D_{5}.
\end{align}
The operators $\bar D_{0},\ldots,\bar D_{5},\bar\square$ generate
the algebra $\Diff(\mathbf{H}^{n}(\mathbb{C})_{\mathbb{S}})$.

\subsection{Relations in algebras
$\Diff(\mathbf{P}^{n}(\mathbb{C})_{\mathbb{S}})$ and
$\Diff(\mathbf{H}^{n}(\mathbb{C})_{\mathbb{S}})$}

The commutative relation for the algebra
$\Diff(\mathbf{P}^{n}(\mathbb{C})_{\mathbb{S}})$ are as follows:
\begin{align*}
&[D_{0},D_{1}]=-D_{3},\,[D_{0},D_{2}]=D_{3},\,[D_{0},D_{3}]=\frac12(D_{1}-D_{2}),\,[D_{0},D_{4}]=-D_{5},\,
[D_{0},D_{5}]=D_{4},\nonumber\\
&[D_{0},\square]=0,\,[D_{1},D_{2}]=-\{D_{0},D_{3}\}-\{\square,D_{4}\},\,[D_{1},D_{3}]=-\frac12\{D_{0},D_{1}\}
+\frac12\{\square,D_{5}\}\nonumber\\&+\frac{(n-1)^{2}}4D_{0},\,[D_{1},D_{4}]=\square,\,[D_{1},D_{5}]=0,\,
[D_{1},\square]=-\frac12\{D_{1},D_{4}\}-\frac12\{D_{3},D_{5}\}\nonumber\\&+\frac{(n-1)^{2}}4D_{4},\,
[D_{2},D_{3}]=\frac12\{D_{0},D_{2}\}+\frac12\{\square,D_{5}\}-\frac{(n-1)^{2}}4D_{0},\,
[D_{1},D_{4}]=\square,\\&[D_{2},D_{5}]=0,\,
[D_{2},\square]=\frac12\{D_{2},D_{4}\}-\frac12\{D_{3},D_{5}\}-\frac{(n-1)^{2}}4D_{4},\,[D_{3},D_{4}]=0,\\
&[D_{3},D_{5}]=\square,\,
[D_{3},\square]=-\frac14\{D_{1}+D_{2},D_{5}\}+\frac{(n-1)^{2}}4D_{5},\,[D_{4},D_{5}]=-D_{0},\\
&[D_{4},\square]=\frac12(D_{1}-D_{2}),\,[D_{5},\square]=D_{3}.
\end{align*}

If $n>2$ then there are no relations of the second type. If $n=2$
then there is one relation of the second type due to
(\ref{AddInv}):
\begin{equation}
\label{AddInvNC}
\frac12\{D_{1},D_{2}\}-D_{3}^{2}-\square^{2}-\frac14(D_{0}^{2}+D_{4}^{2}+D_{5}^{2})=0.
\end{equation}
It is easy to verify that the operator
$D^{*}=D_{0}^{2}+D_{1}+D_{2}+D_{4}^{2}+D_{5}^{2}$ lies in the
centre of the algebra
$\Diff(\mathbf{P}^{n}(\mathbb{C})_{\mathbb{S}})$ in accordance
with the section \ref{spes}.

Using substitution (\ref{GTC}) we obtain analogous relations for
the algebra $\Diff(\mathbf{H}^{n}(\mathbb{C})_{\mathbb{S}})$.

The commutative relation are now as follows:
\begin{align*}
&[\D_{0},\D_{1}]=\D_{3},\,[\D_{0},\D_{2}]=\D_{3},\,[\D_{0},\D_{3}]=\frac12(\D_{2}+\D_{1}),\,[\D_{0},\D_{4}]=\D_{5},\,
[\D_{0},\D_{5}]=\D_{4},\nonumber\\
&[\D_{0},\bar\square]=0,\,[\D_{1},\D_{2}]=-\{\D_{0},\D_{3}\}-\{\bar\square,\D_{4}\},\,[\D_{1},\D_{3}]=-\frac12\{\D_{0},
\D_{1}\}-\frac12\{\bar\square,\D_{5}\}\nonumber\\&-\frac{(n-1)^{2}}4\D_{0},\,[\D_{1},\D_{4}]=-\bar\square,\,
[\D_{1},\D_{5}]=0,\,[\D_{1},\bar\square]=-\frac12\{\D_{1},\D_{4}\}+\frac12\{\D_{3},\D_{5}\}\nonumber\\&-
\frac{(n-1)^{2}}4\D_{4},\,
[\D_{2},\D_{3}]=\frac12\{\D_{0},\D_{2}\}+\frac12\{\bar\square,\D_{5}\}-\frac{(n-1)^{2}}4\D_{0},\,
[\D_{1},\D_{4}]=-\bar\square,\\&[\D_{2},\D_{5}]=0,\,
[\D_{2},\bar\square]=\frac12\{\D_{2},\D_{4}\}-\frac12\{\D_{3},\D_{5}\}-\frac{(n-1)^{2}}4\D_{4},\,[\D_{3},\D_{4}]=0,\\
&[\D_{3},\D_{5}]=\bar\square,\,
[\D_{3},\bar\square]=-\frac14\{\D_{1}-\D_{2},\D_{5}\}-\frac{(n-1)^{2}}4\D_{5},\,[\D_{4},\D_{5}]=-\D_{0},\\
&[\D_{4},\bar\square]=\frac12(\D_{1}+\D_{2}),\,[\D_{5},\bar\square]=\D_{3}.
\end{align*}

If $n>2$ then there are no relations of the second type. If $n=2$
then there is one relation of the second type analogous to
(\ref{AddInvNC}):
\begin{equation}
\label{AddInvNCH}
\frac12\{\D_{1},\D_{2}\}-\D_{3}^{2}-\bar\square^{2}-\frac14(\D_{0}^{2}+\D_{4}^{2}-\D_{5}^{2})=0.
\end{equation}
The operator
$\D^{*}=\D_{0}^{2}+\D_{1}-\D_{2}+\D_{4}^{2}-\D_{5}^{2}$ lies in
the centre of the algebra
$\Diff(\mathbf{H}^{n}(\mathbb{C})_{\mathbb{S}})$.

\section{Algebras $\Diff(\mathbf{P}^{n}(\mathbb{R})_{\mathbb{S}}),\,
\Diff(\mathbf{S}^{n}_{\mathbb{S}})$ and
$\Diff(\mathbf{H}^{n}(\mathbb{R})_{\mathbb{S}})$}
\label{Real}\markright{\ref{Real} Algebras
$\Diff(\mathbf{P}^{n}(\mathbb{R})_{\mathbb{S}}),\,
\Diff(\mathbf{S}^{n}_{\mathbb{S}})$ and
$\Diff(\mathbf{H}^{n}(\mathbb{R})_{\mathbb{S}})$}

\subsection{Generators of algebras
$\Diff({\bf S}^{n}_{\mathbb{S}})$ and
$\Diff(\mathbf{H}^{n}(\mathbb{R})_{\mathbb{S}})$}

Let now $\langle\cdot,\cdot\rangle$ be the standard scalar product
in the space $\mathbb{R}^{n+1}$. The equation \\
$\langle\x,\x\rangle=R>0$ defines the sphere ${\bf
S}^{n}\cong\SO(n+1)/\SO(n)\subset \mathbb{R}^{n+1}$ of the radius
$R$ with the induced metric on it. The space
$\mathbf{P}^{n}(\mathbb{R})$ is the factor space of ${\bf S}^{n}$
w.r.t. the relation: $\x\sim -\x$. Below we will show that
algebras $\Diff(\mathbf{P}^{n}(\mathbb{R})_{\mathbb{S}})$ and
$\Diff({\bf S}^{n}_{{\mathbb{S}}})$ are isomorphic.

The spaces ${\bf
S}^{1}_{\mathbb{S}},\,\mathbf{P}^{1}(\mathbb{R})_{\mathbb{S}}$ are
one dimensional and the algebra of invariant differential
operators on them is generated by one differential operator of the
first order. Therefore we again suppose that $n\geq2$.

Let
\begin{align}
\label{BasisSO} \Psi_{kj}=\frac12(E_{kj}-E_{jk}),\,1\leqslant
k<j\leqslant n+1
\end{align}
be the base of the algebra $\mathfrak{so}(n+1)$. The commutative
relations for them are contained in (\ref{commutators1}).

Consider the space ${\bf S}^{n}_{\mathbb{S}}$. Let
$\tilde\z_{0}=(1,0,\ldots,0)\in \mathbb{R}^{n+1}$, an element
$\xi_{0}\in T_{\tilde\z_{0}}\mathbb{R}^{n+1}\cong
\mathbb{R}^{n+1}$ has coordinates $(0,1,0,\ldots,0)$. Put
$\z_{0}=\pi\tilde\z_{0},\,\zeta_{0}=\pi_{*}\xi_{0}\in
T_{\z_{0}}{\bf S}^{n}_{\mathbb{S}}$.

The stationary subgroup $K_{0}$ of the group $\SO(n+1)$,
corresponding to the point $(\z_{0},\zeta_{0})\in{\bf
S}^{n}_{\mathbb{S}}$, is the group $\SO(n-1)$, acting onto the
last $(n-1)$-th coordinates.

The group $\SO(n+1)$ is a group covering of the identity component
$G$ of the isometry group for $\mathbf{P}^{n}(\mathbb{R})$. The
group $K_{0}=\SO(n-1)\subset \SO(n+1)$ is a group covering of its
analog $K_{0}'\subset G$. The kernel of such covering is a normal
subgroup of $\SO(n-1)$ that lies in the centre of $\SO(n-1)$
\cite{Post} (lecture 9). Therefore the orbits of $\Ad_{K_{0}}$ and
$\Ad_{K_{0}'}$ actions on $\mathfrak{p}\subset\mathfrak{g}$
coincide and the construction from section \ref{InvOperators}
implies the isomorphism
$\Diff(\mathbf{P}^{n}(\mathbb{R})_{\mathbb{S}})\cong\Diff({\bf
S}^{n}_{{\mathbb{S}}})$.

The algebra $\mathfrak{k}_{0}$ of the group $K_{0}$ is the linear
hull of elements $\Psi_{kj}$ as $3\leqslant k<j\leqslant n+1$.
Choose the complimentary subspace $\tilde{\mathfrak{p}}$ to the
subalgebra $\mathfrak{k}_{0}$ in the algebra
$\mathfrak{g}=\mathfrak{so}(n+1)$ as the linear hull of elements:
\begin{align}\label{subspaceS}
\Psi_{1k},\,2\leqslant k\leqslant n+1,\, \Psi_{2k},\,3\leqslant
k\leqslant n+1.
\end{align}
The expansion $\mathfrak{so}(n+1)=\tilde{\mathfrak{p}}\oplus
\mathfrak{so}(n-1)$ is reductive.

We will obtain the particular case of Proposition \ref{prop1} for
$q_{1}=0,q_{2}=n-1$ setting:
\begin{align}\label{reper3}
\Lambda=-2\Psi_{12},\,e_{2\lambda,k-2}=2\Psi_{1k},\,
f_{2\lambda,k-2}=-2\Psi_{2k},\, k=3,\dots,n+1.
\end{align}

Now we have the expansion $\tilde{\mathfrak{p}}=\mathfrak{a}\oplus
\mathfrak{k}_{2\lambda}\oplus \mathfrak{p}_{2\lambda}$, which is
invariant w.r.t. the $\Ad_{K_{0}}$-action. It is easy to see that
in the space $\mathfrak{a}$ the $K_{0}$-action is trivial and in
the spaces $\mathfrak{k}_{2\lambda}$ and $\mathfrak{p}_{2\lambda}$
it is equivalent to the standard action of the group $\SO(n-1)$ in
the space $\mathbb{R}^{n-1}$. The $K_{0}$-action in the space
$\mathfrak{a}$ has the invariant $D_{0}=\Lambda$. The description
of base $K_{0}$-invariants in the space
$\mathfrak{p}_{2\lambda}\oplus\mathfrak{k}_{2\lambda}$ is
different in cases $n=2,\,n=3$ and $n\geqslant 4$.

\subsubsection{The case $n\geqslant 4$}
The $\SO(n-1)$-action in $\mathbb{R}^{n-1}$ has one independent
real invariant: $\langle\z,\z\rangle,\, \z\in \mathbb{R}^{n-1}$,
and the diagonal action of $\SO(n-1)$ in the space
$\mathfrak{p}_{2\lambda}\oplus\mathfrak{k}_{2\lambda}\cong
\mathbb{R}^{n-1}\oplus\mathbb{R}^{n-1}$ has three independent real
invariants:
\begin{equation}\label{ComInvR}
\langle\z_{1},\z_{1}\rangle,\,\langle\z_{2},\z_{2}\rangle,\,\langle\z_{1},\z_{2}\rangle,
\,\z_{1},\z_{2}\in\mathbb{R}^{n-1}.
\end{equation}

Denote the corresponding elements from
$\mu(\tilde{\mathfrak{p}}^{K_{0}})\in U(\mathfrak{g})^{K_{0}}$ in
the following way:
\begin{align}
D_{1}=4\sum\limits_{k=3}^{n+1}\Psi_{1k}^{2},\,
D_{2}=4\sum\limits_{k=3}^{n+1}\Psi_{2k}^{2},\,
D_{3}=-2\sum\limits_{k=3}^{n+1}\left\{\Psi_{1k},\Psi_{2k}\right\}
\nonumber.
\end{align}
All these invariants were found in section \ref{spes} for the
general situation.

Thus operators $D_{0},D_{1},D_{2},D_{3}$ generate the algebra
$\Diff({\bf S}^{n}_{\mathbb{S}})$.

The degrees of the generators are as follows:
\begin{align}\label{degreesR}
\degr(D_{0})=1,\,\degr(D_{1})=\degr(D_{2})=\degr(D_{3})=2.
\end{align}

The operator $D_{3}$ is symmetric and the operators $D_{0}$ is
skew symmetric w.r.t. the transposition of coordinates $z_{1}$ and
$z_{2}$. The operators $D_{1}$ and $D_{2}$ turn into each other
under this transposition.

\subsubsection{The case $n=2$}
In this case $K_{0}$ is the trivial group and the independent
invariants are
$D_{0},\,D_{1}=e_{2\lambda,1},\,D_{2}=f_{2\lambda,1}$. Thus the
algebra $\Diff({\bf S}^{2}_{\mathbb{S}})$ is isomorphic to
$U(\mathfrak{so}(3))$. The centre of this algebra is generated by
the operator $D_{0}^{2}+D_{1}^{2}+D_{2}^{2}$.

\subsubsection{The case $n=3$}
In this case $K_{0}=\mathfrak{so}(2)$ and we have the additional
(with respect to the case $n\geqslant 4$) invariant of the second
order
$$
\square=2(\{\Psi_{13},\Psi_{24}\}-\{\Psi_{14},\Psi_{23}\}).
$$
It is algebraically connected with operators
$D_{0},D_{1},D_{2},D_{3}$ which are defined as in the case
$n\geqslant 4$.

\subsubsection{Generators of the algebra $\Diff(\mathbf{H}^{n}(\mathbb{R})_{\mathbb{S}})$}

First, let $n\geqslant 4$. In order to get the generators of the
algebra $\Diff(\mathbf{H}^{n}(\mathbb{R})_{\mathbb{S}})$ we can
use the formal substitution:
\begin{align*}
\Lambda \rightarrow
\ii\Lambda,\,\Psi_{1k}\rightarrow\ii\Psi_{1k},\,
\Psi_{2k}\rightarrow\Psi_{2k},\,k=3,\ldots,n+1.
\end{align*}
This substitution produces the following substitution for the
generators $D_{0},\ldots,D_{3}$:
\begin{align}\label{GTR}
D_{0}&\rightarrow\ii \bar D_{0},\,D_{1}\rightarrow -\bar
D_{1},\,D_{2}\rightarrow \bar D_{2},\,D_{3}\rightarrow\ii \bar
D_{3}.
\end{align}
The operators $\bar D_{0},\ldots,\bar D_{3}$ generate the algebra
$\Diff(\mathbf{H}^{n}(\mathbb{R})_{\mathbb{S}})$.

In the case $n=3$ we have the additional substitution
$\square\rightarrow\ii\bar\square$ and the operators $\bar
D_{0},\ldots,\bar D_{3},\bar\square$ generate the algebra
$\Diff(\mathbf{H}^{3}(\mathbb{R})_{\mathbb{S}})$.

In the case $n=2$ we obtain the substitution
\begin{align*}
D_{0}&\rightarrow\ii \bar D_{0},\,D_{1}\rightarrow \ii\bar
D_{1},\,D_{2}\rightarrow \bar D_{2}.
\end{align*}
The algebra $\Diff(\mathbf{H}^{2}(\mathbb{R})_{\mathbb{S}})$ is
isomorphic to $U(\mathfrak{so}(2,1))$ and its centre is generated
by the operator $D_{0}^{2}+D_{1}^{2}-D_{2}^{2}$.

\subsection{Relations in algebras
$\Diff({\bf S}^{n}_{\mathbb{S}})$ and
$\Diff(\mathbf{H}^{n}(\mathbb{R})_{\mathbb{S}})$}

Here we shall consider only the case $n\geqslant 3$, since
$\Diff({\bf S}^{2}_{\mathbb{S}})\cong U(\mathfrak{so}(3))$ and
$\Diff(\mathbf{H}^{2}(\mathbb{R})_{\mathbb{S}})\cong
U(\mathfrak{so}(2,1))$.

The commutative relation for the algebra $\Diff({\bf
S}^{n}_{\mathbb{S}})$ are as follows:
\begin{align*}
&[D_{0},D_{1}]=-2D_{3},\,[D_{0},D_{2}]=2D_{3},\,[D_{0},D_{3}]=D_{1}-D_{2},\,
[D_{1},D_{2}]=-2\{D_{0},D_{3}\},\\ &[D_{1},D_{3}]=-\{D_{0},D_{1}\}
+\frac{(n-1)(n-3)}2D_{0},\,[D_{2},D_{3}]=\{D_{0},D_{2}\}-\frac{(n-1)(n-3)}2D_{0}.
\end{align*}
For $n=3$ the additional operator $\square$ lies in the centre of
the algebra $\Diff({\bf S}^{3}_{\mathbb{S}})$.

If $n>3$ then there are no relations of the second type. If $n=3$
then there is one relation of the second type:
\begin{equation}
\label{AddInvNR}
\frac12\{D_{1},D_{2}\}-D_{0}^{2}=D_{3}^{2}+\square^{2}.
\end{equation}
It is easy to verify that the operators
$D^{*}_{1}=D_{0}^{2}+D_{1}+D_{2}$ and
$$D_{2}^{*}=\frac12\{D_{1},D_{2}\}-D_{3}^{2}+\left(1-\frac{(n-3)(n-1)}4\right)(D_{1}+D_{2})$$
lie in the centre of the algebra $\Diff({\bf
S}^{n}_{\mathbb{S}})$. If $n=3$ it holds
$D^{*}_{2}=\square^{2}+D_{1}^{*}$ due to (\ref{AddInvNR}). In this
case two operators $D_{1}^{*}$ and $\square$ generate the centre
of the algebra $\Diff({\bf S}^{3}_{\mathbb{S}})$.

Using substitution (\ref{GTR}) we obtain analogous relations for
the algebra $\Diff(\mathbf{H}^{n}(\mathbb{R})_{\mathbb{S}})$.

The commutative relation are now as follows:
\begin{align*}
&[\D_{0},\D_{1}]=2\D_{3},\,[\D_{0},\D_{2}]=2\D_{3},\,[\D_{0},\D_{3}]=\D_{2}+\D_{1},\,
[\D_{1},\D_{2}]=-2\{\D_{0},\D_{3}\},\\
&[\D_{1},\D_{3}]=-\{\D_{0},\D_{1}\}
-\frac{(n-1)(n-3)}2\D_{0},\,[\D_{2},\D_{3}]=\{\D_{0},\D_{2}\}-\frac{(n-1)(n-3)}2\D_{0},
\end{align*}
and for $n=3$ also
$$
[\D_{0},\square]=[\D_{1},\square]=[\D_{2},\square]=[\D_{3},\square]=0.
$$
The first three relations were found in \cite{Reimann}, but the
other relations were not calculated there.

If $n>3$ then there are no relations of the second type. If $n=3$
then there is one relation of the second type analogous to
(\ref{AddInvNR}):
\begin{equation}\label{AddInvNR1}
\frac12\{\D_{1},\D_{2}\}-\D_{0}^{2}=\D_{3}^{2}+\bar\square^{2}.
\end{equation}

The operators $\D^{*}_{1}=\D_{0}^{2}+\D_{1}-\D_{2}$ and
$$\D_{2}^{*}=\frac12\{\D_{1},\D_{2}\}-\D_{3}^{2}+\left(1-\frac{(n-3)(n-1)}4\right)(\D_{1}-\D_{2})$$
lie in the centre of the algebra
$\Diff(\mathbf{H}^{n}(\mathbb{R})_{\mathbb{S}})$ and if $n=3$ it
holds $\D^{*}_{2}=\bar\square^{2}+\D_{1}^{*}$ due to
(\ref{AddInvNR1}). In this case the operators $\D_{1}^{*}$ and
$\bar\square$ generate the centre of the algebra $\Diff({\bf
H}^{3}_{\mathbb{S}})$.
%
%

\section{The model of the space $\mathbf{P}^{2}(\mathbb{C}a)$} \label{octo}
\markright{\ref{octo}}

Our description of Caley algebra $\mathbb{C}a$ and the octonionic
projective plane $\mathbf{P}^{2}(\mathbb{C}a)$ in this section is
based on \cite{Baez}, \cite{Post}, \cite{Adams}.

\subsection{The algebra $\mathbb{C}a$}\label{Ca}

According to Frobenius theorem there are only four finite
dimensional division algebras over $\mathbb{R}$: $\mathbb{R}$
itself and algebras $\mathbb{C},\mathbb{H},\mathbb{C}a$. The
latter is an eight dimensional normed division algebra of
octonions. It is noncommutative and nonassociative, but
alternative, i.e. for any two elements $\xi,\eta\in\mathbb{C}a$ it
holds $(\xi\eta)\eta=\xi(\eta\eta)$ and
$\xi(\xi\eta)=(\xi\xi)\eta$. The group of all automorphisms of
$\mathbb{C}a$ is the exceptional simple compact $14$-dimensional
Lie group $G_{2}$. The standard base of $\mathbb{C}a$ over
$\mathbb{R}$ is $\{e_{i}\}_{i=0}^{7}$, where
$e_{0}=1\in\mathbb{R}$ and
$e_{i}^{2}=-1,\,e_{i}e_{j}=-e_{j}e_{i},\,i,j=1,\dots,7,\,i\ne j$.
The elements $\{e_{i}\}_{i=1}^{7}$ are multiplied according to the
following scheme:

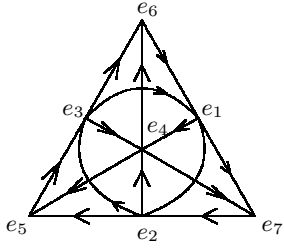
\begin{picture}(40,36)
\qbezier(5,5)(12.5,18)(20,31) \qbezier(5,5)(20,5)(35,5)
\qbezier(35,5)(27.5,18)(20,31) \qbezier(12.5,18)(20,14)(35,5)
\qbezier(20,5)(20,14)(20,31) \qbezier(5,5)(20,14)(27.5,18)
\qbezier(20,5)(9,9)(12.5,18) \qbezier(20,5)(31,9)(27.5,18)
\qbezier(12.5,18)(20,26)(27.5,18) \put(20.7,16.2){{\small
$e_{4}$}} \put(28,18.5){{\small $e_{1}$}} \put(19.4,32){{\small
$e_{6}$}} \put(9.5,18.5){{\small $e_{3}$}} \put(19.4,2){{\small
$e_{2}$}} \put(2,3){{\small $e_{5}$}} \put(36,3){{\small $e_{7}$}}
\qbezier(17.2,26)(15,24)(15,24)
\qbezier(17,25.8)(16.6,23.3)(16.6,23.3)
\qbezier(9,12)(6.8,10)(6.8,10) \qbezier(9,12)(8.4,9.3)(8.4,9.3)
\qbezier(20,11)(19,9)(19,9)\qbezier(20,11)(21,9)(21,9)
\qbezier(20,25)(19,23)(19,23)\qbezier(20,25)(21,23)(21,23)
\qbezier(11,5)(13,4)(13,4)\qbezier(11,5)(13,6)(13,6)
\qbezier(28,5)(30,6)(30,6)\qbezier(28,5)(30,4)(30,4)
\qbezier(10,8)(12,10)(12,10)\qbezier(10,8)(12.6,8.7)(12.6,8.7)
\qbezier(23.3,21.3)(22,22.2)(22,22.2)\qbezier(23.3,21.3)(21.7,21.3)(21.7,21.3)
\qbezier(16,6.95)(17.6,6.7)(17.6,6.7)\qbezier(16,6.95)(16.85,5.65)(16.85,5.65)
\qbezier(31.74,10.74)(31.4,12.2)(31.4,12.2)\qbezier(31.7,10.7)(30.5,11.5)(30.5,11.5)
\qbezier(23.3,25.4)(22.9,27)(22.9,27)\qbezier(23.3,25.3)(22,26.4)(22,26.4)
\qbezier(24.3,16.3)(25.8,17.7)(25.8,17.7)\qbezier(24.3,16.3)(26.2,16.6)(26.2,16.6)
\qbezier(31.1,7.3)(28.7,8)(28.7,8)\qbezier(31.2,7.3)(29.4,9)(29.4,9)
\qbezier(16.5,15.8)(14.2,16.4)(14.2,16.4)\qbezier(16.5,15.9)(14.7,17.5)(14.7,17.5)
\end{picture}

Here $e_{i}e_{j}=e_{k}$ if these elements lie on one line or on
the circle and are ordered by arrows as $e_{i},e_{j},e_{k}$. The
conjugation $\iota: \mathbb{C}a\mapsto\mathbb{C}a$ acts as
$\iota(e_{0})\equiv\overline{e}_{0}=e_{0},\,\iota(e_{i})\equiv\overline{e}_{i}=-e_{i},i=1,\dots,7$
and is extended by linearity over whole $\mathbb{C}a$. Let
$\RE\xi=\frac12(\xi+\bar\xi),\,\im\xi=\frac12(\xi-\bar\xi),\,\xi\in\mathbb{C}a$.
Define the scalar product in $\mathbb{C}a$ by the formula:
$\langle\eta,\xi\rangle=\dfrac12(\bar\eta\xi+\bar\xi\eta)=\RE(\bar\xi\eta)=\RE(\bar\eta\xi)\in\mathbb{R}$
and the norm by the formula
$\|\eta\|=\langle\eta,\eta\rangle^{1/2}$. In the algebra
$\mathbb{C}a$ every two elements generate an associative
subalgebra and the following {\it central Moufang identity} is
valid:
\begin{equation}\label{Moufang}
u\cdot xy\cdot u=ux\cdot yu,\, u,x,y\in \mathbb{C}a.
\end{equation}
Here we use the notations $u\cdot xy:=u(xy),\, xy\cdot u:=(xy)u$.

 There are the following description of spinor and
vector representations (all 8-dimensional) of the group $\Spin(8)$
in $\mathbb{C}a$ \cite{Post}, \cite{Onishchik}, which will be used
later. Define linear operators in $\mathbb{C}a$:
\begin{align*}
L_{\alpha}:\,&\xi\mapsto\frac12e_{\alpha}\xi,\,\alpha=1,\dots,7,\,\xi\in\mathbb{C}a,\\
L_{\alpha,\beta}:\,&\xi\mapsto\frac12e_{\alpha}(e_{\beta}\xi),\,1\leqslant
\alpha<\beta\leqslant 7,\,\xi\in\mathbb{C}a.
\end{align*}
This operators are generators of the left spinor representation of
the group $\Spin(8)$, i.e. they are the images of some base of the
Lie algebra $\mathfrak{spin}(8)$ under this representation.
Similarly, operators
\begin{align*}
R_{\alpha}:\,&\xi\mapsto\frac12\xi
e_{\alpha},\,\alpha=1,\dots,7,\,\xi\in\mathbb{C}a,\\
R_{\alpha,\beta}:\,&\xi\mapsto\frac12(\xi
e_{\beta})e_{\alpha},\,1\leqslant \alpha<\beta\leqslant
7,\,\xi\in\mathbb{C}a
\end{align*}
are generators of the right spinor representation of the group
$\Spin(8)$. All these operators are skew symmetric w.r.t. the
scalar product in $\mathbb{C}a$.

Formulae above define operators
$L_{\alpha,\beta},R_{\alpha,\beta}$ also for $1\leqslant
\beta<\alpha\leqslant 7$. If $\mathbb{C}a'$ is the space of pure
imaginary octonions, $u\in\mathbb{C}a',\,\xi\in\mathbb{C}a$, then
due to the alternativity of $\mathbb{C}a$:
$$
\xi u\cdot u=\xi u^{2}=-\xi|u|^{2}=-|u|^{2}\xi=u\cdot u\xi.
$$
For $u=e_{\alpha}+e_{\beta},\,1\leqslant\alpha<\beta\leqslant 7$
it holds
\begin{align*}
-2\xi&=-\xi|e_{\alpha}+e_{\beta}|^{2}=\xi(e_{\alpha}+e_{\beta})\cdot(e_{\alpha}+e_{\beta})=
\xi e_{\alpha}\cdot e_{\alpha}+\xi e_{\alpha}\cdot e_{\beta}+\xi
e_{\beta}\cdot e_{\alpha}+\xi e_{\beta}\cdot e_{\beta}\\
&=-\xi+\xi e_{\alpha}\cdot e_{\beta}+\xi e_{\beta}\cdot
e_{\alpha}-\xi
\end{align*}
and $\xi e_{\alpha}\cdot e_{\beta}+\xi e_{\beta}\cdot
e_{\alpha}=0$. Similarly, $e_{\alpha}\cdot e_{\beta}\xi +
e_{\beta}\cdot e_{\alpha}\xi=0$. For $0\leqslant i,j\leqslant
7,\,i\ne j$ we can write more general formulae, useful in the
following:
\begin{equation}
\label{important} e_{i}\cdot e_{j}\xi=-\bar  e_{j}\cdot \bar
e_{i}\xi,\,\xi e_{i}\cdot e_{j}=-\xi\bar e_{j}\cdot \bar
e_{i},\,\xi\in\mathbb{C}a.
\end{equation}

In particular, we have
$L_{\alpha,\beta}=-L_{\beta,\alpha},\,R_{\alpha,\beta}=-R_{\beta,\alpha},\,1\leqslant\alpha,\beta\leqslant
7,\alpha\ne\beta$.

For the element $g\in\Spin(8)$ denote by $g^{L},g^{R}$ and $g^{V}$
its images under left spinor, right spinor and vector
representation respectively. The following proposition is a
version of the {\it triality principle} for the group
$\Spin(8)$\footnote{Other versions of this principle are in
\cite{Baez}.}.
\begin{proposit}[\cite{Post}]\label{triality}
For any element $g\in\Spin(8)$ it holds
\begin{equation}\label{trial}
g^{V}(\xi\eta)=g^{L}(\xi)\cdot
g^{R}(\eta),\,\xi,\eta\in\mathbb{C}a.
\end{equation}
Conversely, if $A,B,C$ are orthogonal operators
$\mathbb{C}a\mapsto\mathbb{C}a$ such that
$$
A(\xi\eta)=B(\xi)\cdot C(\eta),
$$
for any $\xi,\eta\in\mathbb{C}a$, then there exist unique
$g\in\Spin(8)$ such that $A=g^{V},B=g^{L},C=g^{R}$.
\end{proposit}
From equation (\ref{trial}) we obtain its infinitesimal analogs:
\begin{align}\label{IZ1}
V_{i}(\xi\eta)&=L_{i}(\xi)\cdot\eta+\xi\cdot
R_{i}(\eta),\,i=1,\dots,7,\\\label{IZ2}
V_{i,j}(\xi\eta)&=L_{i,j}(\xi)\cdot\eta+\xi\cdot
R_{i,j}(\eta),\,1\leqslant i<j\leqslant 7,\xi,\eta\in\mathbb{C}a,
\end{align}
where $V_{i}$ and $V_{i,j}$ are generators of the vector
representation of the group $\Spin(8)$.

\subsection{The Jordan algebra $\mathfrak{h}_{3}(\mathbb{C}a)$}

The {\it Hermitian conjugation} $A\mapsto A^{*}$ for square matrix
with octonion entries is defined as the composition of octonionic
conjugation and transposition of $A$, similar to complex or
quaternion cases. A matrix $A$ is called {\it Hermitian} iff
$A^{*}=A$. {\it The simple exceptional Jordan algebra}
$\mathfrak{h}_{3}(\mathbb{C}a)$ consists of all Hermitian $3\times
3$ matrices with octonion entries. It is endowed with the Jordan
commutative multiplication: $$X\circ
Y=\dfrac12(XY+YX),\,X,Y\in\mathfrak{h}_{3}(\mathbb{C}a).
$$
This multiplication satisfies to the identity $(X^{2}\circ Y)\circ
X=X^{2}\circ(Y\circ X)$ which is the condition for an algebra with
commutative (but not necessarily associative) multiplication to be
Jordan. The Jordan algebra $\mathfrak{h}_{3}(\mathbb{C}a)$ is
$27$-dimensional over $\mathbb{R}$. Every its element can be
represented in the form:
\begin{equation}\label{X}
X=a_{1}E_{1}+a_{2}E_{2}+a_{3}E_{3}+X_{1}(\xi_{1})+X_{2}(\xi_{2})+X_{3}(\xi_{3}),
\end{equation}
where
\begin{align*}
E_{1}&=\left(\begin{array}{ccc} 1 & 0 & 0 \\ 0 & 0 & 0 \\ 0 & 0 &
0
\end{array}\right),\,
E_{2}=\left(\begin{array}{ccc} 0 & 0 & 0 \\ 0 & 1 & 0 \\ 0 & 0 & 0
\end{array}\right),\,
E_{3}=\left(\begin{array}{ccc} 0 & 0 & 0 \\ 0 & 0 & 0 \\ 0 & 0 & 1
\end{array}\right),\\
X_{1}(\xi)&=\left(\begin{array}{ccc} 0 & 0 & 0 \\ 0 & 0 & \xi \\ 0
& \bar\xi & 0 \end{array}\right),\,
X_{2}(\xi)=\left(\begin{array}{ccc} 0 & 0 & \bar\xi \\ 0 & 0 & 0
\\ \xi & 0 & 0 \end{array}\right),\,
X_{3}(\xi)=\left(\begin{array}{ccc} 0 & \xi & 0 \\ \bar\xi & 0 & 0
\\ 0 & 0 & 0
\end{array}\right),
\end{align*}
$a_{i}\in \mathbb{R},\xi_{i}\in\mathbb{C}a,\,i=1,2,3.$ It is easy
to show that
\begin{align}\label{circ}
E_{i}\circ E_{j}&=\left\{E_{i},\,\mbox{if}\;i=j,\atop
0,\,\mbox{if}\; i\ne j,\right.\nonumber\\ E_{i}\circ
X_{j}(\xi)&=\left\{0,\,\mbox{if}\;i=j,\atop
\frac12X_{j}(\xi),\,\mbox{if}\;i\ne j,\right.\\ X_{i}(\xi)\circ
X_{j}(\eta)&=\left\{(\xi,\eta)(E-E_{i}),\,\mbox{if}\;i=j,\atop
\frac12 X_{i+j}(\overline{\xi\eta}),\,\mbox{if}\;j\equiv i+1\mod
3,\right.\nonumber
\end{align}
where $E=E_{1}+E_{2}+E_{3}$ is the unit matrix. In the last
formula all indices are considered modulo $3$.

The group of all automorphisms of the Jordan algebra
$\mathfrak{h}_{3}(\mathbb{C}a)$ is the exceptional simple compact
$52$-dimensional real Lie group $F_{4}$. This group conserves the
following bilinear and trilinear functionals:
$\mathcal{A}(X,Y)=\Tr(X\circ
Y),\,\mathcal{B}(X,Y,Z)=\mathcal{A}(X\circ Y,Z)$. Conversely every
linear operator $\mathfrak{h}_{3}(\mathbb{C}a)\mapsto
\mathfrak{h}_{3}(\mathbb{C}a)$, conserving these two functionals,
lies in $F_{4}$.

Define the norm of the element (\ref{X}) as
$\|X\|^{2}=\mathcal{A}(X,X)=\sum_{i=1}^{3}(a_{i}^{2}+2|\xi|^{2})$.
The last equality is the consequence of (\ref{circ}).
\begin{theore}[Freudenthal]
For any $X\in\mathfrak{h}_{3}(\mathbb{C}a)$ there exists an
automorphism $\Phi\in F_{4}$, such that
\begin{equation}\label{Fr}
\Phi X=\lambda_{1}E_{1}+\lambda_{2}E_{2}+\lambda_{3}E_{3},
\end{equation}
where $\lambda_{1}\geqslant\lambda_{2}\geqslant\lambda_{3}$, and
the form (\ref{Fr}) is uniquely determined by $X$. Two elements
from $\mathfrak{h}_{3}(\mathbb{C}a)$ lie on the same orbit of
$F_{4}$ iff their diagonal forms (\ref{Fr}) are the same.
\end{theore}
\subsection{The octonionic projective plane
$\mathbf{P}^{2}(\mathbb{C}a)$} Elements
$X\in\mathfrak{h}_{3}(\mathbb{C}a)$ satisfying conditions
\begin{equation}\label{cond}
X^{2}=X,\; \Tr X=1
\end{equation}
form the octonionic projective plane
$\mathbf{P}^{2}(\mathbb{C}a)$, which is a $16$-dimensional real
manifold. Automorphisms of $\mathfrak{h}_{3}(\mathbb{C}a)$
conserves equations (\ref{cond}) and the group $F_{4}$ acts on
$\mathbf{P}^{2}(\mathbb{C}a)$. From the Freudenthal theorem and
equations (\ref{cond}) it follows that every element of
$\mathbf{P}^{2}(\mathbb{C}a)$ can be transformed by an appropriate
element of $F_{4}$ to the element $E_{1}$. Thus
$\mathbf{P}^{2}(\mathbb{C}a)$ is a homogeneous space of the group
$F_{4}$ and calculations in \cite{Post} (lecture 16) shows, that
the stationary subgroup of every point
$X\in\mathbf{P}^{2}(\mathbb{C}a)$ is isomorphic to the group
$\Spin (9)$.

Let
\begin{equation*}
X=(1+a_{1})E_{1}+a_{2}E_{2}+a_{3}E_{3}+X_{1}(\xi_{1})+X_{2}(\xi_{2})+X_{3}(\xi_{3})
\in\mathbf{P}^{2}(\mathbb{C}a),
\end{equation*}
where $a_{i},|\xi_{i}|,i=1,2,3$ are tending to zero. Then due to
(\ref{circ}) we have
\begin{equation*}
X\circ
X=(1+2a_{1})E_{1}+X_{2}(\xi_{2})+X_{3}(\xi_{3})+o\left(\sum_{i=1}^{3}(a_{i}^{2}+|\xi|^{2})\right)
\end{equation*}
and the equality $X\circ X=X$ implies
$a_{1}=a_{2}=a_{3}=0,\xi_{1}=0$. It means that we can identify the
tangent space $T_{E_{1}}\mathbf{P}^{2}(\mathbb{C}a)$ with the set
$\{X_{2}(\xi_{2})+X_{3}(\xi_{3})\,|\,\xi_{1},\xi_{2}\in\mathbb{C}a\}$.

Let $K\subset F_{4}$ be the stationary subgroup corresponding to
the point $E_{1}$ and acting by automorphisms in the space
$T_{E_{1}}\mathbf{P}^{2}(\mathbb{C}a)\simeq\{X_{2}(\xi_{2})+X_{3}(\xi_{3})\,|\,\xi_{1},\xi_{2}\in\mathbb{C}a\}$.
Let $K_{0}$ be the stationary subgroup of $K$, corresponding to
the element $X_{3}(1)\in T_{E_{1}}\mathbf{P}^{2}(\mathbb{C}a)$.

According to the section \ref{spes} we are to calculate the
$K_{0}$-action on $T_{E_{1}}\mathbf{P}^{2}(\mathbb{C}a)$. For any
element $X\in\mathfrak{h}_{3}(\mathbb{C}a)$ let $\Ann
X:=\left\{Y\in\mathfrak{h}_{3}(\mathbb{C}a)\,|\,Y\circ
X=0\right\}$. Being an automorphism of the algebra
$\mathfrak{h}_{3}(\mathbb{C}a)$, an element $\Phi\in K_{0}$
conserves the space $\Ann X_{3}(1)$. It follows from (\ref{circ})
that
$$
\Ann X_{3}(1)=\left\{a(E_{1}-E_{2})+bE_{3}+X_{3}(\xi)\,|\,a,b\in
\mathbb{R},\,\xi\in\mathbb{C}a'\right\}.
$$
Let $\Phi(E_{1}-E_{2})=a(E_{1}-E_{2})+bE_{3}+X_{3}(\xi)$, then we
have
$$1=\mathcal{A}\left(E_{1}-E_{2},E_{1}\right)=\mathcal{A}\left(\Phi(E_{1}-E_{2}),\Phi(E_{1})\right)=
\mathcal{A}\left(a(E_{1}-E_{2})+bE_{3}+X_{3}(\xi),E_{1}\right)=a.
$$
This implies $\Phi(E_{1}-E_{2})=E_{1}-E_{2}+bE_{3}+X_{3}(\xi)$ and
the equality $\|E_{1}-E_{2}\|=\|\Phi(E_{1}-E_{2})\|$ gives
$b=0,\xi=0$. It means that $\Phi(E_{2})=E_{2}$ and therefore
$\Phi(E_{3})=\Phi(E-E_{1}-E_{2})=E-E_{1}-E_{2}=E_{3}$. Thus the
group $K_{0}$ conserves elements $E_{1},E_{2},E_{3}$.

Let $K'$ be the subgroup of $F_{4}$ conserving element
$E_{1},E_{2},E_{3}$. We see that $K_{0}\subset K'\subset K$. Since
$\Ann
E_{1}=\left\{a_{2}E_{2}+a_{3}E_{3}+X_{1}(\xi),\,a_{1},a_{2}\in
\mathbb{R},\,\xi\in\mathbb{C}a\right\}$, then the group $K'$ maps
$X_{1}(\xi)\mapsto X_{1}(\tilde\xi)$ and similarly
$X_{i}(\xi_{i})\mapsto X_{i}(\tilde\xi_{i}),\,i=1,2,3$.

Let $\Phi_{i}:\mathbb{C}a\mapsto\mathbb{C}a,\,i=1,2,3$ be
orthogonal operators such that $\Phi
X_{i}(\xi_{i})=X_{i}\left(\Phi_{i}(\xi_{i})\right)$ for $\Phi\in
K'$. The last formula in (\ref{circ}) implies
\begin{align*}
X_{3}\left(\Phi_{3}(\overline{\xi\eta})\right)&=\Phi\left(X_{3}(\overline{\xi\eta})\right)=
2\Phi\left(X_{1}(\xi)\circ
X_{2}(\eta)\right)=2\Phi\left(X_{1}(\xi)\right)\circ\Phi\left(X_{2}(\eta)\right)\\
&=2 X_{1}\left(\Phi_{1}(\xi)\right)\circ
X_{2}\left(\Phi_{2}(\eta)\right)=X_{3}\left(\overline{\Phi_{1}(\xi)\Phi_{2}(\eta)}\right).
\end{align*}
It gives
\begin{equation}\label{fundument}
\Phi_{1}(\xi)\Phi_{2}(\eta)=\overline{\Phi_{3}\left(\overline{\xi\eta}\right)},
\end{equation}
for $\Phi\in K',\,\xi,\eta\in\mathbb{C}a$.

Denote by $\mathbb{C}a_{i},\,i=1,2,3$ the domains for the
operators $\Phi_{i},\,i=1,2,3$. Then
$T_{E_{1}}\mathbf{P}^{2}(\mathbb{C}a) \\ \simeq
\mathbb{C}a_{2}\oplus \mathbb{C}a_{3}$.

The formula (\ref{fundument}) and the proposition \ref{triality}
imply
\begin{proposit}\label{phi}
Operators $\Phi_{1}$ and $\Phi_{2}$ are respectively left and
right spinor representations of the group $\Spin(8)\simeq K'$ and
the composition $\iota\circ\Phi_{3}\circ\iota$ is the vector
representation of $\Spin(8)$.
\end{proposit}
The group $\Spin(8)$ is the universal (double) covering of the
group $\SO(8)$ and their Lie algebras $\mathfrak{spin}(8)$ and
$\mathfrak{so}(8)$ are isomorphic.

Now consider representations of the Lie algebra $\mathfrak{k}'$ of
the group $K'$ in $\mathbb{C}a_{i},\,i=1,2,3$. All these
representations are faithful. For $A\in\mathfrak{k}'$ denote by
$A^{(i)}$ the corresponding skew symmetric operator in
$\mathbb{C}a_{i},\,i=1,2,3$. From (\ref{fundument}) we obtain the
following infinitesimal analog of the triality principle:
\begin{equation}\label{trinf}
A^{(1)}(\xi)\cdot\eta+\xi\cdot
A^{(2)}(\eta)=\overline{A^{(3)}(\overline{\xi\eta})}.
\end{equation}
From (\ref{IZ1}) and (\ref{IZ2}) we obtain that if $A^{(1)}=L_{i}$
(respectively $A^{(1)}=L_{i,j}$) then
$A^{(2)}=R_{i},\,A^{(3)}=\iota\circ V_{i}\circ\iota$ (respectively
$A^{(2)}=R_{i,j},\,A^{(3)}=\iota\circ V_{i,j}\circ\iota$).

Let us identify the algebra $\mathfrak{k}'$ with its vector
representation in $\mathbb{C}a_{3}$, in particular we put $A\equiv
A^{(3)}$ for $A\in\mathfrak{k}'$. By $\varkappa$ denote the
inclusion $\mathfrak{k}'$ into the Lie algebra $\mathfrak{f}_{4}$
corresponding to the group $F_{4}$.

By the definition, the Lie algebra $\mathfrak{k}_{0}$ of the group
$K_{0}\subset K'$ consists of the skew symmetric operators in
$\mathbb{C}a_{3}$, transforming $1\in\mathbb{C}a_{3}$ into $0$.
The group $K_{0}$ is isomorphic to $\Spin(7)$, acting in
$\mathbb{C}a_{1}$ by the left spinor representation, in
$\mathbb{C}a_{2}$ by the right spinor representation (equivalent
for $\Spin(7)$ to the left one, see (\ref{phi12}) below), and in
$\mathbb{C}a_{3}'$ by the vector representation, which are
restrictions of analogous representations of $K'\simeq\Spin(8)$.

Let $\mathfrak{m}$ be the space of $3\times 3$ skew-Hermitian
matrices with octonion entries and the zero trace. Let
\begin{align*}
Y_{1}(\xi)&=\left(\begin{array}{ccc} 0 & 0 & 0 \\ 0 & 0 & \xi \\ 0
& -\bar\xi & 0 \end{array}\right),\,
Y_{2}(\xi)=\left(\begin{array}{ccc} 0 & 0 & -\bar\xi \\ 0 & 0 & 0
\\ \xi & 0 & 0 \end{array}\right),\,
Y_{3}(\xi)=\left(\begin{array}{ccc} 0 & \xi & 0 \\ -\bar\xi & 0 &
0
\\ 0 & 0 & 0
\end{array}\right),\xi\in\mathbb{C}a
\end{align*}
be elements from $\mathfrak{m}$ and the linear subspace
$\mathfrak{m}_{0}\subset\mathfrak{m}$ consists of elements of the
form
$$
\sum_{i=1}^{3}Y_{i}(\xi_{i}),\,\xi_{i}\in\mathbb{C}a.
$$
From \cite{Post} (lecture 16) we can extract the following
proposition:
\begin{proposit}\label{Post}
For $Y\in\mathfrak{m}$ the linear operator $\ad
Y:\mathfrak{h}_{3}(\mathbb{C}a)\mapsto\mathfrak{h}_{3}(\mathbb{C}a)$,
acting according to the formula $\ad Y(X)=YX-X
Y,\,X\in\mathfrak{h}_{3}(\mathbb{C}a)$ is the differentiation of
the algebra $\mathfrak{h}_{3}(\mathbb{C}a)$. Thus the space
$\mathfrak{m}$ is contained in $\mathfrak{f}_{4}$. There is the
expansion into the direct sum of linear spaces
$$
\mathfrak{f}_{4}\simeq\mathfrak{k}'\oplus\mathfrak{m}_{0}
$$
with the following commutator relations
\begin{align}\label{OctCom1}
&\left[\varkappa A,\ad Y_{i}(\xi)\right]=\ad
Y_{i}(A^{(i)}\xi),\,i=1,2,3
\\\label{OctCom2} &\left[\ad Y_{i}(\xi),\ad
Y_{j}(\eta)\right]=\left\{\varkappa
C_{i,\xi,\eta},\,\mbox{if}\;j=i,\atop \ad
Y_{i+2}\left(-\overline{\xi\eta}\right),\,\mbox{if}\;j=i+1,\right.
\end{align}
where $A\equiv A^{(3)}\in\mathfrak{k}',\,\xi,\eta\in\mathbb{C}a$,
operators $A^{(i)}$ are from (\ref{trinf}), the indices in the
last equation are considered modulo $3$ and skew-Hermitian
operators
$C_{i,\xi,\eta}:\mathbb{C}a_{3}\mapsto\mathbb{C}a_{3},\,i=1,2,3$
are given by the following formulas:
\begin{align}\label{C3}
&C_{1,\xi,\eta}:\zeta\mapsto
\zeta\xi\cdot\bar\eta-\zeta\eta\cdot\bar\xi,\nonumber\\
&C_{2,\xi,\eta}:\zeta\mapsto
\bar\eta\cdot\xi\zeta-\bar\xi\cdot\eta\zeta,\,\zeta\in\mathbb{C}a\\
&C_{3,\xi,\eta}:\zeta\mapsto
4(\xi,\zeta)\eta-4(\eta,\zeta)\xi.\nonumber
\end{align}
\end{proposit}
The action of operators $\varkappa C_{i,\xi,\eta}$ on the spaces
$\mathbb{C}a_{1}$ and $\mathbb{C}a_{2}$ is obtained from
(\ref{C3}) by the cyclic permutation of indices:
\begin{align}\label{C12}
&\left.\varkappa
C_{1,\xi,\eta}\right|_{\mathbb{C}a_{1}}:\zeta\mapsto
4(\xi,\zeta)\eta-4(\eta,\zeta)\xi,\nonumber\\ &\left.\varkappa
C_{2,\xi,\eta}\right|_{\mathbb{C}a_{1}}:\zeta\mapsto
\zeta\xi\cdot\bar\eta-\zeta\eta\cdot\bar\xi,\nonumber\\
&\left.\varkappa
C_{3,\xi,\eta}\right|_{\mathbb{C}a_{1}}:\zeta\mapsto
\bar\eta\cdot\xi\zeta-\bar\xi\cdot\eta\zeta,\\ &\left.\varkappa
C_{1,\xi,\eta}\right|_{\mathbb{C}a_{2}}:\zeta\mapsto
\bar\eta\cdot\xi\zeta-\bar\xi\cdot\eta\zeta,\nonumber\\
&\left.\varkappa
C_{2,\xi,\eta}\right|_{\mathbb{C}a_{2}}:\zeta\mapsto
4(\xi,\zeta)\eta-4(\eta,\zeta)\xi,\nonumber\\ &\left.\varkappa
C_{3,\xi,\eta}\right|_{\mathbb{C}a_{2}}:\zeta\mapsto
\zeta\xi\cdot\bar\eta-\zeta\eta\cdot\bar\xi.\nonumber
\end{align}
Note that in \cite{Post} (lecture 16) analogs of formulae
(\ref{fundument}), (\ref{OctCom1}) and the last formula
(\ref{circ}) contain errors.

\section{Generators of algebras
$\Diff\left(\mathbf{P}^{2}(\mathbb{C}a)_{\mathbb{S}}\right)$ and
$\Diff\left(\mathbf{H}^{2}(\mathbb{C}a)_{\mathbb{S}}\right)$}
\label{CaS}\markright{\ref{CaS} Generators of algebras
$\Diff\left(\mathbf{P}^{2}(\mathbb{C}a)_{\mathbb{S}}\right)$ and
$\Diff\left(\mathbf{H}^{2}(\mathbb{C}a)_{\mathbb{S}}\right)$}

Now we are to specify the construction from section \ref{spes} for
the space $M=\mathbf{P}^{2}(\mathbb{C}a)_{\mathbb{S}}$.

\subsection{The special base in $\mathfrak{a}\oplus\mathfrak{p}_{\lambda}\oplus\mathfrak{k}_{\lambda}
\oplus\mathfrak{p}_{2\lambda}\oplus\mathfrak{k}_{2\lambda}$}

It is easily seen that
$$
[Y_{1}(\xi),E_{1}]=0,\,[Y_{2}(\xi),E_{1}]=X_{2}(\xi),\,[Y_{3}(\xi),E_{1}]=-X_{3}(\xi),\,\xi\in\mathbb{C}a,
$$
so we can identify the space
$T_{E_{1}}\mathbf{P}^{2}(\mathbb{C}a)$ with the space
$\left\{Y_{2}(\xi)+Y_{3}(\eta)|\,\xi,\eta\in\mathbb{C}a\right\}\subset
\mathfrak{m}_{0}$. From (\ref{OctCom1}) we obtain that the
expansion
$$
\left\{Y_{2}(\xi)+Y_{3}(\eta)|\,\xi,\eta\in\mathbb{C}a\right\}=
\left\{Y_{3}(\xi)|\,\xi\in\mathbb{R}\right\}\oplus
\left\{Y_{2}(\xi)|\,\xi\in\mathbb{C}a\right\}\oplus
\left\{Y_{3}(\xi)|\,\xi\in\mathbb{C}a'\right\}
$$
is $\Ad_{K_{0}}$-invariant and by the comparison with sections
\ref{twopoint} and \ref{spes} we can put:
$$
\mathfrak{a}:=\left\{Y_{3}(\xi)|\,\xi\in\mathbb{R}\right\},\,
\mathfrak{p}_{\lambda}:=\left\{Y_{2}(\xi)|\,\xi\in\mathbb{C}a\right\},\,
\mathfrak{p}_{2\lambda}:=\left\{Y_{3}(\xi)|\,\xi\in\mathbb{C}a'\right\}.
$$

Let
$y=\left(E_{1},\frac12X_{3}(1)\right)\in\mathbf{P}^{2}(\mathbb{C}a)_{\mathbb{S}}$,
where $\frac12X_{3}(1)\in\mathbb{S}_{E_{1}}$. We have
$T_{y}\mathbf{P}^{2}(\mathbb{C}a)_{\mathbb{S}}=T_{E_{1}}\mathbf{P}^{2}(\mathbb{C}a)\oplus
T_{\frac12X_{3}(1)}\mathbb{S}_{E_{1}}$ and
$$
T_{\frac12X_{3}(1)}\mathbb{S}_{E_{1}}\simeq\left\{X_{2}(\xi)|\,\xi\in\mathbb{C}a\right\}\oplus
\left\{X_{3}(\xi)|\,\xi\in\mathbb{C}a'\right\}.
$$
Since $\ad
Y_{1}(\xi)\left(X_{3}(1)\right)=-X_{2}(\bar\xi),\,\xi\in\mathbb{C}a$,
the space $\left\{X_{2}(\xi)|\,\xi\in\mathbb{C}a\right\}\subset
T_{\frac12X_{3}(1)}\mathbb{S}_{E_{1}}$ is identified with the
space
$\left\{Y_{1}(\xi)|\,\xi\in\mathbb{C}a\right\}\subset\mathfrak{m}_{0}$.
Since
$\dim_{\mathbb{R}}\left\{Y_{1}(\xi)|\,\xi\in\mathbb{C}a\right\}=8=\dim_{\mathbb{R}}\mathfrak{p_{\lambda}}$
we are to denote
$\mathfrak{k_{\lambda}}:=\left\{Y_{1}(\xi)|\,\xi\in\mathbb{C}a\right\}$.
Thus
$\mathfrak{m}_{0}=\mathfrak{a}\oplus\mathfrak{p}_{\lambda}\oplus\mathfrak{p}_{2\lambda}\oplus
\mathfrak{k}_{\lambda}$.

Denote by $A_{ij}\in\mathfrak{k}',\,i\ne j$ the generators of the
rotation in the $2$-dimensional plane, containing elements
$e_{i},e_{j}\in\mathbb{C}a_{3}$, such that
$A_{ij}e_{j}=e_{i},\,A_{ij}e_{i}=-e_{j}$. The operators
$A_{ij},\,1\leqslant i<j\leqslant 7$ are the base of the algebra
$\mathfrak{k}_{0}$. Similar to the quaternion case the subspace
$\mathfrak{q}$ of the algebra $\mathfrak{k}'$ with the base
$A_{0\alpha}=:A_{\alpha},\,\alpha=1,\dots,7$ is
$\Ad_{K_{0}}$-invariant and is identified through the
$K_{0}$-action on $T_{y}\mathbf{P}^{2}(\mathbb{C}a)_{\mathbb{S}}$
with the space
$\left\{X_{3}(\xi)|\,\xi\in\mathbb{C}a'\right\}\subset
T_{\frac12X_{3}(1)}\mathbb{S}_{E_{1}}$. Therefore we define
$\mathfrak{k}_{2\lambda}:=\mathfrak{q}$.

\begin{Lem}\label{lem1} It holds
\begin{align*}
A^{(1)}_{\alpha}&=L_{\alpha},\,A^{(2)}_{\alpha}=R_{\alpha},
\,A^{(1)}_{\alpha\beta}=L_{\beta,\alpha},\,A^{(2)}_{\alpha\beta}=R_{\beta,\alpha},\\
C_{3,e_{\alpha},e_{\beta}}&=4A_{\beta,\alpha},\,C_{3,e_{0},e_{\alpha}}=-4A_{\alpha},\,
\alpha,\beta=1,\dots,7,\alpha\neq\beta
\end{align*}
\end{Lem}
\begin{proof}
From (\ref{trinf}) we have
\begin{equation*}
A^{(3)}(\xi)=\overline{A^{(1)}(\overline{\xi})+\overline{\xi}A^{(2)}(1)}.
\end{equation*}
Let $A^{(1)}=L_{\alpha}$, then $A^{(2)}=R_{\alpha}$ and
$A^{(3)}(e_{k})=\frac12\overline{\left(e_{\alpha}\bar e_{k}+\bar
e_{k}e_{\alpha}\right)}=-\frac12\left(e_{k}e_{\alpha}+e_{\alpha}e_{k}\right)$.
If $1\leqslant k\neq\alpha$, then
$e_{k}e_{\alpha}=-e_{\alpha}e_{k}$ and $A^{(3)}(e_{k})=0$.
Therefore $A^{(3)}=A_{\alpha}$, since
$A^{(3)}(1)=-e_{\alpha},\,A^{(3)}(e_{\alpha})=1$. This proves
$A^{(1)}_{\alpha}=L_{\alpha},\,A^{(2)}_{\alpha}=R_{\alpha}$.

Let now $A^{(1)}=L_{\beta,\alpha}$, then
$A^{(2)}=R_{\beta,\alpha}$ and
$A^{(3)}(e_{k})=\frac12\overline{\left(e_{\beta}\cdot
e_{\alpha}\bar e_{k}+\bar e_{k}\cdot
e_{\alpha}e_{\beta}\right)}=\frac12\left(e_{k}e_{\alpha}\cdot
e_{\beta} +e_{\beta}e_{\alpha}\cdot e_{k}\right)$. It is easy to
verify by direct computation that if $\alpha=1,\beta=2$ then
$A^{(3)}(e_{k})=0$, for $k\neq 1,2$ and
$A^{(3)}(e_{1})=-e_{2},\,A^{(3)}(e_{2})=e_{1}$. Thus
$L^{(3)}_{\beta,\alpha}=A_{12}$. Therefore
$L^{(3)}_{\beta,\alpha}=A_{\alpha\beta}$ for any other pair of
$e_{\alpha},e_{\beta}$, since the group $G_{2}$ of automorphisms
of $\mathbb{C}a$ acts transitively on any pair of imaginary units
\cite{Post} (lecture 15). This proves
$A^{(1)}_{\alpha\beta}=L_{\beta,\alpha},\,A^{(2)}_{\alpha\beta}=R_{\beta,\alpha}$.

The last two equalities of this lemma are obvious.
\end{proof}

Let summarize these reasoning in the following proposition:
\begin{proposit}\label{octbase}
Let \begin{align*} \Lambda&:=\dfrac12\ad
Y_{3}(e_{1}),\,e_{2\lambda,\alpha}:=\dfrac12\ad
Y_{3}(e_{\alpha}),\,f_{2\lambda,\alpha}:=\varkappa A_{\alpha},\\
e_{\lambda,i}&:=-\dfrac12 \ad Y_{2}(\bar
e_{i}),\,f_{\lambda,i}:=\dfrac12\ad Y_{1}(e_{i}),\,\tilde
A_{\alpha\beta}:=\varkappa A_{\alpha\beta},
\end{align*}
where latin indices vary from $0$ to $7$ and greek ones (except
$\lambda$) vary from $1$ to $7$. We have the following commutator
relations:
\begin{align*}
[\Lambda,e_{2\lambda,\alpha}]&=-f_{2\lambda,\alpha},\,[\Lambda,f_{2\lambda,\alpha}]=e_{2\lambda,\alpha},\,
[\Lambda,e_{\lambda,i}]=-\frac12f_{\lambda,i},\,[\Lambda,f_{\lambda,i}]=\frac12e_{\lambda,i},\,
[\Lambda,\tilde A_{\alpha\beta}]=0,\\
[e_{2\lambda,\alpha},e_{2\lambda,\beta}]&=\tilde
A_{\beta\alpha},\,
[e_{2\lambda,\alpha},f_{2\lambda,\beta}]=-\delta_{\alpha\beta}\Lambda,\,
[f_{2\lambda,\alpha},f_{2\lambda,\beta}]=\tilde A_{\beta\alpha},\,
[e_{2\lambda,\alpha},e_{\lambda,j}]=\frac12f_{\lambda,e_{\alpha}e_{j}},\\
[e_{2\lambda,\alpha},f_{\lambda,j}]&=\frac12e_{\lambda,e_{\alpha}e_{j}},\,
[f_{2\lambda,\alpha},e_{\lambda,j}]=-\frac12e_{\lambda,e_{\alpha}e_{j}},\,
[f_{2\lambda,\alpha},f_{\lambda,j}]=\frac12f_{\lambda,e_{\alpha}e_{j}},\\
[e_{\lambda,i},e_{\lambda,j}]&=\frac14 \varkappa C_{2,\bar
e_{i},\bar e_{j}}=\frac12 f_{2\lambda,e_{i}\bar e_{j}}+\frac12
\varkappa\tilde C_{2,i,j},\,i\ne j,\\
[f_{\lambda,i},f_{\lambda,j}]&=\frac14 \varkappa
C_{1,e_{i},e_{j}}=-\frac12 f_{2\lambda,e_{i}\bar e_{j}}+\frac12
\varkappa\tilde C_{1,i,j},\,i\ne j,\\
[e_{\lambda,i},f_{\lambda,j}]&=\left\{-\frac12\Lambda,\,i=j \atop
-\frac12e_{2\lambda,e_{i}\bar e_{j}},\,i\neq j\right.,
\end{align*}
where we denote $f_{\lambda,e_{\alpha}e_{j}}:=f_{\lambda,i}$ if
$e_{\alpha}e_{j}=e_{i}$ and
$f_{\lambda,e_{\alpha}e_{j}}:=-f_{\lambda,i}$ if
$e_{\alpha}e_{j}=-e_{i}$. Analogous notation we use for
$e_{\lambda,i},e_{2\lambda,\gamma},f_{2\lambda,\gamma}$. Here
operators $\tilde C_{l,i,j},\,l=1,2,i\ne j$ are in $
\mathfrak{k}_{0}$ and act as:
\begin{align*}
\tilde C_{1,i,j}(e_{k})&=e_{k}e_{i}\cdot\bar e_{j},\,e_{k}\ne
1,\pm e_{i}\bar e_{j},\,\tilde C_{1,i,j}(e_{k})=0,\,e_{k}=1,\pm
e_{i}\bar e_{j},\\ \tilde C_{2,i,j}(e_{k})&=e_{j}\cdot\bar
e_{i}e_{k},\,e_{k}\ne 1,\pm e_{i}\bar e_{j},\,\tilde
C_{2,i,j}(e_{k})=0,\,e_{k}=1,\pm e_{i}\bar e_{j}.
\end{align*}
The chosen bases
$\Lambda,e_{\lambda,i},e_{2\lambda,\alpha},f_{\lambda,i},f_{2\lambda,\alpha}$
in spaces
$\mathfrak{a},\mathfrak{p}_{\lambda},\mathfrak{p}_{2\lambda},\mathfrak{k}_{\lambda},\mathfrak{k}_{2\lambda}$
correspond to proposition \ref{prop1}.
\end{proposit}
\begin{proof}
The commutator relations are easy consequences of
(\ref{OctCom1}),(\ref{OctCom2}), lemma \ref{lem1} and relations in
the algebra $\mathfrak{k}'\simeq \mathfrak{so}(8)$. For example,
let us calculate the commutator
$[f_{2\lambda,\alpha},e_{\lambda,j}]$. Actually, from
(\ref{OctCom1}) and lemma \ref{lem1} we obtain:
\begin{align*}
[f_{2\lambda,\alpha},e_{\lambda,j}]&=-\left[\varkappa
A_{\alpha},\frac12\ad Y_{2}(\bar e_{j})\right]=-\frac12\ad
Y_{2}\left(A^{(2)}_{\alpha}\bar e_{j}\right)=-\frac12\ad
Y_{2}\left(R_{\alpha}\bar e_{j}\right)\\&=-\frac14\ad
Y_{2}\left(\bar e_{j}e_{\alpha}\right)=\frac14\ad Y_{2}\left(\bar
e_{j}\bar e_{\alpha}\right)=-\frac12e_{\lambda,e_{\alpha}e_{j}}.
\end{align*}
Similar calculations are also valid for
$[f_{2\lambda,\alpha},f_{\lambda,j}]$.

Now, let us calculate $[e_{\lambda,i},e_{\lambda,j}],\,i\ne j$.
From (\ref{OctCom2}) we obtain:
$$
[e_{\lambda,i},e_{\lambda,j}]=\frac14\left[\ad Y_{2}(\bar
e_{i}),\ad Y_{2}(\bar e_{j})\right]=\frac14\varkappa C_{2,\bar
e_{i},\bar e_{j}},\,i\ne j.
$$
From (\ref{important}) and (\ref{C3}) we obtain
$$
\frac12C_{2,\bar e_{i},\bar
e_{j}}(e_{k})=\frac12\left(e_{j}\cdot\bar
e_{i}e_{k}-e_{i}\cdot\bar e_{j}e_{k}\right)=-e_{i}\cdot\bar
e_{j}e_{k}.
$$
In particular,
$$
\frac12C_{2,\bar e_{i},\bar e_{j}}(1)=-e_{i}\bar
e_{j},\,\frac12C_{2,\bar e_{i},\bar e_{j}}(e_{i}\bar
e_{j})=-\left(e_{i}\cdot\bar e_{j}\right)^{2}=1,
$$
so
$$
\frac12\varkappa C_{2,\bar e_{i},\bar e_{j}}=\varkappa
A_{e_{i}\bar e_{j}}+\varkappa \tilde
C_{2,i,j}=f_{2\lambda,e_{i}\bar e_{j}}+\varkappa \tilde C_{2,i,j},
$$
where $\tilde C_{2,i,j}\in\mathfrak{k}_{0}$ and
$$
\tilde C_{2,i,j}(e_{k})=e_{j}\cdot\bar e_{i}e_{k},\,e_{k}\ne 1,\pm
e_{i}\bar e_{j},\,\tilde C_{2,i,j}(e_{i}\bar e_{j})=\tilde
C_{2,i,j}(1)=0.
$$
The similar calculations are also valid for
$[f_{\lambda,i},f_{\lambda,j}]$.
\end{proof}

\subsection{Invariants in $S(\mathfrak{a}\oplus\mathfrak{p}_{\lambda}\oplus\mathfrak{k}_{\lambda}
\oplus\mathfrak{p}_{2\lambda}\oplus\mathfrak{k}_{2\lambda})$}
Invariant operators $D_{0},\dots,D_{6}$, corresponding to some
$K_{0}$-invariants in
$S(\mathfrak{a}\oplus\mathfrak{p}_{\lambda}\oplus\mathfrak{k}_{\lambda}
\oplus\mathfrak{p}_{2\lambda}\oplus\mathfrak{k}_{2\lambda})$ are
already constructed in section \ref{spes}. Here we shall construct
other independent invariants of $K_{0}$-action on
$S(\mathfrak{a}\oplus\mathfrak{p}_{\lambda}\oplus\mathfrak{k}_{\lambda}
\oplus\mathfrak{p}_{2\lambda}\oplus\mathfrak{k}_{2\lambda})$ or
equivalently from
$S(\mathfrak{p}_{\lambda}\oplus\mathfrak{k}_{\lambda}
\oplus\mathfrak{p}_{2\lambda}\oplus\mathfrak{k}_{2\lambda})$,
since $\mathfrak{a}$ is an invariant one dimensional space, and
corresponding invariant differential operators.

An element $\Phi\in K'$ is from $K_{0}\subset K'$ iff
$\Phi_{3}(1)=1$ and then $\Phi_{3}(\xi)=\xi$ for any
$\xi\in\mathbb{R}\subset\mathbb{C}a_{3}$. Below in this sections
$\Phi\in K_{0}$. The orthogonality of $\Phi_{i}$ means that
\begin{equation}
\label{orthogon}
\RE\left(\Phi_{i}(\xi)\overline{\Phi_{i}(\eta)}\right)=\RE(\xi\bar\eta),\,\xi,\eta\in\mathbb{C}a_{i}.
\end{equation}
In particular, $\Phi_{i}(\xi)\overline{\Phi_{i}(\xi)}=|\xi|^{2}$
and
\begin{equation}
\label{phinv}
\overline{\Phi_{i}(\xi)^{-1}}=\Phi_{i}(\xi)/|\xi|^{2}.
\end{equation}
For $\eta=\bar\xi$ from (\ref{fundument}) we obtain
$\Phi_{1}(\xi)\Phi_{2}(\bar\xi)=\overline{\Phi_{3}(|\xi|^{2})}=|\xi|^{2}$,
so from (\ref{phinv})
$\Phi_{1}(\xi)=|\xi|^{2}\Phi_{2}(\bar\xi)^{-1}=\overline{\Phi_{2}(\bar\xi)}$
and
\begin{equation}
\label{phi12} \Phi_{1}=\iota\circ\Phi_{2}\circ\iota.
\end{equation}
Let
$Q_{1}(\xi,\eta)=\RE(\xi\eta),\,\xi\in\mathbb{C}a_{1},\eta\in\mathbb{C}a_{2}$.
From (\ref{orthogon}) and (\ref{phi12}) we get:
$$
Q_{1}(\Phi_{1}(\xi),\Phi_{2}(\eta))=\RE(\Phi_{1}(\xi)\Phi_{2}(\eta))=\RE(\Phi_{1}(\xi)\overline{\Phi_{1}(\bar\eta)})=
\RE(\xi\overline{\bar\eta})=Q_{1}(\xi,\eta).
$$
Thus, $Q_{1}(\xi,\eta)$ is invariant under the $K_{0}$-action.

From proposition \ref{phi} it follows that
$\Phi_{1}=g^{L},\Phi_{2}=g^{R},\Phi_{3}=\iota\circ
g^{V}\circ\iota=g^{V}$, where $g^{L},g^{R},g^{V}$ are respectively
left spinor, right spinor and vector representation of the group
$K_{0}\simeq\Spin(7)$, since $\iota\,|_{\mathbb{C}a_{3}'}=-\id$.
Besides, the $K_{0}$-action on
$\im(\xi\eta),\,\xi\in\mathbb{C}a_{1},\eta\in\mathbb{C}a_{2}$
equals $g^{V}$, so
$Q_{2}(\xi,\eta,\zeta):=\RE\left(\im(\xi\eta)\zeta\right)$ is
invariant under $K_{0}$-action for $\zeta\in\mathbb{C}a_{3}'$.

According to section \ref{spes} the $K_{0}$-action on
$\mathfrak{p}_{\lambda}$ is equivalent to the $K_{0}$-action on
$\mathfrak{k}_{\lambda}$ and the $K_{0}$-action on
$\mathfrak{p}_{2\lambda}$ is equivalent to the $K_{0}$-action on
$\mathfrak{k}_{2\lambda}$. These equivalence is establishes by the
correspondence of bases $e_{\lambda,i}\leftrightarrow
f_{\lambda,i}$ and $e_{2\lambda,\alpha}\leftrightarrow
f_{2\lambda,\alpha}$. It is also confirmed by the formulas
$[Y_{2}(\xi),E_{1}]=X_{2}(\xi),\,\ad
Y_{1}(\xi)\left(X_{3}(1)\right)=-X_{2}(\bar\xi),\,\xi\in\mathbb{C}a$
and (\ref{phi12}). Therefore the analog of
$\im(\xi\eta),\,\xi\in\mathbb{C}a_{1},\eta\in\mathbb{C}a_{2}$ in
$S(\mathfrak{p}_{\lambda}\oplus\mathfrak{k}_{\lambda}
\oplus\mathfrak{p}_{2\lambda}\oplus\mathfrak{k}_{2\lambda})\otimes\mathbb{C}a$
is
$$
\sum_{i\ne j}f_{\lambda,i}e_{\lambda,\bar e_{j}}\otimes
e_{i}e_{j}.
$$

Thus, after the identification
$\mathfrak{p}_{2\lambda}\simeq\mathbb{C}a_{3}'$ the invariant
$Q_{2}$ gives the invariant from
$S(\mathfrak{p}_{\lambda}\oplus\mathfrak{k}_{\lambda}
\oplus\mathfrak{p}_{2\lambda}\oplus\mathfrak{k}_{2\lambda})$
$$
\sum_{i\ne j}f_{\lambda,i}e_{\lambda,\bar
e_{j}}e_{2\lambda,e_{i}e_{j}}=\sum_{i\ne
j}f_{\lambda,i}e_{\lambda,j}e_{2\lambda,e_{i}\bar e_{j}}.
$$
Therefore we can define the invariant differential operator:
\begin{equation*}
D_{7}=-\frac14\sum_{i\ne
j}\left\{\left\{f_{\lambda,i},e_{\lambda,j}\right\},e_{2\lambda,e_{i}\bar
e_{j}}\right\}=\frac14\sum_{i\ne
j}\left\{\left\{f_{\lambda,j},e_{\lambda,i}\right\},e_{2\lambda,e_{i}\bar
e_{j}}\right\}.
\end{equation*}
Similarly, the identification
$\mathfrak{k}_{2\lambda}\simeq\mathbb{C}a_{3}'$  gives the
invariant differential operator:
\begin{equation*}
D_{8}=-\frac14\sum_{i\ne
j}\left\{\left\{f_{\lambda,i},e_{\lambda,j}\right\},f_{2\lambda,e_{i}\bar
e_{j}}\right\}=\frac14\sum_{i\ne
j}\left\{\left\{f_{\lambda,j},e_{\lambda,i}\right\},f_{2\lambda,e_{i}\bar
e_{j}}\right\}.
\end{equation*}

It is clear that the equation (\ref{fundument}) remains valid
after the cyclic permutation of indices $1,2,3$:
\begin{align}
\label{fund1}
\Phi_{3}(\zeta)\Phi_{1}(\xi)=\overline{\Phi_{2}(\overline{\zeta\xi})},\,
\Phi_{2}(\eta)\Phi_{3}(\zeta)=\overline{\Phi_{1}(\overline{\eta\zeta})},\,
\xi\in\mathbb{C}a_{1},\eta\in\mathbb{C}a_{2},\zeta\in\mathbb{C}a_{3}'.
\end{align}
Define $$P(\xi,\eta,\zeta_{1},\zeta_{2}):=\RE(\zeta_{1}\xi\cdot
\eta\zeta_{2}),\,\zeta_{1},\zeta_{2}\in\mathbb{C}a_{3}'.$$ The
function $P(\xi,\eta,\zeta_{1},\zeta_{2})$ is invariant w.r.t. the
$K_{0}$-action, since due to (\ref{fund1}), (\ref{phi12}) and
(\ref{orthogon}):
\begin{align*}
P&\left(\Phi_{1}(\xi),\Phi_{2}(\eta),\Phi_{3}(\zeta_{1}),\Phi_{3}(\zeta_{2})\right)=
\RE\left(\Phi_{3}(\zeta_{1})\Phi_{1}(\xi)\cdot\Phi_{2}(\eta)\Phi_{3}(\zeta_{2})\right)\\
&=
\RE\left(\overline{\Phi_{2}\left(\overline{\zeta_{1}\xi}\right)}\;
\overline{\Phi_{1}\left(\overline{\eta\zeta_{2}}\right)}\right)=\RE\left(\Phi_{1}\left(\zeta_{1}\xi\right)
\overline{\Phi_{1}\left(\overline{\eta\zeta_{2}}\right)}\right)
=\RE\left(\zeta_{1}\xi\cdot\overline{\overline{\eta\zeta_{2}}}\right)=P(\xi,\eta,\zeta_{1},\zeta_{2}).
\end{align*}

Functions $P(\xi,\eta,\zeta_{1},\zeta_{2})$ and
$P(\xi,\eta,\zeta_{2},\zeta_{1})$ are not independent. Indeed, the
corollary 15.12 in \cite{Adams} gives:
$$
\RE(ab\cdot c)=\RE(bc\cdot a)=\RE(ca\cdot b)=\RE(a\cdot
bc)=\RE(b\cdot ca)=\RE(c\cdot ab),\,a,b,c\in\mathbb{C}a.
$$
Therefore using the Moufang identity (\ref{Moufang}) we obtain:
\begin{align}
\label{dependence} P(\xi,\eta,\zeta,\zeta)=\RE(\zeta\cdot
\xi\eta\cdot\zeta)=\RE(\zeta^{2}\cdot\xi\eta)=-\RE(|\zeta|^{2}\xi\eta)=-|\zeta|^{2}\RE(\xi\eta)=
-|\zeta|^{2}Q_{1}(\xi,\eta),
\end{align}
which means that for $\zeta_{1}=\zeta_{2}=\zeta$ invariant
$P(\xi,\eta,\zeta,\zeta)$ is expressed through invariants of the
second order. Using the polarization of (\ref{dependence}) w.r.t.
$\zeta$, which means the substitution $\zeta=\zeta_{1}+\zeta_{2}$,
we get:
$$
P(\xi,\eta,\zeta_{1},\zeta_{2})+P(\xi,\eta,\zeta_{2},\zeta_{1})=-2\langle\zeta_{1},\zeta_{2}\rangle
Q_{1}(\xi,\eta).
$$
It means the dependence of two invariants
$P(\xi,\eta,\zeta_{1},\zeta_{2})$,
$P(\xi,\eta,\zeta_{2},\zeta_{1})$ and invariants
$Q_{1}(\xi,\eta),\,\langle\zeta_{1},\zeta_{2}\rangle$ of the
second order. The last two invariants correspond to operators
$D_{3}$ and $D_{6}$.

For constructing the invariant differential operator $D_{9}$ we
shall use the invariant function
$$
P(\xi,\eta,\zeta_{1},\zeta_{2})-P(\xi,\eta,\zeta_{2},\zeta_{1}).
$$
Using $\sum_{k}e_{\lambda,\bar e_{k}}\otimes e_{k}$ as the analog
of $\eta$ we get the corresponding expression from
$S(\mathfrak{p}_{\lambda}\oplus\mathfrak{k}_{\lambda}
\oplus\mathfrak{p}_{2\lambda}\oplus\mathfrak{k}_{2\lambda})$:
\begin{align*}
&\phantom{=}\sum_{i\ne j\atop j\ne k}\left(f_{2\lambda,e_{j}\bar
e_{i}}f_{\lambda,i}e_{\lambda,\bar e_{k}}e_{2\lambda,\bar
e_{k}\bar e_{j}}-e_{2\lambda, e_{j}\bar
e_{i}}f_{\lambda,i}e_{\lambda,\bar e_{k}}f_{2\lambda,\bar
e_{k}\bar e_{j}}\right)\\&=\sum_{i\ne j\atop j\ne
k}\left(e_{2\lambda,e_{i}\bar
e_{j}}f_{\lambda,i}e_{\lambda,k}f_{2\lambda,e_{k}\bar
e_{j}}-f_{2\lambda,e_{i}\bar
e_{j}}f_{\lambda,i}e_{\lambda,k}e_{2\lambda,e_{k}\bar
e_{j}}\right),
\end{align*}
since for $i\ne j$ it holds $e_{j}\bar e_{i}=-\overline{e_{j}\bar
e_{i}}=-e_{i}\bar e_{j}$.

Define the corresponding invariant differential operator as
\begin{equation*}
D_{9}=\frac18\sum_{i\ne j\atop j\ne
k}\left(\left\{\left\{e_{2\lambda,e_{i}\bar
e_{j}},f_{\lambda,i}\right\},\left\{f_{2\lambda,e_{k}\bar
e_{j}},e_{\lambda,k}\right\}\right\}-\left\{\left\{e_{2\lambda,e_{i}\bar
e_{j}},e_{\lambda,i}\right\},\left\{f_{2\lambda,e_{k}\bar
e_{j}},f_{\lambda,k}\right\}\right\} \right).
\end{equation*}

Let us show that there are exactly $9$ independent
$K_{0}$-invariants in
$S(\mathfrak{p}_{\lambda}\oplus\mathfrak{k}_{\lambda}
\oplus\mathfrak{p}_{2\lambda}\oplus\mathfrak{k}_{2\lambda})$.

Indeed, $\dim(\mathfrak{p}_{\lambda}\oplus\mathfrak{k}_{\lambda}
\oplus\mathfrak{p}_{2\lambda}\oplus\mathfrak{k}_{2\lambda})=8+8+7+7=30$
and $\dim K_{0}=\dim\Spin(7)=21$. Therefore the codimension of
$K_{0}$-orbits in
$\mathfrak{p}_{\lambda}\oplus\mathfrak{k}_{\lambda}
\oplus\mathfrak{p}_{2\lambda}\oplus\mathfrak{k}_{2\lambda}$ is at
least $30-21=9$ and there should be at least $9$ independent
$K_{0}$-invariants in
$S(\mathfrak{p}_{\lambda}\oplus\mathfrak{k}_{\lambda}
\oplus\mathfrak{p}_{2\lambda}\oplus\mathfrak{k}_{2\lambda})$.

From another hand, it is obvious that the stationary subgroup,
corresponding to a point in a general position, of the group
$\Spin(7)$, acting in
$\mathfrak{p}_{2\lambda}\oplus\mathfrak{k}_{2\lambda}$ by
$g^{V}\oplus g^{V}$, is $\Spin(5)$. Therefore the dimension of
general $\Spin(7)$-orbits in
$\mathfrak{p}_{2\lambda}\oplus\mathfrak{k}_{2\lambda}$ is
$\dim\Spin(7)-\dim\Spin(5)=11$. The group $\Spin(5)$ is isomorphic
to $U_{\mathbb{H}}(2)$, see \cite{Adams}, Proposition 5.1. In the
section \ref{Quat} the six independent invariants of the diagonal
$U_{\mathbb{H}}(2)$-action in
$\mathbb{H}^{2}\oplus\mathbb{H}^{2}\simeq
\mathfrak{p}_{\lambda}\oplus\mathfrak{k}_{\lambda}$ were found, so
general orbits of the last action are $10$-dimensional, since
$\dim_{\mathbb{R}}(\mathbb{H}^{2}\oplus\mathbb{H}^{2})-6=10$.
Thus, general $\Spin(7)$-orbits in
$\mathfrak{p}_{\lambda}\oplus\mathfrak{k}_{\lambda}
\oplus\mathfrak{p}_{2\lambda}\oplus\mathfrak{k}_{2\lambda}$ are
$11+10=21$-dimensional, their codimension is $9$ and there are
exactly $9$ functionally independent invariants of
$\Spin(7)$-action in
$\mathfrak{p}_{\lambda}\oplus\mathfrak{k}_{\lambda}
\oplus\mathfrak{p}_{2\lambda}\oplus\mathfrak{k}_{2\lambda}$.

It is not known if there are any other invariants of this action,
which are polynomial in
$e_{\lambda,i},f_{\lambda,i},e_{2\lambda,\alpha},f_{2\lambda,\alpha}$
and are not polynomial in $D_{1},\dots,D_{9}$. Such invariants
should be connected with $D_{1},\dots,D_{9}$ by algebraic equation
of a degree more than one. In the case of
$\mathbf{P}^{n}(\mathbb{H})_{\mathbb{S}},\,n\geqslant 3$ there is
such invariant $D_{10}$ and $D_{10}^{2}$ is polynomial in
$D_{1},\dots,D_{9}$. The operator $D_{10}$ arises in commutative
relations of $D_{1},\dots,D_{9}$.

In the next section it is found that all commutators of operators
$D_{1},\dots,D_{9}$ in the octonionic case are polynomial in
$D_{1},\dots,D_{9}$. Therefore it seems probable that there is no
an analog of $D_{10}$ in the octonionic case.

It is easily verified that automorphisms $\zeta_{\alpha},\sigma$
acts on $D_{7},D_{8},D_{9}$ as
\begin{align*}
\zeta_{\alpha}(D_{7})&=\cos(\alpha) D_{7}-\sin(\alpha)
D_{8},\,\zeta_{\alpha}(D_{8})=\sin(\alpha)
D_{7}+\cos(\alpha)D_{8},\,\zeta_{\alpha}(D_{9})=D_{9},\\
\sigma(D_{7})&=D_{7}, \,\sigma(D_{8})=-D_{8}, \sigma(D_{9})=D_{9}.
\end{align*}

Similarly to the previous sections, in order to get the generators
of the algebra
\newline $\Diff(\mathbf{H}^{2}(\mathbb{C}a)_{\mathbb{S}})$ one can use
Remark \ref{transform} and make the formal substitution:
\begin{align*}
\Lambda &\rightarrow \ii\Lambda,e_{\lambda,i}\rightarrow\ii
e_{\lambda,i},f_{\lambda,i}\rightarrow
f_{\lambda,i},e_{2\lambda,\alpha}\rightarrow\ii
e_{\lambda,\alpha},f_{\lambda,\alpha}\rightarrow
f_{\lambda,\alpha}.
\end{align*}
This substitution produces the following substitution for the
generators $D_{0},\ldots,D_{10}$:
\begin{align}\label{GToct}
D_{0}&\rightarrow\ii \bar D_{0},\,D_{1}\rightarrow -\bar
D_{1},\,D_{2}\rightarrow \bar D_{2},\,D_{3}\rightarrow\ii \bar
D_{3},\,D_{4}\rightarrow -\bar D_{4},\,D_{5}\rightarrow \bar
D_{5},\nonumber\\ D_{6}&\rightarrow \ii\bar
D_{6},\,D_{7}\rightarrow -\bar D_{7},\,D_{8}\rightarrow \ii\bar
D_{8},\,D_{9}\rightarrow
 -\bar D_{9}.
\end{align}
The operators $\bar D_{0},\ldots,\bar D_{9}$ generate the algebra
$\Diff(\mathbf{H}^{2}(\mathbb{C}a)_{\mathbb{S}})$.

\section{Relations in algebras
$\Diff\left(\mathbf{P}^{2}(\mathbb{C}a)_{\mathbb{S}}\right)$ and
$\Diff\left(\mathbf{H}^{2}(\mathbb{C}a)_{\mathbb{S}}\right)$}
\label{RelCa}\markright{\ref{RelCa} Relations in algebras
$\Diff\left(\mathbf{P}^{2}(\mathbb{C}a)_{\mathbb{S}}\right)$ and
$\Diff\left(\mathbf{H}^{2}(\mathbb{C}a)_{\mathbb{S}}\right)$}

Below there are all $45$ commutative relations of operators
$D_{0},\dots,D_{9}$. An example of a calculation of such relation
is in appendix A. All methods described in section \ref{Relations}
for calculating commutative relations were used in this case.
Besides, the numeration the base elements
$e_{\lambda,i},f_{\lambda,i},e_{2\lambda,\alpha},f_{\lambda,\alpha}$
by octonionic units $e_{i},\,i=0,\dots,7$ is very convenient.

\begin{align}\label{octonions}
[D_{0},D_{1}]&=-D_{3},\,[D_{0},D_{2}]=D_{3},\,[D_{0},D_{3}]=\frac12
(D_{1}-D_{2}),\,[D_{0},D_{4}]=-2D_{6},\,\nonumber\\
[D_{0},D_{5}]&=2D_{6},\,[D_{0},D_{6}]=D_{4}-D_{5},
[D_{0},D_{7}]=-D_{8}, [D_{0},D_{8}]=D_{7},\,[D_{0},D_{9}]=0,\,
\nonumber\\ [D_{1},D_{2}]&=-\{D_{0},D_{3}\}-2D_{7},\,
[D_{1},D_{3}]=-\frac12\{D_{0},D_{1}\}+D_{8}+10D_{0},\,[D_{1},D_{4}]=2D_{7},\nonumber\\
[D_{1},D_{5}]&=0,\,[D_{1},D_{6}]=D_{8},\,
[D_{1},D_{7}]=\frac12\{D_{1},D_{2}-D_{4}\}-D_{9}-\frac12\{D_{3},D_{6}\}-D_{3}^{2}-5D_{0}^{2}\nonumber\\
&-\frac3{32}D_{1}-\frac{283}{32}D_{2}+\frac{19}2D_{4}-\frac12D_{5},\,
[D_{1},D_{8}]=-\frac12\{D_{3},D_{5}\}-\frac12\{D_{1},D_{6}\}+10D_{6}\nonumber\\
&+\frac{35}4D_{3},\,
[D_{1},D_{9}]=\frac12\{D_{5},D_{7}\}-\frac12\{D_{6},D_{8}\}-\frac{189}{32}\{D_{0},D_{3}\}-\frac{169}{16}D_{7},
\nonumber\\ [D_{2},D_{3}]&=\frac12\{D_{0},D_{2}\}+D_{8}-10D_{0},\,
[D_{2},D_{4}]=-2D_{7},\,[D_{2},D_{5}]=0,[D_{2},D_{6}]=-D_{8},\,\nonumber\\
[D_{2},D_{7}]&=-\frac12\{D_{2},D_{1}-D_{4}\}+D_{9}-\frac12\{D_{3},D_{6}\}+D_{3}^{2}+5D_{0}^{2}
+\frac3{32}D_{2}+\frac{283}{32}D_{1}-\frac{19}2D_{4}\nonumber\\
&+\frac12D_{5},\,[D_{2},D_{8}]=\frac12\{D_{2},D_{6}\}-\frac12\{D_{3},D_{5}\}+\frac{35}4
D_{3}-10D_{6},\nonumber\\
[D_{2},D_{9}]&=-\frac12\{D_{5},D_{7}\}+\frac12\{D_{6},D_{8}\}+\frac{189}{32}\{D_{0},D_{3}\}+\frac{169}{16}D_{7},
[D_{3},D_{4}]=0,\nonumber\\ [D_{3},D_{5}]&=2D_{8},\,
[D_{3},D_{6}]=D_{7},\,[D_{3},D_{7}]=-\frac14\{D_{1}+D_{2},D_{6}\}
+10D_{6},\nonumber\\
[D_{3},D_{8}]&=\frac12\{D_{1},D_{2}\}-\frac14\{D_{1}+D_{2},D_{5}\}
-D_{9}-D_{3}^{2}-5D_{0}^{2}-\frac{143}{32}(D_{1}+D_{2})-\frac12D_{4}\nonumber\\&+\frac{19}2D_{5},\,
[D_{3},D_{9}]=\frac12\{D_{4},D_{8}\}-\frac12\{D_{6},D_{7}\}+\frac{189}{64}\{D_{0},D_{1}-D_{2}\}
-\frac{169}{16}D_{8},\,\nonumber\\
[D_{4},D_{5}]&=-2\{D_{0},D_{6}\},\,[D_{4},D_{6}]=-\{D_{0},D_{4}\}+\frac{35}2D_{0},\,
[D_{4},D_{7}]=\frac12\{D_{1}-D_{2},D_{4}\}\nonumber\\&+\frac{35}4(D_{2}-D_{1}),\,
[D_{4},D_{8}]=\frac12\{D_{1}-D_{2},D_{6}\}-\{D_{0},D_{7}\},\nonumber\\
[D_{4},D_{9}]&=-9\{D_{0},D_{6}\},\,
[D_{5},D_{6}]=\{D_{0},D_{5}\}-\frac{35}2D_{0},\,
[D_{5},D_{7}]=\{D_{3},D_{6}\}+\{D_{0},D_{8}\},\nonumber\\
[D_{5},D_{8}]&=\{D_{3},D_{5}\}-\frac{35}2D_{3},\,
[D_{5},D_{9}]=9\{D_{0},D_{6}\},\,
[D_{6},D_{7}]=\frac14\{D_{1}-D_{2},D_{6}\}\nonumber\\
&+\frac12\{D_{3},D_{4}\}+\frac12\{D_{0},D_{7}\}-\frac{35}4D_{3},\,
[D_{6},D_{8}]=\frac14\{D_{1}-D_{2},D_{5}\}+\frac12\{D_{3},D_{6}\}\nonumber\\
&-\frac12\{D_{0},D_{8}\} +\frac{35}8(D_{2}-D_{1}),\,
[D_{6},D_{9}]=\frac92\{D_{0},D_{4}-D_{5}\},\nonumber\\
[D_{7},D_{8}]&=-\frac14\{D_{0},\{D_{1},D_{2}\}\}+\frac12\{D_{0},D_{3}^{2}\}+\frac12\{D_{0},D_{9}\}
+\frac14\{D_{1}-D_{2},D_{8}\}+\frac14\{D_{0},D_{5}\}\nonumber\\&+\frac{283}{64}\{D_{0},D_{1}+D_{2}\}-
\frac{175}{2}D_{0}-\frac12\{D_{3},D_{7}\}+5D_{0}^{3}+\frac14\{D_{0},D_{4}\},\,\nonumber\\
[D_{7},D_{9}]&=\frac14\{\{D_{0},D_{7}\},D_{6}\}+\frac18\{D_{2}-D_{1},\{D_{4},D_{5}\}\}-
\frac14\{\{D_{0},D_{4}\},D_{8}\}+\frac14\{D_{1}-D_{2},D_{6}^{2}\}\nonumber\\
&-\frac12\{D_{0},D_{8}\}+\frac{25}{32}\{D_{3},D_{6}\}+\frac{185}{64}\{D_{1}-D_{2},D_{4}\}
+\frac{17}{8}\{D_{1}-D_{2},D_{5}\}\nonumber\\&+\frac{35\cdot181}{128}(D_{2}-D_{1}),\,
[D_{8},D_{9}]=-\frac14\{\{D_{0},D_{6}\},D_{8}\}-\frac14\{D_{3},\{D_{4},D_{5}\}\}\nonumber\\&+
\frac14\{\{D_{0},D_{7}\},D_{5}\}+\frac12\{D_{3},D_{6}^{2}\}
+\frac{169}{32}\{D_{3},D_{5}\}+\frac{45}{64}\{D_{1}-D_{2},D_{6}\}\nonumber\\&+\frac{37}{8}\{D_{3},D_{4}\}
+\frac{5}{8}\{D_{0},D_{7}\}-\frac{35\cdot177}{64}D_{3}.\nonumber
\end{align}

Using these relations it is not difficult to verify that the
operator $D^{*}=D_{0}^{2}+D_{1}+D_{2}+D_{4}+D_{5}$ lies in the
centre of the algebra
$\Diff(\mathbf{P}^{n}(\mathbb{H})_{\mathbb{S}})$ in accordance
with the section \ref{spes}.

Using substitution (\ref{GToct}) one can obtains from relations
above the commutative relations for the algebra
$\Diff(\mathbf{H}^{2}(\mathbb{C}a)_{\mathbb{S}})$.

\section{Connection of algebras $\Diff(M_{\mathbb{S}})$ with the two-body problem}
\label{Connect}\markright{\ref{Connect} Connection with the
two-body problem}

In the paper \cite{Shchep2002} there was found an expression of
the quantum two-body Hamiltonian with a central potential
$V(\rho)$ on an arbitrary two-point homogeneous space $M$ through
radial differential operators and generators of an isometry group.
Using the notation of the present paper we can write these
expressions in the following way:
$$
\hat
H=L_{2}+\{L_{1},D_{0}\}+a_{0}D_{0}^{2}+\sum_{i=1}^{6}a_{i}D_{i}+V(\rho),
$$
for $M=\mathbf{P}^{n}(\mathbb{H})$ and
$M=\mathbf{P}^{2}(\mathbb{C}a)$;
$$
\hat
H=L_{2}+\{L_{1},D_{0}\}+a_{0}D_{0}^{2}+\sum_{i=1}^{3}a_{i}D_{i}+a_{4}D_{4}^{2}+a_{5}D_{5}^{2}
+a_{6}\{D_{4},D_{5}\}+V(\rho),
$$
for $M=\mathbf{P}^{n}(\mathbb{C})$;
$$
\hat
H=L_{2}+\{L_{1},D_{0}\}+a_{0}D_{0}^{2}+\sum_{i=1}^{3}a_{i}D_{i}+V(\rho),
$$
for $M=\mathbf{P}^{n}(\mathbb{R}),\,{\bf S}^{n},n\geqslant 3$ and
$$
\hat
H=L_{2}+\{L_{1},D_{0}\}+a_{0}D_{0}^{2}+a_{1}D_{1}^{2}+a_{2}D_{2}^{2}+\frac12a_{3}\{D_{1},D_{2}\}+V(\rho),
$$
for $M=\mathbf{P}^{2}(\mathbb{R}),\,{\bf S}^{2}$.

Here $\rho$ is the distance between particles, $L_{i},\,i=1,2$ is
some ordinary differential operator of the $i$-th order w.r.t.
$\rho$, $a_{0}=\text{const},\, a_{i},\,i=1,\dots,6$ are some
functions of $\rho$ and masses of particles. The analogous
expressions for noncompact spaces can be obtained by substitutions
$D_{i}\to \D_{i}$ from above.

The main difference of these expressions from Euclidean case is
the presence of noncommutative operators with coefficients
depending on $\rho$. This difference makes the two-body problem on
$M$ quite difficult. However, every common eigenfunction of
generators $D_{i}$ gives an isolated ordinary differential
equation for a radial part of an eigenfunction for $\hat H$. Using
this approach some exact spectral series for two-body problem on
$\bf S^{n}$ were found for several potentials in \cite{ShchStep}.
For other two-point compact homogeneous spaces similar
calculations should be more difficult.

\appendix
\section{Calculation of some commutative relations}
\label{Appendix A}\markright{\ref{Appendix A} Calculation of some
commutative relations}

In this appendix we shall illustrate the main ideas of calculating
some commutative relations. Let start from commutative relations
(\ref{CompComm}) from section \ref{Relations}. We shall obtain
some relations requiring the minimal calculations.

Let operators $D_{0},\ldots,D_{10}$ are defined as in section
\ref{Quat}. First let us consider the commutator $[D_{1},D_{4}]$.
It is not difficult to verify the following equalities for
elements $A,B,C$ of an arbitrary associative algebra:
\begin{equation}\label{A1}
 [A,\{B,C\}]=\{[A,B],C\}+\{B,[A,C]\},
\end{equation}
\begin{equation}\label{A2}
\{\{A,B\},C\}-\{A,\{B,C\}\}=[B,[A,C]],
\end{equation}
\begin{equation}\label{A3}
\{\{A,B\},C\}=2\{B,C\}A+\{[A,B],C\}+\{[A,C],B\}+[B,[A,C]].
\end{equation}
In particular, when $C=B$, from (\ref{A1}) we have
$[A,B^{2}]=\{[A,B],B\}$. This implies:
$$
[D_{1},D_{4}]=\{[D_{1},\Upsilon_{12}],\Upsilon_{12}\}+\{[D_{1},\Omega_{12}],\Omega_{12}\}+
\{[D_{1},\Theta_{12}],\Theta_{12}\}.
$$
Using (\ref{commutators1}) and (\ref{A1}) again, we obtain
$$
[D_{1},\Upsilon_{12}]=\square_{1},\,[D_{1},\Omega_{12}]=\square_{2},\,[D_{1},\Theta_{12}]=\square_{3}.
$$
Thus

\begin{equation}\label{A4}
[D_{1},D_{4}]=\{\square_{1},\Upsilon_{12}\}+\{\square_{2},\Omega_{12}\}+\{\square_{3},\Theta_{12}\}=2D_{7}.
\end{equation}

Using the permutation of coordinates $z_{1}$ and $z_{2}$ (or
equivalently the automorphism $\sigma\circ\zeta_{\pi}$, see
section \ref{Quat}), we obtain from (\ref{A4}):
$$
[D_{2},D_{4}]=-2D_{7}.
$$

Suppose now we already know the expressions for commutators
\begin{align*}
&[D_{0},D_{1}],\,[D_{0},D_{3}],\,[D_{0},D_{7}],\,[D_{1},D_{2}],\,[D_{1},D_{4}],\,[D_{1},D_{5}],\,[D_{1},D_{6}],\,
[D_{1},D_{7}],\\
&[D_{1},D_{8}],\,[D_{2},D_{6}],\,[D_{3},D_{4}],\,[D_{3},D_{6}],\,[D_{4},D_{5}],\,[D_{4},D_{6}],\,[D_{4},D_{8}].
\end{align*}
Then from the Jacobi identity and (\ref{A1}) we have
\begin{align*}
0&=[D_{1},[D_{8},D_{4}]]+[D_{4},[D_{1},D_{8}]]+[D_{8},[D_{4},D_{1}]]=[D_{1},\frac12\{D_{2}-D_{1},D_{6}\}+\{D_{0},D_{7}\}]
\\ &+[D_{4},n(n-1)D_{6}-\frac12\{D_{3},D_{5}\}+\frac34D_{3}-\frac12\{D_{1},D_{6}\}]-2[D_{8},D_{7}]\\&=
\frac12\{[D_{1},D_{2}],D_{6}\}-\frac12\{D_{1}-D_{2},[D_{1},D_{6}]\}+\{[D_{1},D_{0}],D_{7}\}+\{D_{0},[D_{1},D_{7}]\}
\\&+n(n-1)[D_{4},D_{6}]-\frac12\{[D_{4},D_{3}],D_{5}\}-\frac12\{D_{3},[D_{4},D_{5}]\}+\frac34[D_{4},D_{3}]-
\frac12\{D_{1},[D_{4},D_{6}]\}\\&-\frac12\{[D_{4},D_{1}],D_{6}\}-2[D_{8},D_{7}]=-\frac12\{\{D_{3},D_{0}\},D_{6}\}
-\{D_{4},D_{6}\}+\frac12\{D_{2}-D_{1},D_{8}\}\\&+\{D_{3},D_{7}\}+\{D_{0},n(n-1)D_{4}-\frac12\{D_{3},D_{6}\}-
\frac12\{D_{1},D_{4}\}+\frac38(D_{1}-D_{2})+D_{9}+D_{10}\}\\&-n(n-1)\{D_{0},D_{4}\}+\frac32n(n-1)D_{0}+
\{D_{3},\{D_{6},D_{0}\}\}+\frac12\{D_{1},\{D_{0},D_{4}\}\}-\frac34\{D_{1},D_{0}\}\\&+\{D_{7},D_{6}\}-2[D_{8},D_{7}]=
\frac12\{D_{2}-D_{1},D_{8}\}+\{D_{3},D_{7}\}+\frac38\{D_{0},D_{1}-D_{2}\}\\&+\{D_{0},D_{9}+D_{10}\}+\frac32n(n-1)D_{0}-
\frac34\{D_{0},D_{1}\}-2[D_{8},D_{7}].
\end{align*}
In the last equality we took into account the formulas
\begin{align*}
&\{\{D_{6},D_{0}\},D_{3}\}-\{D_{6},\{D_{0},D_{3}\}\}+\{\{D_{0},D_{6}\},D_{3}\}-\{D_{0},\{D_{6},D_{3}\}\}=
[D_{0},[D_{6},D_{3}]]\\&+[D_{6},[D_{0},D_{3}]]=-[D_{0},D_{7}]-\frac12[D_{6},D_{2}-D_{1}]=D_{8}+\frac12(-D_{8}-D_{8})=0,\\
&\{\{D_{0},D_{4}\},D_{1}\}-\{D_{0},\{D_{4},D_{1}\}\}=[D_{4},[D_{0},D_{1}]]=[D_{4},D_{3}]=0,
\end{align*}
which are consequences of (\ref{A2}).

Thus we get:
\begin{align*}
[D_{7},D_{8}]&=\frac14\{D_{1}-D_{2},D_{8}\}-\frac12\{D_{3},D_{7}\}+\frac3{16}\{D_{0},D_{1}+D_{2}\}-
\frac12\{D_{0},D_{9}+D_{10}\}\nonumber\\&-\frac34n(n-1)D_{0}.
\end{align*}

Now let us demonstrate the calculation modulo
$(U(\mathfrak{g})\mathfrak{k}_{0})^{K_{0}}$. Let
$D_{0},\ldots,D_{3}$ are generators of $\Diff({\bf
S}^{n}_{\mathbb{S}}),\,n\geqslant 3$,
$\mathfrak{g}=\mathfrak{so}(n+1),\,\mathfrak{k}_{0}=\mathfrak{so}(n-1),\,K_{0}=\SO(n-1)$.
Then from (\ref{A1}) we obtain:
\begin{align}\label{star}
[D_{1},D_{3}]&=-8\sum\limits_{k,l=3}^{n+1}\left(\{\{\Psi_{1k},[\Psi_{1k},\Psi_{1l}]\},\Psi_{2l}\}+
\{\Psi_{1l},\{\Psi_{1k},[\Psi_{1k},\Psi_{2l}]\}\}\right)\nonumber\\
&=4\sum\limits_{k,l=3}^{n+1}\left(\{\{\Psi_{1k},\Psi_{kl}\},\Psi_{2l}\}+\delta_{kl}
\{\Psi_{1l},\{\Psi_{1k},\Psi_{12}\}\}\right)\nonumber\\
&=4\sum\limits_{k,l=3\atop k\neq
l}^{n+1}\{\{\Psi_{kl},\Psi_{1k}\},\Psi_{2l}\}+4\sum\limits_{k=3}^{n+1}\{\Psi_{1k},\{\Psi_{1k},\Psi_{12}\}\}.
\end{align}
From formula (\ref{A3}) and commutative relations
(\ref{commutators1}) one obtains:
\begin{align*}
&\sum\limits_{k,l=3\atop k\neq
l}^{n+1}\{\{\Psi_{kl},\Psi_{1k}\},\Psi_{2l}\}=\sum\limits_{k,l=3\atop
k\neq
l}^{n+1}(2\{\Psi_{1k},\Psi_{2l}\}\Psi_{kl}+\{[\Psi_{kl},\Psi_{1k}],\Psi_{2l}\}+
\{[\Psi_{kl},\Psi_{2l}],\Psi_{1k}\}\\
&+[\Psi_{1k},[\Psi_{kl},\Psi_{2l}]])\equiv\sum\limits_{k,l=3\atop
k\neq
l}^{n+1}(-\frac12\{\Psi_{1l},\Psi_{2l}\}+\frac12\{\Psi_{2k},\Psi_{1k}\}+\frac12[\Psi_{1k},\Psi_{2k}])
\mod{\left(U(\mathfrak{g})}\mathfrak{k}_{0}\right)^{K_{0}}\\
&=-\frac14\sum\limits_{k,l=3\atop k\neq
l}^{n+1}\Psi_{12}=-\frac{(n-1)(n-2)}4\Psi_{12}=\frac{(n-1)(n-2)}8D_{0}.
\end{align*}
Formula (\ref{A2}) gives:
\begin{align*}
&\sum\limits_{k=3}^{n+1}\{\Psi_{1k},\{\Psi_{1k},\Psi_{12}\}\}=
\sum\limits_{k=3}^{n+1}(\{\{\Psi_{1k},\Psi_{1k}\},\Psi_{12}\}-[\Psi_{1k},[\Psi_{1k},\Psi_{12}]])=
2\left\{\sum\limits_{k=3}^{n+1}\Psi_{1k}^{2},\Psi_{12}\right\}\\
&-\frac12\sum\limits_{k=3}^{n+1}[\Psi_{1k},\Psi_{2k}]=-\frac14\{D_{0},D_{1}\}+
\frac14\sum\limits_{k=3}^{n+1}\Psi_{12}=-\frac14\{D_{0},D_{1}\}-\frac{n-1}8D_{0}.
\end{align*}
Finally, from (\ref{star}) we obtain:
\begin{align*}
[D_{1},D_{3}]=-\{D_{0},D_{1}\}-\frac{n-1}2D_{0}+\frac{(n-1)(n-2)}2D_{0}=
-\{D_{0},D_{1}\}+\frac{(n-1)(n-3)}2D_{0}.
\end{align*}

Calculation of $[D_{1},D_{3}]$ for algebras
$\Diff(\mathbf{P}^{n}(\mathbb{H})_{\mathbb{S}})$ and
$\Diff(\mathbf{P}^{n}(\mathbb{C})_{\mathbb{S}})$ are analogous,
but much longer.

Let us demonstrate calculations in octonionic case by one example.
Below indices $i,j$ vary from $0$ to $7$. Let $D_{0},\ldots,D_{9}$
are generators of
$\Diff\left(\mathbf{P}^{2}(\mathbb{C}a)_{\mathbb{S}}\right)$,
$\mathfrak{g}=\mathfrak{f}_{4},\,\mathfrak{k}_{0}=\mathfrak{spin}(7),\,K_{0}=\Spin(7)$.
Then from (\ref{A1}) and proposition \ref{octbase} we obtain:
\begin{align*}
&\phantom{=}[D_{1},D_{3}]=\frac12\sum_{i,j}\left(\left\{\{[e_{\lambda,i},e_{\lambda,j}],e_{\lambda,i}\},f_{\lambda,j}\right\}
+\left\{e_{\lambda,j},\{e_{\lambda,i},[e_{\lambda,i},,f_{\lambda,j}]\}\right\}\right)\\&=
\frac18\sum_{i\ne j}\left\{\{\varkappa C_{2,\bar e_{i},\bar
e_{j}},e_{\lambda,i}\},f_{\lambda,j}\right\}-\frac14\sum_{i}\left\{e_{\lambda,i},\{e_{\lambda,i},\Lambda\}\right\}
-\frac14\sum_{i\ne
j}\left\{e_{\lambda,j},\{e_{\lambda,i},e_{2\lambda,e_{i}\bar
e_{j}}\}\right\}.
\end{align*}
Formulas (\ref{A3}), (\ref{C12}), (\ref{important})  and
proposition \ref{octbase} imply:
\begin{align*}
&\phantom{=}\frac18\sum_{i\ne j}\left\{\{\varkappa C_{2,\bar
e_{i},\bar
e_{j}},e_{\lambda,i}\},f_{\lambda,j}\right\}=\frac18\sum_{i,j}\left(2\{e_{\lambda,i},f_{\lambda,j}\}
\varkappa C_{2,\bar e_{i},\bar e_{j}}+\left\{[\varkappa C_{2,\bar
e_{i},\bar
e_{j}},e_{\lambda,i}],f_{\lambda,j}\right\}\right.\\&+\left(\left\{[\varkappa
C_{2,\bar e_{i},\bar
e_{j}},f_{\lambda,j}],e_{\lambda,i}\right\}+\left[e_{\lambda,i},\left[\varkappa
C_{2,\bar e_{i},\bar e_{j}},f_{\lambda,j}\right]\right]\right)\\
&\equiv\sum_{i\ne j}\left(-\frac1{16}\left\{\ad
Y_{2}\left(\left.\varkappa C_{2,\bar e_{i},\bar
e_{j}}\right|_{\mathbb{C}a_{2}}\bar
e_{i}\right),f_{\lambda,j}\right\}+\frac1{16}\left\{\ad
Y_{1}\left(\left.\varkappa C_{2,\bar e_{i},\bar
e_{j}}\right|_{\mathbb{C}a_{1}}e_{j}\right),e_{\lambda,i}\right\}\right.\\
&+\left.\frac1{16}\left[e_{\lambda,i},\ad
Y_{1}\left(\left.\varkappa C_{2,\bar e_{i},\bar
e_{j}}\right|_{\mathbb{C}a_{1}}e_{j}\right)\right]+\frac12
\{e_{\lambda,i},f_{\lambda,j}\}f_{2\lambda,e_{i}\bar
e_{j}}\right)\mod (U(\mathfrak{g})\mathfrak{k}_{0})^{K_{0}}\\
&=\sum_{i\ne j}\left(-\frac14\left\{\ad Y_{2}\left(\bar
e_{j}\right),f_{\lambda,j}\right\}-\frac18\left\{\ad
Y_{1}\left(e_{i}\right),e_{\lambda,i}\right\}-\frac18\left[e_{\lambda,i},\ad
Y_{1}\left(e_{i}\right)\right]\right.\\&+\left.\frac14\left(\{f_{2\lambda,e_{i}\bar
e_{j}},\{e_{\lambda,i},f_{\lambda,j}\}\}-[f_{2\lambda,e_{i}\bar
e_{j}},\{e_{\lambda,i},f_{\lambda,j}\}]\right)\right)\\&=\sum_{i\ne
j}\left(\frac12\{e_{\lambda,j},f_{\lambda,j}\}-\frac14\{f_{\lambda,i},e_{\lambda,i}\}-
\frac14[e_{\lambda,i},f_{\lambda,i}]\right)+D_{8}\\&-\frac14\sum_{i\ne
j}\left(\left\{[f_{2\lambda,e_{i}\bar
e_{j}},e_{\lambda,i}],f_{\lambda,j}\right\}+\{e_{\lambda,i},[f_{2\lambda,e_{i}\bar
e_{j}},f_{\lambda,j}]\}\right)\\&=D_{8}+\sum_{i\ne
j}\left(\frac14\{e_{\lambda,j},f_{\lambda,j}\}+\frac18\Lambda+\frac18\{e_{\lambda,e_{i}\bar
e_{j}\cdot
e_{i}},f_{\lambda,j}\}-\frac18\{e_{\lambda,i},f_{\lambda,e_{i}\bar
e_{j}\cdot e_{j}}\}\right)\\&=D_{8}+\sum_{i\ne
j}\left(\frac14\{e_{\lambda,j},f_{\lambda,j}\}-\frac18\{e_{\lambda,e_{j}\bar
e_{i}e_{i}},f_{\lambda,j}\}-\frac18\{e_{\lambda,i},f_{\lambda,i}\}\right)+7\Lambda=D_{8}+7D_{0}.
\end{align*}
Similarly, from (\ref{A2}) and proposition \ref{octbase} we get:
\begin{align*}
&-\frac14\sum_{i}\{\{\Lambda,e_{\lambda,i}\},e_{\lambda,i}\}=
-\frac14\sum_{i}\left(\{\Lambda,\{e_{\lambda,i},e_{\lambda,i}\}\}+[e_{\lambda,i},[\Lambda,e_{\lambda,i}]]\right)\\&=
-\frac12\left\{\Lambda,\sum_{i}e_{\lambda,i}^{2}\right\}+\frac18\sum_{i}[e_{\lambda,i},f_{\lambda,i}]=
-\frac12\{D_{0},D_{1}\}-\frac12D_{0}.
\end{align*}
Also
\begin{align*}
&-\frac14\sum_{i\ne j}\left\{\{e_{2\lambda,e_{i}\bar
e_{j}},e_{\lambda,i}\},e_{\lambda,j}\right\}=-\frac14\sum_{i\ne
j}\left(\left\{e_{2\lambda,e_{i}\bar
e_{j}},\{e_{\lambda,i},e_{\lambda,j}\}\right\}+[e_{\lambda,i},[e_{2\lambda,e_{i}\bar
e_{j}},e_{\lambda,j}]]\right)\\&=-\frac18\sum_{i\ne
j}[e_{\lambda,i},f_{\lambda,e_{i}\bar e_{j}\cdot
e_{j}}]=-\frac18\sum_{i\ne
j}[e_{\lambda,i},f_{\lambda,i}]=\frac{7\cdot 8}{2\cdot
8}\Lambda=\frac72 D_{0},
\end{align*}
since $e_{2\lambda,e_{i}\bar e_{j}}$ is antisymmetric and
$\{e_{\lambda,i},e_{\lambda,j}\}$ is symmetric w.r.t. $i,j$.

Thus
$$
[D_{1},D_{3}]=D_{8}+(7-\frac12+\frac72)D_{0}-\frac12\{D_{0},D_{1}\}=D_{8}-\frac12\{D_{0},D_{1}\}+10D_{0}.
$$

\section{}
\label{Appendix B}\markright{\ref{Appendix B} }

In this appendix we will prove the following theorem:
\begin{theore}
Let $M$ be a two point $G$-homogeneous Riemannian space, where $G$
is the identity component of the isometry group for $M$. For every
smooth vector field $v$ on $M$ define a function $f_{v}$ on $M_{
\mathbb{S}}$ by the following formula:
$$
f_{v}(y)=\hat g (v(x),\xi)\equiv \langle v(x),\xi\rangle,
$$
where $x\in M,\,\hat
g(\cdot,\cdot)\equiv\langle\cdot,\cdot\rangle$ is the Riemannian
metric on $M$, $\xi\in
T_{x}M,\,\langle\xi,\xi\rangle=1,\,y=(x,\xi)\in M_{\mathbb{S}}$.
Let $D_{0}\in\Diff(M_{\mathbb{S}})$ be the differential operator
constructed in section \ref{spes} (for the noncompact case, see
Remark \ref{transform}). For every element $X\in\mathfrak{g}$ we
denote by $\tilde X$ the corresponding Killing vector field on
$M$. Then the condition $D_{0}f_{v}\equiv 0$ on $M_{\mathbb{S}}$
is equivalent to the equality $v=\tilde X$ for some $X\in
\mathfrak{g}$. In other words, the kernel of the operator $D_{0}$
consist of functions $f_{\tilde X}$, where $X$ runs over the
algebra $\mathfrak{g}$.
\end{theore}
This theorem for the case $M=\mathbf{H}^{n}(\mathbb{R})$ was
formulated and proved in \cite{Reimann} by the explicit coordinate
calculations. Here we will prove it in the general case in a more
conceptual way.
\begin{proof}
Let $K$ be the stationary subgroup corresponding to the point
$x_{0}\in M$, $e_{0}=\frac1R\tilde \Lambda(x_{0})\in
T_{x_{0}}M,\,\langle e_{0},e_{0}\rangle=1$, where $\Lambda$ and
$R$ are from Proposition 1. The space $M_{\mathbb{S}}$ is the
$G$-orbit $Gy_{0}$, where $y_{0}=(x_{0},e_{0})\in M_{\mathbb{S}}$.

The action of $D_{0}$ on $f_{v}$ can be written in the following
way \cite{Hel84}(theorem 4.3):
$$
(D_{0}f_{v})(gy_{0})=\left.\d1{t}\right|_{t=0}f_{v}\left(g\exp(t\Lambda)y_{0}\right),\,g\in
G.
$$
Therefore,
\begin{align*}
(D_{0}f_{v})(gy_{0})&=\left.\d1{t}\right|_{t=0}\left\langle
v\left(g\exp(t\Lambda)x_{0}\right),g\exp(t\Lambda)e_{0}\right\rangle\\
&=\left.\d1{t}\right|_{t=0}\left\langle
v\left(g\exp(t\Lambda)g^{-1}gx_{0}\right),g\exp(t\Lambda)\left.\d1{\mu}\right|_{\mu=0}
\exp(\mu\Lambda)x_{0}\right\rangle\\
&=\left.\d1{t}\right|_{t=0}\left\langle
v\left(\exp(t\Ad_{g}\Lambda)gx_{0}\right),\left.\d1{\mu}\right|_{\mu=0}\exp(t\Ad_{g}\Lambda)
\exp(\mu\Ad_{g}\Lambda)gx_{0}\right\rangle\\
&=\left.\d1{t}\right|_{t=0}\left\langle
v\left(\exp(t\Ad_{g}\Lambda)gx_{0}\right),\left.\d1{\mu}\right|_{\mu=0}
\exp(\mu\Ad_{g}\Lambda)\exp(t\Ad_{g}\Lambda)gx_{0}\right\rangle.
\end{align*}
Due to the transitivity of $G$-action on $M_{\mathbb{S}}$ the
point
$y:=(x,e):=\left(gx_{0},\left.\widetilde{\Ad_{g}\Lambda}\right|_{gx_{0}}\right)$
can be considered as arbitrary. Denote $W=\Ad_{g}\Lambda$. Then
\begin{align*}
(D_{0}f_{v})(y)=\left.\d1{t}\right|_{t=0}\left\langle
v\left(\exp(tW)x\right),\tilde W
\left(\exp(tW)x\right)\right\rangle=\pounds_{\tilde W}\hat
g(v(x),\tilde W(x)),
\end{align*}
where $\pounds_{X}$ is the Lie derivative along the vector field
$X$. The vector field $\tilde W$ is Killing, so $\pounds_{\tilde
W}\hat g=0$ and
\begin{align}\label{end}
D_{0}f_{v}&=\hat g(\pounds_{\tilde W}v,\tilde W)+\hat
g(v,\pounds_{\tilde W}\tilde W)=\hat g(\pounds_{\tilde W}v,\tilde
W)=-\hat g(\pounds_{v}\tilde W,\tilde W)\\
&=\frac12\left(\pounds_{v}\hat g\right)(\tilde W,\tilde
W)\nonumber-\frac12\pounds_{v}\left(\hat g(\tilde W,\tilde
W)\right)=\frac12\left(\pounds_{v}\hat g\right)(\tilde W,\tilde
W),
\end{align}
due to $\hat g(\tilde W(x),\tilde W(x))=\hat g(g\tilde
\Lambda(x_{0}),g\tilde \Lambda(x_{0}))=\hat g(\tilde
\Lambda(x_{0}),\tilde \Lambda(x_{0}))=R^{2}$ and
$\pounds_{X}Y=[X,Y]_{c}$ where $[X,Y]_{c}$ is the commutator of
vector fields $X$ and $Y$. The element $\tilde W\in T_{x}M$ is
arbitrary, therefore from (\ref{end}) we see that the condition
$D_{0}f_{v}=0$ is equivalent to the equality $\pounds_{v}\hat
g=0$, which means that $v$ is a Killing vector field and has the
form $v=\tilde X$ for some $X\in \mathfrak{g}$ if and only if
$D_{0}f_{v}\equiv 0$.
\end{proof}

\end{document}